\documentclass[12pt]{article}

\usepackage{amsmath,amssymb,amsthm,amscd,a4wide,upref}
\usepackage[toc]{appendix}

\usepackage{hyperref}
\hypersetup{colorlinks,citecolor=blue,filecolor=black,linkcolor=blue,urlcolor=blue}
\usepackage{enumerate}
\usepackage{mathtools}

\begin{document}
	\newcommand{\End}{\mathrm{End}}
	
	\newcommand{\diag}{\mathrm{diag}}
	\newcommand{\CC}{\mathbb{C}}
	\newcommand{\ZZ}{\mathbb{Z}}
	\newcommand{\X}{\mathrm{X}}
	\newcommand{\Y}{\mathrm{Y}}
	\newcommand{\ZX}{\mathrm{ZX}}
	
	\numberwithin{equation}{section}
	
	\newtheorem{thm}{Theorem}[section]
	\newtheorem{lem}[thm]{Lemma}
	\newtheorem{prop}[thm]{Proposition}
	\newtheorem{cor}[thm]{Corollary}
	\newtheorem{conj}[thm]{Conjecture}
	\newtheorem*{mthm}{Main Theorem}
	\newtheorem*{mthma}{Theorem A}
	\newtheorem*{mthmb}{Theorem B}
	\newtheorem*{mthmc}{Theorem C}
	\newtheorem*{mthmd}{Theorem D}

	\theoremstyle{definition}
	\newtheorem{defin}[thm]{Definition}
	
	\theoremstyle{remark}
	\newtheorem{remark}[thm]{Remark}
	\newtheorem{example}[thm]{Example}
	
	\numberwithin{equation}{section}
	
	\renewcommand{\labelenumi}{(\roman{enumi})}
	
	\allowdisplaybreaks[4]
	
	\title{\Large\bf Parabolic presentations of  Yangian in types $B$ and $C$}
	
	\author{Zhihua Chang, Naihuan Jing, Ming Liu\footnote{Corresponding Author: M. Liu (Email: ming.l1984@gmail.com)}, and Haitao Ma}
	
	\date{}
	\maketitle
	
	\begin{abstract}
		We establish a parabolic presentation of the extended Yangian $\X(\mathfrak{g}_{N})$ associated with the Lie algebras $\mathfrak{g}_{N}$ of type $B$ and $C$, parameterized by a symmetric composition $\nu$ of $N$. By formulating a block matrix version of the RTT presentation of $\X(\mathfrak{g}_{N})$, we systematically derive the generators and relations through the Gauss decomposition of the generator matrix in $\nu$-block form. Furthermore, leveraging this parabolic presentation, we obtain a novel formula for the center of $\X(\mathfrak{g}_{N})$, offering new insights into its structure.
		
		\medskip
		
		\medskip
		
		
	\end{abstract}
	
	\section{Introduction}
	\label{sec:int}
	
	Yangians are a family of quantum groups that originated from the study of quantum inverse scattering theory in the work of Faddeev and the St.Petersburg school \cite{ks:qs}. Drinfeld defined them \cite{d:ha} as Hopf algebras to deform the universal enveloping algebras $\mathrm{U}(\mathfrak{g}[x])$ of the current algebras associated with simple finite-dimensional Lie algebras $\mathfrak{g}$. Yangians play a significant role in diverse areas of research in mathematics and physics, including rational solutions of the Yang-Baxter equation in the theory of integrable models \cite{d:ha}, representations of classical Lie algebras \cite{m:yc} and finite W-algebras \cite{bk:walg1}.
	
	The Yangian $\Y(\mathfrak{g})$, which corresponds to a simple finite-dimensional Lie algebra $\mathfrak {g}$ of types $A$, $B$, $C$, or $D$, is generally described in three distinct forms: the $J$-presentation, featuring a finite number of generators as outlined in \cite{d:ha}; the Drinfeld presentation described \cite{d:nr}; and the RTT presentation, which includes a single ternary relation on the generator matrix and draws inspiration from the Faddeev-Reshetikhin-Takhtajan quantum group presentation \cite{frt:rtt}. The relationship between the Drinfeld and RTT presentations has been demonstrated in \cite{bk:pp} for type A and in \cite{jl:ib,jlm:ib} for types $B$, $C$, and $D$. Additionally, the relationship between the $J$-presentation and the Drinfeld presentation was addressed in 	\cite{grw:eb}. All three formats serve as important tools for examining Yangians \cite{m:yc}. Recently, Lu, Wang, Zhang \cite{LWZ1,LWZ2} have employed the method of Gauss decomposition to investigate the Drinfeld presentation of twisted Yangians.
	
	A parabolic presentation of the Yangian $\Y(\mathfrak{gl}_n)$ associated with a composition $\lambda$ of $n$ was initiated in \cite{bk:pp} for type $A$. This presentation
	utilized a block Gauss decomposition of the generator matrix $T(u)$ in the RTT presentation. The parabolic presentation of $\Y(\mathfrak{gl}_n)$ reduces to the RTT presentation if $\lambda=(n)$, and to the Drinfeld presentation if $\lambda=(1,1,\ldots,1)$. With the help of a parabolic presentation, the parabolic subalgebras of $\Y(\mathfrak{gl}_n)$ can be explicitly described. These
	parabolic presentations have played a crucial role in subsequent work \cite{bk:walg1}, where
	it was demonstrated that a finite W-algebra of type A is isomorphic to some truncated shifted Yangian.
	This work has led to a further development of the highest weight module theory for shifted Yangians and finite W-algebras \cite{bk:walg2}. As super analogues of these works, Peng obtained a parabolic presentation of the super Yangian $\Y(\mathfrak{gl}_{m|n})$ \cite{peng11,peng16}, introduced a shifted super Yangian as a subalgebra of $Y(\mathfrak{gl}_{m|n})$ and established an isomorphism from the finite W-superalgebra to a quotient of a shifted super Yangian  \cite{peng14,peng15,peng21}.
	
	Recently, Kamnitzer, Webster, Weekes and Yacobi introduced the general shifted Yangian for an arbitrary simple Lie algebra, which is the quantization of transverse slices to Schubert varieties in the affine Grassmannian \cite{kwwy:yq} and generalized Brundan-Kleshchev's shifted Yangian \cite{bk:walg1}. In \cite{bfn:cg}, the authors established the isomorphism between the quantized Coulomb branch of a framed quiver gauge theory and truncated shifted Yangian for simply laced types.
	More recently, Nakajima and Weekes \cite{nw:cs} proved the isomorphism between the quantized Coulomb branch of a framed quiver gauge theory and truncated shifted Yangians for non-simply-laced types. Also Wang \cite{W} recently announced their forthcoming works on finite W-algebras, shifted twisted Yangians and affine Grassmannian slices.
	
	The close and deep relation between Yangians and finite W-algebras in type $A$ motivates us to explore a similar/alternative theory in type $B$, $C$, and $D$. The primary aim of this paper is to provide a parabolic presentation for the Yangian of type $B$ and $C$.
	Their connection with finite $W$-algebras will be considered in our subsequent works.
	
	We work with the extended Yangian $\X(\mathfrak{g}_{N})$ for $\mathfrak{g}_N=\mathfrak{o}_{2n+1}$ or $\mathfrak{g}_N=\mathfrak{sp}_{2n}$ since it is presented by a single ternary relation. The Yangian $\Y(\mathfrak{g}_N)$ is the quotient of $\X(\mathfrak{g}_N)$ by central elements as shown in \cite{amr:rp}. We follow the general strategy in \cite{bk:pp} to obtain the generators in a parabolic presentation, yet through a new block Gauss decomposition of the generator matrix $T(u)$ of $\X(\mathfrak{g}_N)$ to simplify the computation. The block Gauss decomposition relies on a symmetric composition $\nu=(\nu_1,\nu_2,\ldots,\nu_M)$ of $N$, where a symmetric composition means $\nu_a=\nu_{M+1-a}$ for $a=1,\ldots,M$.
	
	To effectively determine the relations among these generators, we have developed a new block matrix method to
	show the embedding theorem obtained in \cite{jlm:ib}. The new block matrix method has simplified many key computations and provides
	a new perspective to present the relations for the Gauss generators. The problem is then reduced to the cases where $\nu$ is a symmetric $M$-tuple with $M\leqslant5$. Although the relations among RTT generators of $\X(\mathfrak{g}_N)$ are much more complicated than those in type $A$, our new block matrix methodology
	enables us to conduct efficient calculations with matrix blocks to encapsulate parabolic generators.
	As a result, all relations in a parabolic presentation of $\X(\mathfrak{g}_N)$ will be expressed in a block matrix form as stated in Theorems~\ref{thm:odd:para} and ~\ref{thm:even:para}. Specifically, the parabolic presentation of $\X(\mathfrak{g}_N)$ reduces to its RTT presentation if $\nu=(N)$, and to its Drinfeld presentation if $\nu=(1,1,\ldots,1)$.  We believe our new method also applies to determining a parabolic presentation for Yangians of type $D$, but it requires more laborious calculations. We will address the type $D$ case separately in the future.
	
	The rest of this paper is organized as follows: In Section~\ref{sec:Yangian}, we briefly review the RTT presentation and properties of the extended Yangian $\X(\mathfrak{g}_N)$. Based on this, we present a new block matrix version of the RTT presentation and obtain a new set of generators of $\X(\mathfrak{g}_N)$ through a block Gauss decomposition of the generator matrix $T(u)$ of $\X(\mathfrak{g}_N)$. We also obtain a block version of the embedding theorem~\ref{embedding} and formulate part of parabolic relations that are the same as those in type $A$. Section~\ref{sec:para:odd} is dedicated to discussing the parabolic presentation of $\X(\mathfrak{g}_N)$ associated with a symmetric composition $\nu$ of odd length. We compute all the relations among the new generators, and prove that the parabolic presented Yangian is isomorphic to the RTT presented Yangian.  Section~\ref{sec:para:even} focuses on the parabolic presentation of $\X(\mathfrak{g}_N)$ associated with a symmetric composition $\nu$ of even length. The main results are stated in Theorems~\ref{thm:odd:para} and ~\ref{thm:even:para}. Finally, in Section~\ref{sec:center}, we derive a new expression for the center of $\X(\mathfrak{g}_N)$ by using the parabolic presentation.
	
	\section{Yangians of types BC}
	\label{sec:Yangian}
	\subsection{RTT presentation of $\X(\mathfrak{g}_{N})$ }
	
	In this subsection, we briefly review the RTT presentation of the extended Yangian $\X(\mathfrak{g}_{N})$ and $\Y(\mathfrak{g}_{N})$. The Lie algebra $\mathfrak{g}_{N}$ is the Lie subalgebra of $\mathfrak{gl}_{N}$ spanned by all elements
	\begin{equation}
		\label{fij}
		F_{i j}
		=e_{i j}-\varepsilon_i\varepsilon_j e_{j^\prime i^\prime},\quad\text{for }i,j=1,\ldots, N,
	\end{equation}
	where $e_{i j}$ denotes
	the standard basis element of $\mathfrak{gl}_{N}$,
	$i^\prime=N+1-i$, $\varepsilon_i=1$, $i=1,\cdots N$ for type $B_n$ and
	$$\varepsilon_i=\begin{cases}1,&\qquad\text{for $i=1,\cdots n$},\\
		-1,&\qquad\text{for $i=n+1,\cdots 2n$},
	\end{cases}$$
	for type $C_n$.
	
	In order to introduce the RTT presentation, we need the following notation: For every element
	\begin{equation}
		\label{mata}
		X=\sum_{i,j=1}^{N} X_{ij}\otimes e_{ij}\in \mathcal{A}\otimes\End\CC^{N},
	\end{equation}
	where $\mathcal{A}$ is a unital associative algebra, we denote
	\begin{equation}
		\label{xa}
		X^{[a]}=\sum_{i,j=1}^{N}  X_{ij}\otimes1^{\otimes(a-1)}\otimes e_{ij}\otimes 1^{\otimes(m-a)}
		\in \mathcal{A}\otimes\left(\End\CC^{N}\right)^{\otimes m},
	\end{equation}
	for  each $a\in\{1,\dots,m\}$, where $1$ is the identity endomorphism. Moreover, given an element
	\begin{equation}
		C=\sum_{i,j,k,l=1}^{N} c^{}_{ijkl}\, e_{ij}\otimes e_{kl}\in
		\End \CC^{N}\otimes\End \CC^{N},
	\end{equation}
	and two indices $a,b\in\{1,\dots,m\}$ such that $a<b$,
	we set
	\begin{equation}
		\label{ars}
		C^{ab}=\sum_{i,j,k,l=1}^N c^{}_{ijkl}\,
		1^{\otimes(a-1)}\otimes e_{ij}\otimes 1^{\otimes(b-a-1)}\otimes e_{kl}\otimes 1^{\otimes(m-b)}\in\left(\End\CC^{N}\right)^{\otimes m}.
	\end{equation}
	
	Let $R(u)$ be the following rational function in a complex parameter $u$
	with values in the tensor product algebra
	$\End\CC^{N}\otimes\End\CC^{N}$ given in \cite{zz:rf}:
	\begin{equation}
		\label{zamolr}
		R(u)=1-\frac{\mathsf{P}}{u}+\frac{\mathsf{Q}}{u-\kappa},
	\end{equation}
	where $\kappa=N/2\mp 1$,
	\begin{equation}
		\mathsf{P}=\sum_{i,j=1}^{N}e_{ij}\otimes e_{ji},\text{ and }
		\mathsf{Q}=\sum_{i,j=1}^{N} \varepsilon_i\varepsilon_j\, e_{ij}
		\otimes e_{i'j'}.
	\end{equation}
	The rational function \eqref{zamolr} satisfies the Yang--Baxter equation
	\begin{equation}
		\label{yberep}
		R^{12}(u-v)\, R^{13}(u)\, R^{23}(v)
		=R^{23}(v)\, R^{13}(u)\, R^{12}(u-v).
	\end{equation}
	
	The {\em extended Yangian}
	$\X(\mathfrak{g}_{N})$
	is a unital associative algebra with generators
	$t_{ij}^{(r)}$, where $1\leqslant i,j\leqslant N$ and $r=1,2,\dots$,
	satisfying certain quadratic relations. These generators can be encapsulated into the formal series
	\begin{equation}\label{tiju}
		t_{ij}(u)=\delta_{ij}+\sum_{r=1}^{\infty}t_{ij}^{(r)}\, u^{-r}
		\in\X(\mathfrak{g}_{N})[[u^{-1}]],
	\end{equation}
	and we set
	\begin{equation}
		T(u)=\sum_{i,j=1}^{N} t_{ij}(u)\otimes e_{ij}
		\in \X(\mathfrak{g}_N)[[u^{-1}]]\otimes \End\CC^N.
	\end{equation}
	The defining relations for the algebra $\X(\mathfrak{g}_{N})$ can be written in the matrix form
	\begin{equation}
		\label{RTTbcd}
		R^{12}(u-v)\, T^{[1]}(u)\, T^{[2]}(v)=T^{[2]}(v)\, T^{[1]}(u)\, R^{12}(u-v),
	\end{equation}
	which can also be written in terms of the series $t_{ij}(u)$ as
	\begin{equation}\label{defrel}
		\begin{aligned}
			[t_{ij}(u),t_{kl}(v)]&=\frac{1}{u-v}
			\left(t_{kj}(u) t_{il}(v)-t_{kj}(v) t_{il}(u)\right)\\
			&-\frac{1}{u-v-\kappa}
			\left(\delta_{ik^{\prime}}\sum_{p=1}^N\varepsilon_i\varepsilon_p t_{pj}(u) t_{p^{\prime}l}(v)-
			\delta_{jl^{\prime}}\sum_{p=1}^N\varepsilon_j\varepsilon_p t_{kp^{\prime}}(v) t_{ip}(u)\right).
		\end{aligned}
	\end{equation}
	
	The {\em Yangian}
	$\Y(\mathfrak{g}_{N})$ is defined as the subalgebra of
	$\X(\mathfrak{g}_{N})$ consisting of the elements stable under
	the automorphisms
	\begin{equation}
		\label{muf}
		\mu_f:T(u)\mapsto f(u)\, T(u),
	\end{equation}
	for all series
	$f(u)=1+f_1u^{-1}+f_2 u^{-2}+\cdots$
	with $f_i\in\mathbb{C}$. It is shown in \cite{amr:rp} that the extended Yangian $X(\mathfrak{g}_{N})$ admits the following tensor product decomposition:
	\begin{equation}
		\label{tensordecom}
		\X(\mathfrak{g}_{N})=\ZX(\mathfrak{g}_{N})\otimes \Y(\mathfrak{g}_{N}),
	\end{equation}
	where $\ZX(\mathfrak{g}_{N})$ is the center of the extended Yangian $\X(\mathfrak{g}_{N})$.
	The center $\ZX(\mathfrak{g}_{N})$ is generated by the coefficients of the series
	\begin{equation}
		\label{zn}
		z(u)=1+\sum_{r=1}^{\infty} z^{(r)}\, u^{-r},
	\end{equation}
	defined  by the equation:
	\begin{equation}
		\label{zcenter}
		T^t(u+\kappa)\, T(u)=T(u)\, T^t(u+\kappa)=z(u) 1,
	\end{equation}
	where
	$$T^t(u+\kappa)=\sum\limits_{i,j=1}^N\varepsilon_i\varepsilon_j  t_{ij}(u+\kappa)\otimes e_{j^{\prime}i^{\prime}}
	=\sum\limits_{i,j=1}^N\varepsilon_i\varepsilon_jt_{j^{\prime}i^{\prime}}(u+\kappa)\otimes e_{ij}.$$
	
	Equivalently,
	the Yangian $\Y(\mathfrak{g}_{N})$ is the quotient of $\X(\mathfrak{g}_{N})$
	by the relation $z(u)=1$, that is,
	\begin{equation}
		\label{unita}
		T^t(u+\kappa)\, T(u)=1.
	\end{equation}
	We refer the readers to \cite{aacfr:rp} and \cite{amr:rp} for more details on the structure of Yangians.
	
	\subsection{Block RTT relations}
	We consider the block form of the generator matrix associated to a composition$$\nu=(\nu_1,\nu_2,\ldots,\nu_{M})$$of $N$.
	We always assume that the composition $\nu$ is symmetric, i.e., $\nu_{a}=\nu_{a'}$ for $a'=M+1-i, a=1,\ldots,M$. In this situation, every $N\times N$-matrix can be regarded as a $M\times M$-block matrix in which the $(a,b)$-entry is a $\nu_a\times\nu_b$-submatrix. We use $\mathrm{M}_{\nu_a\times\nu_b}(\mathbb{C})$ to denote the set of all $\nu_a\times\nu_b$-matrices over $\mathbb{C}$, and simply write $\mathrm{M}_{\nu_a}(\mathbb{C})=\mathrm{M}_{\nu_a\times\nu_a}(\mathbb{C})$.
	
	In order to calculate with matrix blocks, we introduce the following notation. Let $\mathcal{A}$ be a unital associative algebra and $A_{ab}=(\alpha_{ij})_{1\leqslant i\leqslant \nu_a\atop 1\leqslant j\leqslant \nu_b}$ be a $\nu_a\times\nu_b$-matrix with entries $\alpha_{ij}\in\mathcal{A}$. We use  $A_{ab}\boxtimes e_{ab}$ to denote the block matrix associated to $\nu$ whose $(a,b)$-block is $A_{ab}$ and other blocks are $0$. Given a $\nu_a\times \nu_b$-matrix $A_{ab}$ for each pair $(a,b)$ with $a,b=1,2,\ldots,M$, we form an $N\times N$-matrix $A$ whose $(a,b)$-block is $A_{ab}$ for $a,b=1,2,\cdots,M$, that is
	$$A=\begin{pmatrix}A_{11}&A_{12}&\cdots&A_{1M}\\
		A_{21}&A_{22}&\cdots&A_{2M}\\
		\vdots&\vdots&&\vdots\\
		A_{M1}&A_{M2}&\cdots&A_{MM}
	\end{pmatrix}
	=\sum\limits_{a,b=1}^MA_{ab}\boxtimes e_{ab}.$$
	We observe that $\sum\limits_{a,b=1}^MA_{ab}\boxtimes e_{ab}=\sum\limits_{a,b=1}^MB_{ab}\boxtimes e_{ab}$ if and only if $A_{ab}=B_{ab}$ for $a,b=1,2,\ldots,M$. For the $\nu_a\times \nu_b$-matrix unit $e_{ij}$, $1\leq a,b\leq M$, $1\leq i\leq \nu_a$ and $1\leq j\leq \nu_b$, we identify
	\begin{equation*}
		e_{ij}\boxtimes e_{ab}=e_{\bar{\nu}_a+i,\bar{\nu}_b+j}\in\End\mathbb{C}^N,
	\end{equation*}
	where $\bar{\nu}_a=\nu_1+\nu_2+\cdots+\nu_{a-1}$. Then the matrix multiplication can be computed in the block form:
	\begin{equation*}
		\left(e_{ij}\boxtimes e_{ab}\right)\left(e_{kl}\boxtimes e_{cd}\right)=\delta_{bc}\delta_{kj}e_{il}\boxtimes e_{ad},
	\end{equation*}
	for $1\leqslant a,b,c,d\leqslant M$, $1\leqslant i\leqslant \nu_a, 1\leqslant j\leqslant\nu_b, 1\leqslant k\leqslant \nu_c$ and $1\leqslant l\leqslant\nu_d$.
	
	Using this notation, the generator matrix $T(u)$ can be written in the following block form:
	\begin{equation}
		\label{eq:Tublock}
		T(u)
		=\sum_{a,b=1}^{M}\sum_{\substack{i=1,\ldots,\nu_a \\ j=1,\ldots,\nu_b}}t_{\bar{\nu}_a+i,\bar{\nu}_b+j}(u)\otimes \left(e_{ij} \boxtimes e_{ab}\right)
		=\sum\limits_{a,b=1}^{M}T_{ab}(u)\boxtimes e_{ab},
	\end{equation}
	where $T_{ab}(u)=\sum\limits_{\substack{i=1,\ldots,\nu_a \\ j=1,\ldots,\nu_b}}t_{\bar{\nu}_a+i,\bar{\nu}_b+j}(u)\otimes e_{ij}\in \X(\mathfrak{g}_N)\otimes \mathrm{M}_{\nu_a\times \nu_b}(\mathbb{C})$ is the $(a,b)$-block in $T(u)$.
	
	We also write tensor matrices in block form. For a $\nu_a\times\nu_b$-matrix $A_{ab}$ and a $\nu_c\times\nu_d$-matrix $B_{cd}$, we denote the tensor product of two block matrices by
	$$(A_{ab}\otimes B_{cd})\boxtimes (e_{ab}\otimes e_{cd})=\left(A_{ab}\boxtimes e_{ab}\right)\otimes\left(B_{cd}\boxtimes e_{cd}\right).$$
	Then
	\begin{equation*}
		\left(e_{ij}\otimes e_{kl}\right)\boxtimes \left(e_{ab}\otimes e_{cd}\right)
		=e_{\bar{\nu}_a+i,\bar{\nu}_b+j}\otimes e_{\bar{\nu}_c+k,\bar{\nu}_d+l},
	\end{equation*}
	and
	$$\sum\limits_{a,b,c,d=1}^M(A_{ab}\otimes B_{cd})\boxtimes (e_{ab}\otimes e_{cd})
	=\sum\limits_{a,b,c,d=1}^M(\tilde{A}_{ab}\otimes \tilde{B}_{cd})\boxtimes (e_{ab}\otimes e_{cd})$$
	if and only if $A_{ab}\otimes B_{cd}=\tilde{A}_{ab}\otimes \tilde{B}_{cd}$ for all $a,b,c,d=1,2,\ldots, M$.
	Since $I_N=\sum\limits_{c=1}^MI_c\boxtimes e_{cc}$, where $I_c=\sum\limits_{k=1}^{\nu_c} e_{kk}\in \mathrm{M}_{\nu_c}(\mathbb{C})$, the tensor matrices $T^{[1]}(u)$ and $T^{[2]}(u)$ in \eqref{RTTbcd} are also accordingly written as block forms:
	\begin{align*}
		T^{[1]}(u)=\sum_{a,b=1}^{M}\sum_{c=1}^{M}T_{ab;c}^{[1]}(u)\boxtimes \left(e_{ab}\otimes e_{cc}\right),\quad
		T^{[2]}(u)=\sum_{a,b=1}^{M}\sum_{c=1}^MT_{ab;c}^{[2]}(u)\boxtimes \left(e_{cc}\otimes e_{ab}\right),
	\end{align*}
	where
	\begin{align*}
		T_{ab;c}^{[1]}(u)=&\sum\limits_{\substack{i=1,\ldots,\nu_a \\ j=1,\ldots,\nu_b}}
		t_{\bar{\nu}_a+i,\bar{\nu_b}+j}(u)\otimes e_{ij}\otimes I_c
		\in\X(\mathfrak{g}_N)\otimes \mathrm{M}_{\nu_a\times\nu_b}(\mathbb{C})\otimes \mathrm{M}_{\nu_c}(\mathbb{C}),\\
		T_{ab;c}^{[2]}(u)=&\sum\limits_{\substack{i=1,\ldots,\nu_a \\ j=1,\ldots,\nu_b}}
		t_{\bar{\nu}_a+i,\bar{\nu_b}+j}(u)\otimes I_c\otimes e_{ij}\in\X(\mathfrak{g}_N)\otimes \mathrm{M}_{\nu_c}(\mathbb{C})\otimes \mathrm{M}_{\nu_a\times\nu_b}(\mathbb{C}).
	\end{align*}
	The multiplication $T_{ab;p}^{[1]}(u)T_{cd;q}^{[2]}(v)$ is admissible if and only if $p=c$ and $q=b$, and
	$T_{cd;q}^{[2]}(v)T_{ab;p}^{[1]}(u)$ is admissible if and only if $p=d$ and $q=a$. Hence, we simply write admissible products
	$$T_{ab}^{[1]}(u)T_{cd}^{[2]}(v)=T_{ab;c}^{[1]}(u)T_{cd;b}^{[2]}(v),\text{ and }
	T_{cd}^{[2]}(v)T_{ab}^{[1]}(u)=T_{cd;a}^{[2]}(v)T_{ab;d}^{[1]}(u).$$
	Thus, we can get
	\begin{equation}\label{eq:T1T2}
		\begin{aligned}
			T^{[1]}(u)T^{[2]}(v)&=\sum_{a,b,c,d=1}^{M}T^{[1]}_{ab}(u)T_{cd}^{[2]}(v)\boxtimes \left(e_{ab}\otimes e_{cd}\right),\\
			T^{[2]}(v)T^{[1]}(u)&  =\sum_{a,b,c,d=1}^{M}T^{[2]}_{cd}(v)T_{ab}^{[1]}(u)\boxtimes \left(e_{ab}\otimes e_{cd}\right).
		\end{aligned}
	\end{equation}
	
	We also write both tensor factors in the R-matrix \eqref{zamolr} in the block form, then
	\begin{equation}
		\label{eq:blockR}
		R(u)=\sum_{a,b=1}^{M}I_{ab}\boxtimes \left(e_{aa}\otimes e_{bb}\right)-\frac{1}{u}\sum_{a,b=1}^{M} P_{ab}\boxtimes \left(e_{ab}\otimes e_{ba}\right)+
		\frac{1}{u-\kappa}\sum_{a,b=1}^{M}Q_{ab}\boxtimes \left(e_{ab}\otimes e_{a'b'}\right),
	\end{equation}
	where
	$$I_{ab}=\sum_{\substack{i=1,\ldots,\nu_a \\ j=1,\ldots,\nu_b}}e_{ii}\otimes e_{jj},\quad P_{ab}=\sum_{\substack{i=1,\ldots,\nu_a \\ j=1,\ldots,\nu_b}} e_{ij}\otimes e_{ji}, \quad
	Q_{ab}=\sum_{\substack{i=1,\ldots,\nu_a \\ j=1,\ldots,\nu_b}} \epsilon_{i}^a\epsilon_{j}^b e_{ij}\otimes e_{\nu_a+1-i,\nu_b+1-j},$$
	and the sign $\epsilon_i^a=1, a=1,\ldots, M, i=1,\ldots,\nu_a$ if  $\mathfrak{g}_N$ is of type B; For type $C$, the sign $\epsilon_i^a$ is set as follows:
	\begin{enumerate}
		\item If $M=2m+1$, then
		$$\epsilon_{i}^a=
		\begin{cases}1,&\text{if  $a< m+1$ or $a=m+1$ and $i\leqslant \frac{\nu_{m+1}}{2}$},\\
			-1&\text{if  $a> m+1$ or $a=m+1$ and $i> \frac{\nu_{m+1}}{2}$}.
		\end{cases}$$
		\item If $M=2m$, then
		$$\epsilon_{i}^a=
		\begin{cases}1,&\text{if $a=1,\cdots m$},\\
			-1,&\text{if $a=m+1,\cdots 2m$}.
		\end{cases}$$
	\end{enumerate}
	
	\begin{thm}
		The RTT-relation \eqref{RTTbcd} is equivalent to the following relation written in terms of submatrices of $T(u)$:
		\begin{equation}
			\label{eq:blockRTT}
			\begin{aligned}
				\left[T_{ab}^{[1]}(u),\right.&\left.T_{cd}^{[2]}(v)\right]=\frac{1}{u-v}\left(P_{ac}T_{cb}^{[1]}(u)T_{ad}^{[2]}(v)-
				T_{cb}^{[2]}(v)T_{ad}^{[1]}(u)P_{db}\right)\\
				&-\frac{1}{u-v-\kappa}\left(\delta_{ac^{\prime}}\sum_{r=1}^{M}Q_{ar}T_{rb}^{[1]}(u)T_{r'd}^{[2]}(v)
				-\delta_{bd'}\sum_{r=1}^{M}T_{cr'}^{[2]}(v)T_{ar}^{[1]}(u)Q_{rb}\right),
			\end{aligned}
		\end{equation}
		where $a,b,c,d=1,2,\ldots,M$.
	\end{thm}
	
	\begin{proof}
		We write all factors in \eqref{RTTbcd} using \eqref{eq:T1T2} and \eqref{eq:blockR}. By comparing block coefficients of $e_{ab}\otimes e_{cd}$ in the both sides of \eqref{RTTbcd}, we have
		\begin{equation*}
			\begin{aligned}
				T_{ab}^{[1]}(u)&T_{cd}^{[2]}(v)-\frac{1}{u-v}P_{ac}T_{cb}^{[1]}(u)T_{ad}^{[2]}(v)+\frac{1}{u-v-\kappa}
				\delta_{ac^{\prime}}\sum_{r=1}^{M}Q_{ar}T_{rb}^{[1]}(u)T_{r'd}^{[2]}(v)\\
				&=T_{cd}^{[2]}(v)T_{ab}^{[1]}(u)-\frac{1}{u-v}T_{cb}^{[2]}(v)T_{ad}^{[1]}(u)P_{db}+\frac{1}{u-v-\kappa}
				\delta_{bd'}\sum_{r=1}^{M}T_{cr'}^{[2]}(v)T_{ar}^{[1]}(u)Q_{rb},
			\end{aligned}
		\end{equation*}
		which is equivalent to the equation \eqref{eq:blockRTT}.
	\end{proof}
	\begin{remark}
		The block RTT-relation \eqref{eq:blockRTT} reduces to the RTT-relation \eqref{RTTbcd} if the composition $\nu=(N)$, and reduces to the relation \eqref{defrel} if $\nu=(1,1,\ldots,1)$.
	\end{remark}
	
	For convenience in computation,
	the RTT-relation \eqref{RTTbcd} is also written as
	\begin{equation}
		T^{[1]}(u)R^{12}(u-v)\widetilde{T}^{[2]}(v)=\widetilde{T}^{[2]}(v)R^{12}(u-v)T^{[1]}(u),\label{TRT}
	\end{equation}
	where $\widetilde{T}(u)=\sum\limits_{i,j=1}^{N}\tilde{t}_{ij}(u)\otimes e_{ij}$ is the inverse matrix of $T(u)$. The matrix $\widetilde{T}(u)$ also has the block form
	$$\widetilde{T}(u)=\sum_{a,b=1}^{M}\widetilde{T}_{ab}(u)\boxtimes e_{ab}, \text{ where }\widetilde{T}_{ab}(u)=\sum_{\substack{i=1,\ldots,\nu_a \\ j=1,\ldots,\nu_b}}\widetilde{T}_{ab;ij}(u)\otimes e_{ij}.$$
	Then the relation \eqref{TRT} is equivalent to
	\begin{equation}
		\label{blockTTtilde}
		\begin{aligned}
			\left[T_{ab}^{[1]}(u),\widetilde{T}_{cd}^{[2]}(v)\right]=&\frac{1}{u-v}\sum_{r=1}^{M}\left(\delta_{cb}T_{ar}^{[1]}(u)P_{rb}\widetilde{T}_{rd}^{[2]}(v)-
			\delta_{ad}\widetilde{T}_{cr}^{[2]}(v)P_{ar}T_{rb}^{[1]}(u)\right)\\
			&+\frac{1}{u-v-\kappa}\left(\widetilde{T}_{ca'}^{[2]}(v)Q_{ad'}T_{d'b}^{[1]}(u)
			-T_{ac'}^{[1]}(u)Q_{c'b}\widetilde{T}_{b'd}^{[2]}(v)\right).
		\end{aligned}
	\end{equation}
	where the products appearing in the above equation are admissible products, for example
	$
	T_{ar}^{[1]}(u)P_{rb}\widetilde{T}_{rd}^{[2]}(v)$ means the admissible product $T_{ar;b}^{[1]}(u)P_{rb}\widetilde{T}_{rd;b}^{[2]}(v)$.
	\begin{lem}
		\label{lem:PQ}
		The tensor matrices $P_{ab}$ and $Q_{ab}$, $a,b=1,2,\ldots, m$, satisfy the following identities:
		\begin{align}
			P_{ab}P_{ba}=&I_{ab},& Q_{ab}Q_{bc}=&\nu_bQ_{ac},\label{eq:PPQQ}\\
			P_{aa}Q_{ab}=&\pm Q_{ab},&Q_{ba}P_{aa}=&\pm Q_{ba},\label{eq:PQQP}
		\end{align}
		where the negative sign in \eqref{eq:PQQP} only appears if $\mathfrak{g}_N$ is of type $C$, $M=2m+1$ and $a=m+1$.
	\end{lem}
	\begin{proof}
		All identities can be verified directly.
	\end{proof}
	
	\begin{lem}
		\label{lem:PQtrs}
		Let $t$ be the linear operator on $\nu_a\times\nu_b$-matrices such that $e_{ij}^t=\epsilon_{i}^a \epsilon_{j}^be_{v_b+1-j, v_a+1-i}$ if $e_{ij}$ is a $\nu_a\times \nu_b$ matrix unit. Then the following identities hold:
		\begin{equation}
			P_{ab}^{t_1}=\pm Q_{ba},\quad
			P_{ab}^{t_2}=Q_{ab},\quad
			Q_{ab}^{t_1}=P_{ba},\quad
			Q_{ab}^{t_2}=\pm P_{ab},\label{eq:PtrQtr}
		\end{equation}
		where $t_1, t_2$ mean that the map $t$ is applied to the first or second tensor factor; the negative sign in \eqref{eq:PtrQtr} merely appears if $\mathfrak{g}_N$ is of type $C$, $M=2m+1$ and either $a$ or $b$ is $m+1$. \qed
	\end{lem}
	
	We will use the following lemma often in Section \ref{sec:para:odd} and Section \ref{sec:para:even}.
	
	\begin{lem}
		\label{lem:APQ}
		Let $A=\sum\limits_{i=1,\ldots,\nu_a\atop j=1,\ldots,\nu_b}\alpha_{ij}\otimes e_{ij}$ be a $\nu_a\times\nu_b$ matrix with entries in an arbitrary unital algebra $\mathcal{A}$. We denote
		$$A^{[1]}_p=\sum_{i=1,\ldots,\nu_a\atop j=1,\ldots,\nu_b}\alpha_{ij}\otimes e_{ij}\otimes I_p,\text{ and }
		A^{[2]}_p
		=\sum_{i=1,\ldots,\nu_a\atop j=1,\ldots,\nu_b}
		\alpha_{ij}\otimes I_p\otimes e_{ij}.$$
		Then the following identities hold:
		\begin{equation}\label{eq:PA1}
			P_{ca}A^{[1]}=A^{[2]}P_{cb},\quad
			P_{ac}A^{[2]}=A^{[1]}P_{bc},
		\end{equation}
		\begin{equation}\label{eq:A1Q}
			Q_{ca}A^{[1]}=Q_{cb}A^{t, [2]},\quad
			A^{[1]}Q_{bc}=A^{t, [2]}Q_{ac}.
		\end{equation}
		\begin{equation}\label{eq:A2Q}
			Q_{ca}A^{[2]}=\pm Q_{cb}A^{t, [1]},\quad
			A^{[2]}Q_{bc}=\pm A^{t, [1]}Q_{ac},
		\end{equation}
		the negative sign in \eqref{eq:PtrQtr} merely appears if $\mathfrak{g}_N$ is of type $C$, $M=2m+1$ and either $a$ or $b$ is $m+1$.
		
		If we assume in addition that $\nu_a=\nu_b$, then
		$$Q_{ca}A^{[1]}Q_{ad}=Q_{ca}A^{[2]}Q_{ad}=\mathrm{tr}(A)Q_{cd}.$$
		where $\mathrm{tr}(A)=\alpha_{11}+\alpha_{22}+\cdots+\alpha_{\nu_a,\nu_a}$.
		The products occurring in the above equations are all admissible. For instance,
		the equation $P_{ca}A^{[1]}=A^{[2]}P_{cb}$ means $P_{ca}A^{[1]}_c=A^{[2]}_cP_{cb}$.
	\end{lem}
	\begin{proof}
		All equalities can be directly verified with a calculation of tensor matrices.
	\end{proof}
	
	\subsection{Block Gauss Decomposition and embedding theorem}
	
	In order to formulate parabolic generators of the extended Yangian $\X(\mathfrak{g}_{N})$, we need to develop a block Gauss decomposition of its generator matrix $T(u)$. For a fixed symmetric composition $\nu=(\nu_1,\nu_2,\ldots,\nu_{M})$ of $N$,  the generator matrix $T(u)$ of $\X(\mathfrak{g}_{N})$ is factored as
	\begin{equation}
		\label{blockGaussDec}
		T(u)=\mathcal{F}(u)\mathcal{D}(u)\mathcal{E}(u),
	\end{equation}
	for block matrices
	\begin{equation}
		\mathcal{F}(u)=\begin{bmatrix}
			I_{\nu_1}&0&\cdots&0\,\\
			\mathcal{F}_{21}(u)&I_{\nu_2}&\cdots&0\\
			\vdots&\vdots&\ddots&\vdots\\
			\mathcal{F}_{M,1}(u)&\mathcal{F}_{M,2}(u)&\cdots&I_{\nu_{M}}
		\end{bmatrix},
		\mathcal{E}(u)=\begin{bmatrix}
			I_{\nu_1}&\mathcal{E}_{12}(u)&\cdots&\mathcal{E}_{1,M}(u)\,\\
			0&I_{\nu_2}&\dots&\mathcal{E}_{2,M}(u)\\
			\vdots&\vdots&\ddots&\vdots\\
			0&0&\cdots&I_{\nu_{M}}
		\end{bmatrix},
	\end{equation}
	and $$\mathcal{D}(u)=\diag\,\left[\mathcal{D}_1(u),\dots,\mathcal{D}_{M}(u)\right],$$
	where $\mathcal{D}_a(u)=\left(\mathcal{D}_{a,ij}(u)\right)$, $\mathcal{E}_{ab}(u)=\left(\mathcal{E}_{ab,ij}(u)\right)$ and $\mathcal{F}_{ab}(u)=\left(\mathcal{F}_{ab,ij}(u)\right)$
	are $\nu_a\times\nu_a, \nu_a\times\nu_b$ and $\nu_a\times\nu_b$-matrices respectively.
	
	In fact, the block matrix $\mathcal{D}_a(u)$, $\mathcal{E}_{ab}(u)$ and $\mathcal{F}_{ab}(u)$ can be explicitly described in terms of quasi-determinants \cite{gr:dm,gr:gauss}:
	\begin{align*}
		\mathcal{D}_a(u)
		=&\begin{vmatrix}
			T_{1,1}(u)&\cdots&T_{1,a-1}(u)&T_{1,a}(u)\\
			\vdots&\ddots&\vdots&\vdots\\
			T_{a-1,1}(u)&\cdots&T_{a-1,a-1}(u)&T_{a-1,a}(u)\\
			T_{a,1}(u)&\cdots&T_{a,a-1}(u)&\boxed{T_{a,a}(u)}
		\end{vmatrix},\\
		\mathcal{E}_{a,b}(u)
		=&\mathcal{D}_a(u)^{-1}\begin{vmatrix}
			T_{1,1}(u)&\cdots&T_{1,a-1}(u)&T_{1,b}(u)\\
			\vdots&\ddots&\vdots&\vdots\\
			T_{a-1,1}(u)&\cdots&T_{a-1,a-1}(u)&T_{a-1,b}(u)\\
			T_{a,1}(u)&\cdots&T_{a,a-1}(u)&\boxed{T_{a,b}(u)}
		\end{vmatrix},\\
		\mathcal{F}_{b,a}(u)
		=&\begin{vmatrix}
			T_{1,1}(u)&\cdots&T_{1,a-1}(u)&T_{1,a}(u)\\
			\vdots&\ddots&\vdots&\vdots\\
			T_{a-1,1}(u)&\cdots&T_{a-1,a-1}(u)&T_{a-1,a}(u)\\
			T_{b,1}(u)&\cdots&T_{b,a-1}(u)&\boxed{T_{b,a}(u)}
		\end{vmatrix}\mathcal{D}_a(u)^{-1},\\
	\end{align*}
	where $a<b$ and the quasi determinant is defined by
	$$\begin{vmatrix}
		A&B\\
		C&\boxed{D}
	\end{vmatrix}
	=D-CA^{-1}B$$
	for matrices $A$, $B$, $C$, $D$ with entries in a ring provided that $A$ is invertible.
	
	We simply denote
	\begin{equation}
		\mathcal{E}_a(u)=\mathcal{E}_{a,a+1}(u),\text{ and }\mathcal{F}_a(u)=\mathcal{F}_{a+1,a}(u)
	\end{equation}
	for $a=1,\ldots,M-1$. A new family of generators for the extended Yangian $\X(\mathfrak{g}_N)$ can be described as follows:
	
	\begin{prop}
		\label{prop:generators}
		Let $m=\left[\frac{M}{2}\right]$. The coefficients of formal series $\mathcal{D}_{a,ij}(u)$ with $a=1,\ldots, m+1, i,j=1,\ldots,\nu_a$, $\mathcal{E}_{b,ij}(u)$ and $\mathcal{F}_{b,ji}(u)$ with $b=1,\ldots, m$, $i=1,\ldots,\nu_b$, $j=1,\ldots,\nu_{b+1}$ generate the extended Yangian $\X(\mathfrak{g}_{N})$ as an associative algebra.
	\end{prop}
	\begin{proof}
		We observe that each diagonal block $\mathcal{D}_a(u)$ is invertible. Hence, $\mathcal{D}_a(u)$ admits a Gauss decomposition
		$$\mathcal{D}_a(u)=\mathcal{D}_a^-(u)\mathsf{H}_a(u)\mathcal{D}_a^+(u),$$
		where $\mathcal{D}_a^-(u)$ (resp. $\mathcal{D}_a^+(u)$) is a lower-triangular (resp. upper-triangular) matrix with all diagonal entries $1$, and $H_a(u)$ is a diagonal matrix. Hence,
		\begin{align*}
			T(u)=\mathcal{F}(u)\mathcal{D}(u)\mathcal{E}(u)=\mathsf{F}(u)\mathsf{H}(u)\mathsf{E}(u),
		\end{align*}
		where
		\begin{equation}
			\label{generating}
			\begin{aligned}
				\mathsf{H}(u)=&\mathrm{diag}\left(\mathsf{H}_1(u),\ldots,\mathsf{H}_{M}(u)\right),\\
				\mathsf{E}(u)=&\mathrm{diag}\left(\mathcal{D}_1^{+}(u),\ldots,\mathcal{D}_{M}^+(u)\right)\mathcal{E}(u),\\
				\mathsf{F}(u)=&\mathcal{F}(u)\mathrm{diag}\left(\mathcal{D}_1^{-}(u),\ldots,\mathcal{D}_{M}^-(u)\right),
			\end{aligned}
		\end{equation}
		are diagonal, upper-triangular, and lower-triangular matrices respectively, and all diagonal entries of $\mathsf{E}(u)$ and $\mathsf{F}(u)$ are $1$. By the uniqueness of Gauss decomposition, we obtain that
		\begin{equation}
			\label{GaussDec}
			T(u)=\mathsf{F}(u)\mathsf{H}(u)\mathsf{E}(u)
		\end{equation}
		coincides with the standard Gauss decomposition of the generator matrix $T(u)$.
		
		If we write the entries of these matrices as
		$$\mathsf{H}(u)=\mathrm{diag}\left[h_1(u),\ldots,h_{N}(u)\right],\quad \mathsf{E}(u)=(e_{ij}(u))_{N\times N},\qquad \mathsf{F}(u)=(f_{ij}(u))_{N\times N},$$
		then it is shown in \cite{jlm:ib} that the extended Yangian $\X(\mathfrak{g}_{N})$ is generated by the coefficients of $h_i(u)$ with $i=1,\ldots,n+1$, $e_{j,j+1}(u)$ and $f_{j+1,j}(u)$ with $j=1,\ldots, n$, where $n=\left[\frac{N}{2}\right]$. These elements are entries of $\mathsf{H}_a(u), \mathcal{D}_b^+(u)\mathcal{E}_{b}(u), \mathcal{F}_b(u)\mathcal{D}_b^-(u)$ for some $a=1,\ldots, m+1$ and $b=1,\ldots,m$. We further observe from the Gauss decomposition $\mathcal{D}_a(u)=\mathcal{D}_a^-(u)\mathsf{H}_a(u)\mathcal{D}_a^+(u)$ that all entries of $\mathcal{D}_a^-(u)$, $\mathsf{H}_a(u)$ and $\mathcal{D}_a^+(u)$ are generated by entries of $\mathcal{D}_a(u)$. Hence, the coefficients of all entries of $\mathcal{D}_a(u)$ with $a=1,\ldots, m+1$, $\mathcal{E}_b(u)$ and $\mathcal{F}_b(u)$, $b=1,\ldots, m$ generate the associative algebra $\X(\mathfrak{g}_{N})$.
	\end{proof}
	
	In order to figure out the relations among the generators obtained in Proposition~\ref{prop:generators}, we need a block version of the embedding theorem obtained in \cite[Theorem~3.7]{jlm:ib}.
	
	Recall that the generator matrix $T(u)$ for $\X(\mathfrak{g}_{N})$ is factored as in \eqref{GaussDec}:
	$$T(u)=\mathsf{F}(u)\mathsf{H}(u)\mathsf{E}(u).$$
	For each $1\leqslant\ell< n$, we consider the sub-matrices of $\mathsf{E}(u), \mathsf{H}(u)$ and $\mathsf{F}(u)$ corresponding to the rows and columns labeled by $\ell+1,\ldots, n,n+1\ldots, (\ell+1)'$, namely,
	\begin{equation}
		\mathsf{F}^{[\ell]}(u)=\begin{bmatrix}
			1&0&\dots&0\,\\
			f_{\ell+2\, \ell+1}(u)&1&\dots&0\\
			\vdots&\ddots&\ddots&\vdots\\
			f_{(\ell+1)' \ell+1}(u)&\dots&f_{(\ell+1)'\, (\ell+2)'}(u)&1
		\end{bmatrix},
	\end{equation}
	\begin{equation}
		\mathsf{E}^{[\ell]}(u)=\begin{bmatrix} 1&e_{\ell+1 \ell+2}(u)&\ldots&e_{\ell+1 (\ell+1)'}(u)\\
			0&1&\ddots &\vdots\\
			\vdots&\vdots&\ddots&e_{(\ell+2)'(\ell+1)'}(u)\\
			0&0&\ldots&1\\
		\end{bmatrix},
	\end{equation}
	and $\mathsf{H}^{[\ell]}(u)=\diag\,\big[h_{\ell+1}(u),\dots,h_{(\ell+1)'}(u)\big]$. It is easy to observe that $\mathsf{F}^{[\ell]}(u)\mathsf{H}^{[\ell]}(u)\mathsf{E}^{[\ell]}(u)$ is not a sub-matrix of $T(u)$. Nonetheless, it was shown in \cite[Theorem~3.7]{jlm:ib} that there exists a unique injective homomorphism
	\begin{equation}
		\label{embed}
		\psi_{\ell}: \X(\mathfrak{g}_{N-2\ell})\rightarrow \X(\mathfrak{g}_{N})
	\end{equation}
	such that
	$$\psi_{\ell}(T(u))=\mathsf{F}^{[\ell]}(u)\mathsf{H}^{[\ell]}(u)\mathsf{E}^{[\ell]}(u),$$
	where $T(u)=\sum\limits_{i,j=\ell+1}^{(\ell+1)'}t_{ij}^{(r)}(u)\otimes e_{ij}$ is the generator matrix of $\X(\mathfrak{g}_{N-2\ell})$ whose rows and columns are labeled by $\ell+1,\ldots, n,n+1,\ldots, (\ell+1)^{\prime}$.
	
	This result can be reformulated in terms of sub-block-matrices. For a fixed symmetric composition $\nu=(\nu_1,\nu_2,\ldots,\nu_M)$ of $N$, the generator matrix $T(u)$ is factored as a product of block matrices as in \eqref{blockGaussDec}:
	$$T(u)=\mathcal{F}(u)\mathcal{D}(u)\mathcal{E}(u).$$
	We consider the sub-block-matrices of $\mathcal{F}(u)$, $\mathcal{D}(u)$, and $\mathcal{E}(u)$ whose block rows and columns are labeled by $p+1,\ldots,m, m+1,\ldots, (p+1)'$, that is,
	\begin{equation}
		\mathcal{F}^{[p]}(u)=\begin{bmatrix}
			I_{\nu_{p+1}}&0&\dots&0\,\\
			\mathcal{F}_{p+2\, p+1}(u)&I_{\nu_{p+2}}&\dots&0\\
			\vdots&\ddots&\ddots&\vdots\\
			\mathcal{F}_{(p+1)' p+1}(u)&\dots&\mathcal{F}_{(p+1)'\, (p+2)'}(u)&I_{\nu_{p+1}}
		\end{bmatrix},
	\end{equation}
	\begin{equation}
		\mathcal{E}^{[p]}(u)=\begin{bmatrix} I_{\nu_{p+1}}&\mathcal{E}_{p+1 p+2}(u)&\ldots&\mathcal{E}_{p+1 (p+1)'}(u)\\
			0&I_{\nu_{p+2}}&\ddots &\vdots\\
			\vdots&\vdots&\ddots&\mathcal{E}_{(p+2)'(p+1)'}(u)\\
			0&0&\ldots&I_{\nu_{p+1}}\\
		\end{bmatrix}
	\end{equation}
	and $\mathcal{D}^{[p]}(u)=\diag\left[\mathcal{D}_{p+1}(u),\ldots,\mathcal{D}_{(p+1)'}(u)\right]$. Note that each block $\mathcal{D}_a(u), a=1,\ldots, M$ is further factored as
	$$\mathcal{D}_a(u)=\mathcal{D}_a^-(u)\mathsf{H}_a(u)\mathcal{D}_a^+(u),$$
	where $\mathsf{H}_a(u)$ is a diagonal matrix, $\mathcal{D}_a^-(u)$ and $\mathcal{D}_a^+(u)$ are lower-triangular and upper-triangular unipotent matrices respectively. Hence,
	\begin{align*}
		\mathsf{F}^{[\ell_p]}(u)=&\mathcal{F}^{[p]}(u)\diag\left[\mathcal{D}_{p+1}^-(u),\ldots,\mathcal{D}_{(p+1)'}^-(u)\right],
		\text{ and }\\
		\mathsf{E}^{[\ell_p]}(u)=&\diag\left[\mathcal{D}_{p+1}^+(u),\ldots,\mathcal{D}_{(p+1)'}^+(u)\right]\mathcal{E}^{[p]}(u),
	\end{align*}
	where $\ell_p=\nu_1+\cdots+\nu_p$. Consequently, we have the following proposition:
	\begin{thm}[A block version of {\cite[Theorem~3.1]{jlm:ib}}]
		\label{embedding}
		Let $\nu=(\nu_1,\nu_2,\ldots,\nu_{M})$ be a symmetric composition of $N$. For each $1\leqslant p\leqslant m$, there exists a unique injective algebra homomorphism
		$$\varPsi_p: \X(\mathfrak{g}_{N-2\ell_p})\rightarrow \X(\mathfrak{g}_{N})$$
		such that
		$$\varPsi_p(T(u))=\mathcal{F}^{[p]}(u)\mathcal{D}^{[p]}(u)\mathcal{E}^{[p]}(u),$$
		where $\ell_p=\nu_1+\cdots+\nu_p$ and
		$$T(u)=\sum\limits_{i,j=\ell_p+1}^{(\ell_p+1)'}t_{ij}^{(r)}(u)\otimes e_{ij}=\sum\limits_{a,b=p+1}^{(p+1)'}T_{ab}(u)\otimes E_{ab}$$
		is the generator matrix of $\X(\mathfrak{g}_{N-2\ell_p})$ whose block rows and columns are labeled by $p+1$, $p+2$, $\ldots$, $(p+1)'$. \qed
	\end{thm}
	
	Proposition~\ref{embedding} guarantees that the entries of $\Psi_p(T(u))$ generate a subalgebra of $\X(\mathfrak{g}_{N})$ that is isomorphic to $\X(\mathfrak{g}_{N-2\ell_p})$. Therefore, $\Psi_p(T(u))$ satisfies the RTT relation in  $\X(\mathfrak{g}_{N-2\ell_p})$:
	\begin{equation}
		\label{embedRTT}
		R^{12}(u-v)\varPsi_p(T^{[1]}(u))\varPsi_p(T^{[2]}(v))
		=\varPsi_p(T^{[2]}(v))\varPsi_p(T^{[1]}(u))R^{12}(u-v),
	\end{equation}
	Moreover, it follows from \cite[Corollary~3.10]{jlm:ib} that
	\begin{equation}
		\label{eq:TTpcomm}
		\left[\varPsi_{q}\left(T_{ab}^{[1]}(u)\right),\varPsi_{p}\left(T_{cd}^{[2]}(v)\right)\right]=0
	\end{equation}
	for all $q+1\leqslant a,b\leqslant p$, $p+1\leqslant c,d\leqslant (p+1)^{\prime}$ and $q<p$.
	
	\begin{lem}
		\label{lem:EFTp}
		Suppose that $p+1\leqslant a,b,c\leqslant (p+1)^{\prime}$ and $a\neq c'$.
		Then the following relations hold in the extended Yangian
		$\X(\mathfrak{g}_{N})$,
		\begin{align}
			\label{eq:ETp}
			\left[\mathcal{E}_{pa}^{[1]}(u), \varPsi_p\left(T_{bc}^{[2]}(v)\right)\right]
			=&\varPsi_p\left(T_{ba}^{[2]}(v)\right)\left(\mathcal{E}_{pc}^{[1]}(v)-\mathcal{E}_{pc}^{[1]}(u)\right)\frac{P_{ca}}{u-v},\\
			\label{eq:FTp}
			\left[\mathcal{F}_{a p}^{[1]}(u), \varPsi_p\left(T_{c b}^{[2]}(v)\right)\right]
			=&\frac{P_{ac}}{u-v}\left(\mathcal{F}_{c p}^{[1]}(u)-\mathcal{F}_{c p}^{[1]}(v)\right)\varPsi_p\left(T_{a b}^{[2]}(v)\right).
		\end{align}
	\end{lem}
	\begin{proof}
		It is sufficient to verify the relations for $p = 1$ since the general case will follow by applying the homomorphism $\varPsi_{p-1}$. Both relations follow from similar
		arguments so we only verify \eqref{eq:ETp}.
		
		By the block Gauss decomposition \eqref{blockGaussDec}, we have
		\begin{equation}
			T_{bc}(u)=\mathcal{F}_{b1}(u)\mathcal{D}_1(u)\mathcal{E}_{1c}(u)+\varPsi_1\left(T_{bc}(u)\right).
		\end{equation}
		Note that $a\neq c'$, we consider the block RTT relation
		\begin{equation}
			\left[T_{1a}^{[1]}(u),T_{bc}^{[2]}(v)\right]
			=\frac{1}{u-v}\left(P_{1b}T_{ba}^{[1]}(u)T_{1c}^{[2]}(v)-T_{ba}^{[2]}(v)T_{1c}^{[1]}(u)P_{ca}\right),
		\end{equation}
		which yields that
		\begin{align*}
			&\left[T_{1a}^{[1]}(u),\varPsi_1\left(T_{bc}^{[2]}(v)\right)\right]
			+\left[T_{1a}^{[1]}(u),\mathcal{F}_{b1}^{[2]}(v)\mathcal{D}_1^{[2]}(v)\mathcal{E}_{1c}^{[2]}(v)\right]\\
			=&\frac{1}{u-v}\left(P_{1b}\varPsi_1\left(T_{ba}^{[1]}(u)\right)T_{1c}^{[2]}(v)-\varPsi_1\left(T_{ba}^{[2]}(v)\right)T_{1c}^{[1]}(u)P_{ca}\right)\\
			&+\frac{1}{u-v}\left(P_{1b}\mathcal{F}_{b1}^{[1]}(u)\mathcal{D}_{1}^{[1]}(u)\mathcal{E}_{1a}^{[1]}(u)T_{1c}^{[2]}(v)
			-\mathcal{F}_{b1}^{[2]}(v)\mathcal{D}_1^{[2]}(v)\mathcal{E}_{1a}^{[2]}(v)T_{1c}^{[1]}(u)P_{ca}\right).
		\end{align*}
		We also compute that
		\begin{align*}
			&\left[T_{1a}^{[1]}(u),\mathcal{F}_{b1}^{[2]}(v)\mathcal{D}_1^{[2]}(v)\mathcal{E}_{1c}^{[2]}(v)\right]\\
			=&\left[T_{1a}^{[1]}(u),T_{b1}^{[2]}(v)\right]\mathcal{E}_{1c}^{[2]}(v)
			-\mathcal{F}_{b1}^{[2]}(v)\left[T_{1a}^{[1]}(u),T_{11}^{[2]}(v)\right]\mathcal{E}_{1c}^{[2]}(v)
			+\mathcal{F}_{b1}^{[2]}(v)\left[T_{1a}^{[1]}(u),T_{1c}^{[2]}(v)\right]\\
			=&\frac{1}{u-v}\left(P_{1b}T_{ba}^{[1]}(u)T_{11}^{[2]}(v)-T_{ba}^{[2]}(v)T_{11}^{[1]}(u)P_{1a}\right)\mathcal{E}_{1c}^{[2]}(v)\\
			&+\frac{1}{u-v}\mathcal{F}_{b1}^{[2]}(v)\left(T_{1a}^{[2]}(v)T_{11}^{[1]}(u)P_{1a}-P_{11}T_{1a}^{[1]}(u)T_{11}^{[2]}(v)\right)\mathcal{E}_{1c}^{[2]}(v)\\
			&+\frac{1}{u-v}\mathcal{F}_{b1}^{[2]}(v)\left(P_{11}T_{1a}^{[1]}(u)T_{1c}^{[2]}(v)-T_{1a}^{[2]}(v)T_{1c}^{[1]}(u)P_{ca}\right)\\
			=&\frac{1}{u-v}\left(P_{1b}T_{ba}^{[1]}(u)T_{1c}^{[2]}(v)
			-\varPsi_1\left(T_{ba}^{[2]}(v)\right)T_{11}^{[1]}(u)P_{1a}\mathcal{E}_{1c}^{[2]}(v)
			-\mathcal{F}_{b1}^{[2]}(v)T_{1a}^{[2]}(v)T_{1c}^{[1]}(u)P_{ca}\right).
		\end{align*}
		This yields
		\begin{equation*}
			\left[T_{1a}^{[1]}(u),\varPsi_1\left(T_{bc}^{[2]}(v)\right)\right]
			=\frac{1}{u-v}\varPsi_1\left(T_{ba}^{[2]}(v)\right)
			\left(T_{11}^{[1]}(u)P_{1a}\mathcal{E}_{1c}^{[2]}(v)-T_{1c}^{[1]}(u)P_{ca}\right).
		\end{equation*}
		By \eqref{eq:TTpcomm} and Lemma~\ref{lem:PQ}, we obtain the relation~\eqref{eq:ETp}. The relation ~\eqref{eq:FTp} can be verified similarly.
	\end{proof}
	
	\begin{prop}
		\label{prop:typeA}
		Suppose that the generator matrix $T(u)$ of $\X(\mathfrak{g}_{N})$ is factored as in \eqref{blockGaussDec} associated to the symmetric composition $(\nu_1,\nu_2,\ldots,\nu_{M})$ of $N$, and $m=\left[\frac{M}{2}\right]$. For $1\leqslant r,s \leqslant m$ and $1\leqslant a,b<m$, the following relations hold:
		\begin{align}
			&\left[\mathcal{D}_r^{[1]}(u),\mathcal{D}_s^{[2]}(v)\right]
			=\delta_{rs}\frac{P_{ss}}{u-v}
			\left(\mathcal{D}_s^{[1]}(u)\mathcal{D}_s^{[2]}(v)-\mathcal{D}_s^{[1]}(v)\mathcal{D}_s^{[2]}(u)\right),
			\label{eq:A:DaDb}\\
			&\begin{aligned}
				\left[\mathcal{D}_r^{[1]}(u),\mathcal{E}_a^{[2]}(v)\right]
				=&\delta_{ra}\mathcal{D}_a^{[1]}(u)\frac{P_{aa}}{u-v}\left(\mathcal{E}_a^{[2]}(v)-\mathcal{E}_a^{[2]}(u)\right)\\
				&-\delta_{r,a+1}
				\mathcal{D}_{a+1}^{[1]}(u)\left(\mathcal{E}_a^{[2]}(v)-\mathcal{E}_a^{[2]}(u)\right)\frac{P_{a+1,a+1}}{u-v},
			\end{aligned}
			\label{eq:A:DaEb}\\
			&\begin{aligned}
				\left[\mathcal{D}_r^{[1]}(u),\mathcal{F}_a^{[2]}(v)\right]
				=&\delta_{ra}
				\left(\mathcal{F}_a^{[2]}(u)-\mathcal{F}_a^{[2]}(v)\right)\frac{P_{aa}}{u-v}\mathcal{D}_a^{[1]}(u)\\
				&-\delta_{r,a+1}
				\frac{P_{a+1,a+1}}{u-v}\left(\mathcal{F}_a^{[2]}(u)-\mathcal{F}_a^{[2]}(v)\right)\mathcal{D}_{a+1}^{[1]}(u),
			\end{aligned}
			\label{eq:A:DaFb}\\
			&\left[\mathcal{E}_a^{[1]}(u),\mathcal{F}_b^{[2]}(v)\right]
			=\frac{\delta_{ab}}{u-v}\left(\widetilde{D}_a^{[1]}(u)P_{a,a+1}\mathcal{D}_{a+1}^{[1]}(u)
			-\mathcal{D}_{a+1}^{[2]}(v)P_{a,a+1}\widetilde{D}_a^{[2]}(v)\right),
			\label{eq:A:EaFb}\\
			&\left[\mathcal{E}_a^{[1]}(u),\mathcal{E}_b^{[2]}(v)\right]
			=\left[\mathcal{F}_a^{[1]}(u),\mathcal{F}_b^{[2]}(v)\right]=0,\text{ if }|a-b|>1,
			\label{eq:A:EaEb}\\
			&\left[\mathcal{E}_a^{[1]}(u),\mathcal{E}_a^{[2]}(v)\right]
			=\frac{P_{aa}}{u-v}\left(\mathcal{E}_a^{[1]}(u)-\mathcal{E}_a^{[1]}(v)\right)
			\left(\mathcal{E}_a^{[2]}(u)-\mathcal{E}_a^{[2]}(v)\right),
			\label{eq:A:EaEa}\\
			&u\left[\mathring{\mathcal{E}}_{a-1}^{[1]}(u),\mathcal{E}_a^{[2]}(v)\right]
			-v\left[\mathcal{E}_{a-1}^{[1]}(u),\mathring{\mathcal{E}}_a^{[2]}(v)\right]
			=\mathcal{E}_{a-1}^{[1]}(u)P_{aa}\mathcal{E}_a^{[2]}(v),
			\label{eq:A:Ea1Ea}\\
			&\left[\mathcal{F}_a^{[1]}(u),\mathcal{F}_a^{[2]}(v)\right]
			=-\left(\mathcal{F}_a^{[2]}(u)-\mathcal{F}_a^{[2]}(v)\right)\left(\mathcal{F}_a^{[1]}(u)-\mathcal{F}_a^{[1]}(v)\right)\frac{P_{aa}}{u-v},\label{eq:FaFa}\\
			&u\left[\mathring{\mathcal{F}}_{a-1}^{[1]}(u),\mathcal{F}_a^{[2]}(v)\right]
			-v\left[\mathcal{F}_{a-1}^{[1]}(u),\mathring{\mathcal{F}}_a^{[2]}(v)\right]
			=-\mathcal{F}_a^{[2]}(v)P_{aa}\mathcal{F}_{a-1}^{[1]}(u),
			\label{eq:A:Fa1Fa}\\
			&\left[\mathcal{E}_a^{[1]}(u),\left[\mathcal{E}_a^{[2]}(v),\mathcal{E}_b^{[3]}(w)\right]\right]
			+\left[\mathcal{E}_a^{[2]}(v),\left[\mathcal{E}_a^{[1]}(u),\mathcal{E}_b^{[3]}(w)\right]\right]=0, \text{ if }|a-b|=1,
			\label{eq:A:SerreEaEb}\\
			&\left[\mathcal{F}_a^{[1]}(u),\left[\mathcal{F}_a^{[2]}(v),\mathcal{F}_b^{[3]}(w)\right]\right]
			+\left[\mathcal{F}_a^{[2]}(v),\left[\mathcal{F}_a^{[1]}(u),\mathcal{F}_b^{[3]}(w)\right]\right]=0, \text{if}\,|a-b|=1,
			\label{eq:A:SerreFaFb}
		\end{align}
		where $\mathring{\mathcal{E}}_b(u)=\sum\limits_{r\geqslant2}\mathcal{E}_b^{(r)}\otimes u^{-r}$ and $\mathring{\mathcal{F}}_b(u)=\sum\limits_{r\geqslant2}\mathcal{F}_b^{(r)}\otimes u^{-r}$.
	\end{prop}
	
	\begin{proof}
		Since the subalgebra of $\X(\mathfrak{g}_N)$ generated by the coefficients of all entries in $T_{ab}(u)$ for $1\leqslant a,b\leqslant m$ is isomorphic to $\Y(\mathfrak{gl}_{\ell})$, where $\ell=\nu_1+\ldots+\nu_m$, all the relations \eqref{eq:A:DaDb}-\eqref{eq:A:SerreFaFb} have been given in the parabolic presentation of the Yangian $Y(\mathfrak{gl}_{\ell})$ associated with the composition $(\nu_1,\nu_2, \ldots,\nu_m)$. They have been verified in \cite{bk:pp}. We only present a proof of \eqref{eq:A:EaEa} to illustrate our method to conduct calculation with matrix blocks.
		
		We consider the block RTT relation \eqref{eq:blockRTT}:
		$$\left[T_{a,a+1}^{[1]}(u), T_{a,a+1}^{[2]}(v)\right]
		=\frac{P_{aa}}{u-v}\left(T_{a,a+1}^{[1]}(u)T_{a,a+1}^{[2]}(v)-T_{a,a+1}^{[1]}(v)T_{a,a+1}^{[2]}(u)\right).$$
		By applying the homomorphism $\varPsi_{a-1}$, it yields
		\begin{align*}
			&\varPsi_{a-1}\left(\left[T_{a,a+1}^{[1]}(u), T_{a,a+1}^{[2]}(v)\right]\right)\\
			&=\frac{P_{aa}}{u-v}\left(\mathcal{D}_a^{[1]}(u)\mathcal{E}_a^{[1]}(u)\mathcal{D}_a^{[2]}(v)\mathcal{E}_a^{[2]}(v)
			-\mathcal{D}_a^{[1]}(v)\mathcal{E}_a^{[1]}(v)\mathcal{D}_a^{[2]}(u)\mathcal{E}_a^{[2]}(u)\right),
		\end{align*}
		while we also have
		\begin{align*}
			&\varPsi_{a-1}\left(\left[T_{a,a+1}^{[1]}(u), T_{a,a+1}^{[2]}(v)\right]\right)
			=\mathcal{D}_a^{[2]}(v)\mathcal{D}_a^{[1]}(u)\left[\mathcal{E}_a^{[1]}(u),\mathcal{E}_a^{[2]}(v)\right]\\
			&\qquad+\mathcal{D}_a^{[2]}(v)\left[\mathcal{D}_a^{[1]}(u),\mathcal{E}_a^{[2]}(v)\right]\mathcal{E}_a^{[1]}(u)
			+\varPsi_{a-1}\left(\left[T_{a,a+1}^{[1]}(u),T_{aa}^{[2]}(v)\right]\right)\mathcal{E}_a^{[2]}(v),\\
			&\varPsi_{a-1}\left(\left[T_{a,a+1}^{[1]}(u),T_{aa}^{[2]}(v)\right]\right)
			=\frac{P_{aa}}{u-v}\left(\mathcal{D}_a^{[1]}(u)\mathcal{E}_a^{[1]}(u)\mathcal{D}_a^{[2]}(v)
			-\mathcal{D}_a^{[1]}(v)\mathcal{E}_a^{[1]}(v)\mathcal{D}_a^{[2]}(u)\right).
		\end{align*}
		Using \eqref{eq:A:DaEb}, we compute that
		\begin{align*}
			&\left(I-\frac{P_{aa}}{u-v}\right)\left[\mathcal{E}_a^{[1]}(u),\mathcal{E}_a^{[2]}(v)\right]
			+\frac{P_{aa}}{u-v}\left[\mathcal{E}_a^{[1]}(u),\mathcal{E}_a^{[2]}(u)\right]\\
			=&\frac{P_{aa}}{u-v}\left(I-\frac{P_{aa}}{u-v}\right)
			\left(\mathcal{E}_a^{[1]}(u)-\mathcal{E}_a^{[1]}(v)\right)\left(\mathcal{E}_a^{[2]}(u)-\mathcal{E}_a^{[2]}(v)\right).
		\end{align*}
		Now, $I-\frac{P_{aa}}{u-v}$ has the inverse $I+\frac{P_{aa}}{u-v-1}$. We thus obtain that
		\begin{align*}
			&\left[\mathcal{E}_a^{[1]}(u),\mathcal{E}_a^{[2]}(v)\right]
			+\frac{P_{aa}}{u-v}\left(I+\frac{P_{aa}}{u-v-1}\right)\left[\mathcal{E}_a^{[1]}(u),\mathcal{E}_a^{[2]}(u)\right]\\
			=&\frac{P_{aa}}{u-v}
			\left(\mathcal{E}_a^{[1]}(u)-\mathcal{E}_a^{[1]}(v)\right)\left(\mathcal{E}_a^{[2]}(u)-\mathcal{E}_a^{[2]}(v)\right).
		\end{align*}
		The term $\left[\mathcal{E}_a^{[1]}(u),\mathcal{E}_a^{[2]}(u)\right]$ indeed vanishes, which can be observed by multiplying $u-v-1$ on both sides of the equality and set $u=v+1$. Hence, the relation \eqref{eq:A:EaEa} follows.
	\end{proof}
	
	\begin{remark}
		The relation \eqref{eq:A:EaEa} is equivalent to
		$$u\left[\mathring{\mathcal{E}}_a^{[1]}(u),\mathcal{E}_a^{[2]}(v)\right]
		-v\left[\mathcal{E}_a^{[1]}(u),\mathring{\mathcal{E}}_a^{[2]}(v)\right]
		=-P_{aa}\left(\mathcal{E}_a^{[1]}(u)\mathcal{E}_a^{[2]}(v)+\mathcal{E}_a^{[1]}(v)\mathcal{E}_a^{[2]}(u)\right).$$
	\end{remark}
	
	\begin{prop}
		\label{prop:trans}
		Let $\mathcal{E}_{ab}(u)$ and $\mathcal{F}_{ba}(u)$ be the blocks obtained in the block Gauss fractorization \eqref{blockGaussDec} associated to the symmetric composition $(\nu_1,\nu_2, \ldots,\nu_{M})$ of $N$. Then, for each $a=1, 2, \ldots, \left[\frac{M-1}{2}\right]$,
		\begin{equation}
			\mathcal{E}_{(a+1)',a'}^t(u)=-\mathcal{E}_a(u+\kappa_{a+1}),\quad\text{ and }\quad
			\mathcal{F}_{a',(a+1)'}^t(u)=-\mathcal{F}_a(u+\kappa_{a+1}),\label{eq:trans}
		\end{equation}
		where $\kappa_{a+1}=\kappa-\sum\limits_{j=1}^{a}\nu_j$.
	\end{prop}
	
	\begin{proof}
		According to Theorem~\ref{embedding}, it suffices to consider the case where $a=1$. Since the proof to both identities are quite similar, we only show the first one here.
		
		Recall that $T(u)T^t(u+\kappa)=z(u)1$, which shows that the blocks in $T(u)$ and its inverse $\widetilde{T}(u)$ satisfy that
		$$T_{ab}^t(u+\kappa)= z(u)\widetilde{T}_{b^{\prime}a^{\prime}}(u),$$
		for $a, b=1,\ldots, M$. By the block Gauss decomposition, we have
		\begin{align*}
			\mathcal{D}_1^t(u+\kappa)=&T_{11}^t(u+\kappa)=z(u)\widetilde{T}_{1',1'}(u)=z(u)\widetilde{\mathcal{D}}_{1'}(u).\\
			\left(\mathcal{D}_1(u+\kappa)\mathcal{E}_1(u+\kappa)\right)^t
			=&T_{12}^t(u+\kappa)=\widetilde{T}_{2',1'}(u)=-z(u)\mathcal{E}_{2',1'}(u)\widetilde{\mathcal{D}}_{1'}(u).
		\end{align*}
		It follows that
		\begin{equation}
			\left(\mathcal{D}_1(u+\kappa)\mathcal{E}_1(u+\kappa)\right)^t
			=-\mathcal{E}_{2',1'}(u)\mathcal{D}_1^t(u+\kappa).
			\label{eq:transa}
		\end{equation}
		
		On the other hand, $a=1\leqslant\left[\frac{M-1}{2}\right]$, the relation \eqref{eq:A:DaEb} implies that
		$$\left[\mathcal{D}_1^{[1]}(u),\mathcal{E}_1^{[2]}(v)\right]=\mathcal{D}_1^{[1]}(u)\left(\mathcal{E}_1^{[1]}(v)-\mathcal{E}_1^{[1]}(u)\right)\frac{P_{21}}{u-v}.$$
		It yields, by applying the transpose $t$ to the second tensor factor, that
		\begin{equation}\label{eq:D1E1t}
			\left[\mathcal{D}_1^{[1]}(u),\left(\mathcal{E}_1^{[2]}(v)\right)^{t_2}\right]
			=\mathcal{D}_1^{[1]}(u)\left(\mathcal{E}_1^{[1]}(v)-\mathcal{E}_1^{[1]}(u)\right)\frac{Q_{21}}{u-v}
		\end{equation}
		Since $\left(\mathcal{E}_1^{[2]}(v)\right)^{t_2}Q_{11}=\mathcal{E}_1^{[1]}(v)Q_{21}$ by Lemma \ref{lem:APQ}, we multiply $Q_{11}$ from the right of both sides of \eqref{eq:D1E1t} and obtain that
		$$\frac{u-v-\nu_1}{u-v}\mathcal{D}_1^{[1]}(u)\mathcal{E}_1^{[1]}(v)Q_{21}
		=\left(\mathcal{E}_1^{[2]}(v)\right)^{t_2}\mathcal{D}_1^{[1]}(u)Q_{11}
		-\mathcal{D}_1^{[1]}(u)\mathcal{E}_1^{[1]}(u)\frac{\nu_1 Q_{21}}{u-v},
		$$
		which implies by setting $v=u-\nu_1$ that
		$$\mathcal{D}_1^{[1]}(u)\mathcal{E}_1^{[1]}(u)Q_{21}=\left(\mathcal{E}_1^{[2]}(u-\nu_1)\right)^{t_2}\mathcal{D}_1^{[1]}(u)Q_{11}.$$
		It follows from Lemma~\ref{lem:APQ} that
		$$\left(\mathcal{D}_1^{[2]}(u)\mathcal{E}_1^{[2]}(u)\right)^{t_2}Q_{11}=\left(\mathcal{E}_1^{[2]}(u-\nu_1)\right)^{t_2}\left(\mathcal{D}_1^{[2]}(u)\right)^{t_2}Q_{11},$$
		and hence,
		$$\left(\mathcal{D}_1(u)\mathcal{E}_1(u)\right)^t=\mathcal{E}_1^t(u-\nu_1)\mathcal{D}_1^t(u).$$
		Then \eqref{eq:transa} is simplified as $$-\mathcal{E}_{2',1'}(u)\mathcal{D}_1^t(u+\kappa)=\mathcal{E}_1^t(u+\kappa-\nu_1)\mathcal{D}_1^t(u+\kappa),$$
		which shows $\mathcal{E}_{2',1'}(u)=-\mathcal{E}_1^t(u+\kappa_2)$ since $\mathcal{D}_1(u+\kappa)$ is invertible.
	\end{proof}
	In the subsequent sections, we will figure out a complete set of relations among the generators obtained in Proposition~\ref{prop:generators}.

	\section{Parabolic relations associated with an odd composition}
	\label{sec:para:odd}
	
	This section is devoted to discuss the parabolic presentation of $\X(\mathfrak{g}_N)$ associated with an odd symmetric composition $\nu=(\nu_1,\ldots,\nu_m,\nu_{m+1},\nu_{m+2},\ldots,\nu_M)$, where $M=2m+1$. We will formulate the relations for types $B$ and $C$ simultaneously. The types $B$ and $C$ can be distinguished by various blocks $Q_{ab}$ in the $R$-matrix \eqref{eq:blockR}.
	
	\begin{lem}
		\label{lem:DaDm+1_odd}
		Let $M=2m+1\geqslant3$ be an odd integer and $\nu=(\nu_1,\nu_2,\ldots,\nu_M)$ be a symmetric composition of $N$. Then the following identities hold in $\X(\mathfrak{g}_{N})$:
		\begin{align}\label{eq:DaDm+1_odd}
			\left[\mathcal{D}_a^{[1]}(u),\mathcal{D}_{m+1}^{[2]}(v)\right]
			=&0,\text{ for }a\leqslant m\\\label{eq:Dm+1Dm+1_odd}
			\left[\mathcal{D}_{m+1}^{[1]}(u),\mathcal{D}_{m+1}^{[2]}(v)\right]
			=&\left(\frac{P_{m+1,m+1}}{u-v}-\frac{Q_{m+1,m+1}}{u-v-\kappa_{m+1}}\right)\mathcal{D}_{m+1}^{[1]}(u)\mathcal{D}_{m+1}^{[2]}(v)\\
			&-\mathcal{D}_{m+1}^{[2]}(v)\mathcal{D}_{m+1}^{[1]}(u)
			\left(\frac{P_{m+1,m+1}}{u-v}-\frac{Q_{m+1,m+1}}{u-v-\kappa_{m+1}}\right).\nonumber
		\end{align}
	\end{lem}
	
	\begin{proof}
		Equation \eqref{eq:DaDm+1_odd} can be obtained from \eqref{eq:TTpcomm}. Theorem~\ref{embedding} implies
		\eqref{eq:Dm+1Dm+1_odd}.
	\end{proof}
	
	\begin{lem}
		\label{lem:DEF_odd}
		Let $M=2m+1\geqslant3$ be an odd integer and $\nu=(\nu_1,\nu_2,\ldots,\nu_M)$ be a symmetric composition of $N$. Then the following identities hold in $\X(\mathfrak{g}_{N})$:
		\begin{align}
			\left[\mathcal{D}_a^{[1]}(u),\mathcal{E}_{m}^{[2]}(v)\right]
			=&0, \quad\left[\mathcal{D}_a^{[1]}(u),\mathcal{F}_{m}^{[2]}(v)\right]
			=0,\quad \text{for $a<m$},
			\label{eq:DiEm_odd}\\
			\left[\mathcal{D}_{m}^{[1]}(u),\mathcal{E}_m^{[2]}(v)\right]
			=&\frac{P_{mm}}{u-v}\mathcal{D}_m^{[2]}(u)\left(\mathcal{E}_m^{[2]}(v)-\mathcal{E}_{m}^{[2]}(u)\right),
			\label{eq:DmEm_odd}\\\label{eq:DmFm_odd}
			\left[\mathcal{D}_m^{[1]}(u),\mathcal{F}_m^{[2]}(v)\right]
			=&\frac{P_{m,m+1}}{u-v}\left(
			\mathcal{F}_m^{[1]}(u)-\mathcal{F}_m^{[1]}(v)\right)\mathcal{D}_m^{[1]}(u).
		\end{align}
	\end{lem}
	
	\begin{proof}
		The relations \eqref{eq:DiEm_odd} follows from \eqref{eq:TTpcomm} and \eqref{eq:A:DaDb}.
		
		For relation \eqref{eq:DmFm_odd}, we apply the homomorphism $\varPsi_{m-1}$ to the block RTT relation \eqref{eq:blockRTT} and obtain that
		\begin{align*}
			\varPsi_{m-1}\left(\left[T_{mm}^{[1]}(u), T_{m+1,m}^{[2]}(v)\right]\right)
			=\frac{P_{m,m+1}}{u-v}\varPsi_{m-1}\left(T_{m+1,m}^{[1]}(u)T_{mm}^{[2]}(v)-T_{m+1,m}^{[1]}(v)T_{mm}^{[2]}(u)\right).
		\end{align*}
		By the block Gauss decomposition, it can be rewritten as
		\begin{align*}
			\left[\mathcal{D}_m^{[1]}(u),\mathcal{F}_m^{[2]}(v) \mathcal{D}_{m}^{[2]}(v)\right]
			=
			\frac{P_{m,m+1}}{u-v}\left(\mathcal{F}_m^{[1]}(u)\mathcal{D}_m^{[1]}(u)\mathcal{D}_m^{[2]}(v)
			-\mathcal{F}_m^{[1]}(v)\mathcal{D}_m^{[1]}(v)\mathcal{D}_m^{[2]}(u)\right)
		\end{align*}
		Note that
		\begin{align*}
			\left[\mathcal{D}_m^{[1]}(u),\mathcal{F}_m^{[2]}(v) \mathcal{D}_{m}^{[2]}(v)\right]
			=\left[\mathcal{D}_m^{[1]}(u), \mathcal{F}_{m}^{[2]}(v)\right]\mathcal{D}_m^{[2]}(v)
			+\mathcal{F}_{m}^{[2]}(v)\left[\mathcal{D}_m^{[1]}(u), \mathcal{D}_m^{[2]}(v)\right].
		\end{align*}
		Then, we deduce from \eqref{eq:A:DaDb} and Lemma~\ref{lem:PQ} that
		\begin{align*}
			\left[\mathcal{D}_m^{[1]}(u),\mathcal{F}_m^{[2]}(v) \right]\mathcal{D}_{m}^{[2]}(v)
			=\frac{P_{m,m+1}}{u-v}\left(\mathcal{F}_m^{[1]}(u)-\mathcal{F}_m^{[1]}(v)	\right)\mathcal{D}_m^{[1]}(u)\mathcal{D}_m^{[2]}(v),
		\end{align*}
		which implies \eqref{eq:DmFm_odd} by the invertibility of
		$\mathcal{D}_m^{[2]}(v)$. The relation \eqref{eq:DmEm_odd} can be proved similarly.
	\end{proof}

	\begin{lem}
		\label{prop:EFD2_odd}
		Let $M=2m+1\geqslant3$ be an odd integer and $\nu=(\nu_1,\nu_2,\ldots,\nu_M)$ be a symmetric composition of $N$. Then the following identities hold in $\X(\mathfrak{g}_{N})$:
		\begin{align}
			\left[\mathcal{D}_{m+1}^{[1]}(u), \mathcal{E}_a^{[2]}(v)\right]
			=&\left[\mathcal{D}_{m+1}^{[1]}(u), \mathcal{F}_a^{[2]}(v)\right]
			=0,\text{ for }a<m,\label{eq:Dm+1Ea_odd}\\
			\left[\mathcal{E}_{m}^{[1]}(u),\mathcal{D}_{m+1}^{[2]}(v)\right]
			=&\mathcal{D}_{m+1}^{[2]}(v)\left(\mathcal{E}_{m}^{[1]}(v)-\mathcal{E}_{m}^{[1]}(u)\right)\frac{P_{m+1,m+1}}{u-v}\label{eq:EmDm+1_odd}\nonumber\\
			&-\mathcal{D}_{m+1}^{[2]}(v)\left(\mathcal{E}_{m}^{[1]}(v+\kappa_{m+1})-\mathcal{E}_{m}^{[1]}(u)\right)\frac{Q_{m+1,m+1}}{u-v-\kappa_{m+1}},\\
			\left[\mathcal{F}_{m}^{[1]}(u),\mathcal{D}_{m+1}^{[2]}(v)\right]
			=&\frac{P_{m+1,m+1}}{u-v}\left(\mathcal{F}_m^{[1]}(u)-\mathcal{F}_m^{[1]}(v)\right)\mathcal{D}_{m+1}^{[2]}(v)\nonumber\\
			&-\frac{Q_{m+1,m+1}}{u-v-\kappa_{m+1}}\left(\mathcal{F}_{m}^{[1]}(u)-\mathcal{F}_m^{[1]}(v+\kappa_{m+1})\right)\mathcal{D}_{m+1}^{[2]}(v),
			\label{eq:FmDm+1_odd}
		\end{align}
		where $\kappa_{m+1}$ has the same meaning as in Proposition~\ref{prop:trans}
	\end{lem}
	
	\begin{proof}
		The relation \eqref{eq:Dm+1Ea_odd} is implied by \eqref{eq:TTpcomm} and \eqref{eq:DaDm+1_odd}.
		
		For \eqref{eq:EmDm+1_odd}, it suffices to consider the case where $M=3$ by Theorem~\ref{embedding}. The block RTT relation \eqref{eq:blockRTT} in $\X(\mathfrak{g}_{2\nu_m+\nu_{m+1}})$ yields that
		\begin{align*}
			&\left[T_{12}^{[1]}(u),{T}_{22}^{[2]}(v)\right]
			=\frac{1}{u-v}\left(P_{12}T_{22}^{[1]}(u){T}_{12}^{[2]}(v)-T_{22}^{[2]}(v){T}_{12}^{[1]}(u)P_{22}\right)\\
			&\qquad+\frac{1}{u-v-\kappa_m}\left({T}_{23}^{[2]}(v)T_{11}^{[1]}(u)Q_{12}
			+T_{22}^{[2]}(v){T}_{12}^{[1]}(u)Q_{22}+T_{21}^{[2]}(v){T}_{13}^{[1]}(u)Q_{32}\right).
		\end{align*}
		Using the block Gauss decomposition and relation \eqref{eq:DaDm+1_odd}, its left hand side can be written as
		\begin{align*}
			&\mathcal{D}_1^{[1]}(u)\left[\mathcal{E}_1^{[1]}(u),\mathcal{D}_2^{[2]}(v)\right]
			+\left[T_{12}^{[1]}(u),T_{21}^{[2]}(v)\right]\mathcal{E}_1^{[2]}(v)\\
			&-\mathcal{F}_1^{[2]}(v)\left[T_{12}^{[1]}(u),T_{11}^{[2]}(v)\right]\mathcal{E}_1^{[2]}(v)
			+\mathcal{F}_1^{[2]}(v)\left[T_{12}^{[1]}(u),T_{12}^{[2]}(v)\right].
		\end{align*}
		Then it follows from \eqref{eq:blockRTT} and Lemma \ref{lem:APQ} that
		\begin{equation}
			\label{eq:odd:E1D2prf}
			\begin{aligned}
				&\mathcal{D}_1^{[1]}(u)\left[\mathcal{E}_1^{[1]}(u),\mathcal{D}_2^{[2]}(v)\right]
				=\frac{1}{u-v}\mathcal{D}_2^{[2]}(v) \mathcal{D}_1^{[1]}(u)P_{12}\left(\mathcal{E}_1^{[2]}(v)-\mathcal{E}_1^{[2]}(u)\right)\\
				&\qquad+\frac{1}{u-v-\kappa_m}\left(\mathcal{D}_2^{[2]}(v)\mathcal{E}_{23}^{[2]}(v)\mathcal{D}_1^{[1]}(u)Q_{12}+\mathcal{D}_1^{[1]}(u)\mathcal{D}_2^{[2]}(v)\mathcal{E}_1^{[1]}(u)Q_{22}\right).
			\end{aligned}
		\end{equation}
		
		Now, the relation \eqref{eq:A:DaEb} and Lemma~\ref{lem:PQ} imply that
		$$\mathcal{E}_1^{[2]}(v)\mathcal{D}_1^{[1]}(u)
		=\mathcal{D}_1^{[1]}(u)\mathcal{E}_1^{[2]}(v)
		-\mathcal{D}_1^{[1]}(u)\left(\mathcal{E}_1^{[1]}(v)-\mathcal{E}_1^{[1]}(u)\right)\frac{P_{21}}{u-v}.$$
		Since $\mathcal{E}_1^t(v)=-\mathcal{E}_{23}(v-\kappa_{m+1})$ by Proposition~\ref{prop:trans} and $P_{21}^{t_2}=Q_{21}$ by Lemma~\ref{lem:PQtrs}, we have
		$$\mathcal{E}_{23}^{[2]}(v)\mathcal{D}_1^{[1]}(u)
		=-\mathcal{D}_1^{[1]}(u)\mathcal{E}_1^{t,[2]}(v+\kappa_{m+1})
		+\mathcal{D}_1^{[1]}(u)\left(\mathcal{E}_1^{[1]}(v+\kappa_{m+1})-\mathcal{E}_1^{[1]}(u)\right)\frac{Q_{21}}{u-v-\kappa_{m+1}}.$$
		By equations~\eqref{eq:PPQQ} and \eqref{eq:A1Q}, it further implies that
		\begin{align*}
			\mathcal{E}_{23}^{[2]}(v)\mathcal{D}_1^{[1]}(u)Q_{12}
			=&-\mathcal{D}_1^{[1]}(u)\mathcal{E}_1^{[1]}(v+\kappa_{m+1})Q_{22}\\
			&+\mathcal{D}_1^{[1]}(u)\left(\mathcal{E}_1^{[1]}(v+\kappa_{m+1})-\mathcal{E}_1^{[1]}(u)\right)\frac{\nu_{m}Q_{22}}{u-v-\kappa_{m+1}}.
		\end{align*}
		
		Consequently, the relation \eqref{eq:odd:E1D2prf} can be simplified as:
		\begin{align*}
			\mathcal{D}_1^{[1]}(u)\left[\mathcal{E}_1^{[1]}(u),\mathcal{D}_2^{[2]}(v)\right]
			=&\mathcal{D}_1^{[1]}(u)\mathcal{D}_2^{[2]}(v) \left(\mathcal{E}_1^{[1]}(v)-\mathcal{E}_1^{[1]}(u)\right)\frac{P_{22}}{u-v}\\
			&-\mathcal{D}_1^{[1]}(u)\mathcal{D}_2^{[2]}(v)\left(\mathcal{E}_1^{[1]}(v+\kappa_{m+1})-\mathcal{E}_1^{[1]}(u)\right)\frac{Q_{22}}{u-v-\kappa_{m+1}}.
		\end{align*}
		Since $\mathcal{D}_1(u)$ is invertible, we obtain that
		\begin{align*}
			\left[\mathcal{E}_1^{[1]}(u),\mathcal{D}_2^{[2]}(v)\right]
			=&\mathcal{D}_2^{[2]}(v) \frac{\mathcal{E}_1^{[1]}(v)-\mathcal{E}_1^{[1]}(u)}{u-v}P_{22}
			-\mathcal{D}_2^{[2]}(v)\frac{\mathcal{E}_1^{[1]}(v+\kappa_{m+1})-\mathcal{E}_1^{[1]}(u)}{u-v-\kappa_{m+1}}Q_{22}.
		\end{align*}
		We complete verifying \eqref{eq:EmDm+1_odd}. The relation \eqref{eq:FmDm+1_odd} can be checked similarly.		
	\end{proof}
	
	\begin{lem}
		\label{prop:EF_odd}
		Let $M=2m+1\geqslant3$ be an odd integer and $\nu=(\nu_1,\nu_2,\ldots,\nu_M)$ be a symmetric composition of $N$. Then the following identities hold in $\X(\mathfrak{g}_{N})$:
		\begin{align}
			\left[\mathcal{E}_m^{[1]}(u),\mathcal{F}_a^{[2]}(v)\right]
			=&\left[\mathcal{E}_a^{[1]}(u),\mathcal{F}_m^{[2]}(v)\right]
			=0,~~\text{for  }a<m,
			\label{eq:FjEmFmEj_odd}\\		
			\label{eq:E1F1odd}
			\left[\mathcal{E}_m^{[1]}(u),\mathcal{F}_m^{[2]}(v)\right]
			=&\frac{1}{u-v}\left(\widetilde{\mathcal{D}}_m^{[1]}(u)P_{m,m+1}\mathcal{D}_{m+1}^{[1]}(u)
			-\mathcal{D}_{m+1}^{[2]}(v)P_{m,m+1}\widetilde{\mathcal{D}}_m^{[2]}(v)\right)
		\end{align}
	\end{lem}
	
	\begin{proof}
		The relation \eqref{eq:FjEmFmEj_odd} can be obtained by observing
		Lemma \ref{lem:EFTp}, the relations \eqref{eq:TTpcomm}, \eqref{eq:A:DaEb},
		and \eqref{eq:A:DaFb}.
		
		For \eqref{eq:E1F1odd}, we only need to check the case where $M=3$ due to the embedding theorem~\ref{embedding}. In this situation, the block RTT relation \eqref{eq:blockRTT} yields that
		\begin{align*}
			\left[T_{12}^{[1]}(u),T_{21}^{[2]}(v)\right]
			=\frac{P_{12}}{u-v}T_{22}^{[1]}(u)T_{11}^{[2]}(v)
			-T_{22}^{[2]}(v)T_{11}^{[1]}(u)\frac{P_{12}}{u-v}.
		\end{align*}
		It can be rewritten via the block Gauss decomposition \eqref{blockGaussDec} as following:
		\begin{align*}
			\left[\mathcal{D}_1^{[1]}(u)\mathcal{E}_1^{[1]}(u),\mathcal{F}_1^{[2]}(v)\mathcal{D}_1^{[2]}(v)\right]
			=&\frac{P_{12}}{u-v}\left(\mathcal{D}_2^{[1]}(u)\mathcal{D}_1^{[2]}(v)+\mathcal{F}_1^{[1]}(u)\mathcal{D}_1^{[1]}(u)\mathcal{E}_1^{[1]}(u)\mathcal{D}_1^{[2]}(v)\right)
			\\
			&-\left(\mathcal{D}_2^{[2]}(v)\mathcal{D}_1^{[1]}(u)+\mathcal{F}_1^{[2]}(v)\mathcal{D}_1^{[2]}(v)\mathcal{E}_1^{[2]}(v)\mathcal{D}_1^{[1]}(u)\right)\frac{P_{12}}{u-v},
		\end{align*}
		whose left hand side equals to
		\begin{equation*}
			\mathcal{D}_1^{[1]}(u)\left[\mathcal{E}_1^{[1]}(u),\mathcal{F}_1^{[2]}(v)\right]\mathcal{D}_1^{[2]}(v)+\left[\mathcal{D}_1^{[1]}(u),\mathcal{F}_1^{[2]}(v)\right]\mathcal{E}_1^{[1]}(u)\mathcal{D}_1^{[2]}(v)+
			\mathcal{F}_1^{[2]}(v)\left[T_{12}^{[1]}(u),T_{11}^{[2]}(v)\right].
		\end{equation*}
		Note that the relation \eqref{eq:blockRTT} also implies that
		$$\left[T_{12}^{[1]}(u), T_{11}^{[2]}(v)\right]
		=\frac{1}{u-v}\left(
		P_{11}\mathcal{D}_{1}^{[1]}(u)\mathcal{E}_{1}^{[1]}(u) \mathcal{D}_{1}^{[2]}(v)-
		\mathcal{D}_{1}^{[2]}(v)\mathcal{E}_{1}^{[2]}(v)\mathcal{D}_{1}^{[1]}(u)P_{12}\right).$$
		By the relations \eqref{eq:DaDm+1_odd}, \eqref{eq:DmFm_odd},  and Lemma~\ref{lem:APQ}, we further deduce that
		$$\mathcal{D}_1^{[1]}(u)\left[\mathcal{E}_1^{[1]}(u),\mathcal{F}_1^{[2]}(v)\right]\mathcal{D}_1^{[2]}(v)
		=\frac{1}{u-v}\left(P_{12}\mathcal{D}_2^{[1]}(u)\mathcal{D}_1^{[2]}(v)
		-\mathcal{D}_1^{[1]}(u)\mathcal{D}_2^{[2]}(v)P_{12}\right),$$
		which implies \eqref{eq:E1F1odd} since $\mathcal{D}_1(u)$ and $\mathcal{D}_2(v)$ are invertible.
	\end{proof}

	\begin{lem}
		\label{prop:sp2:EEFF}
		Let $M=2m+1\geqslant3$ be an odd integer and $\nu=(\nu_1,\nu_2,\ldots,\nu_M)$ be a symmetric composition of $N$. Then the following identities hold in $\X(\mathfrak{g}_{N})$:
		\begin{align}
			&\begin{aligned}
				&\left[\mathcal{E}_m^{[1]}(u),\mathcal{E}_m^{[2]}(v)\right]+\left[\mathcal{E}_m^{[2]}(u),\mathcal{E}_m^{[1]}(v)\right]\\
				&=\left\{\mathcal{E}_{m}^{[1]}(v)-\mathcal{E}_{m}^{[1]}(u),\mathcal{E}_{m}^{[2]}(v)-\mathcal{E}_{m}^{[2]}(u)\right\}\frac{P_{m+1,m+1}-\frac{1}{2}Q_{m+1,m+1}}{u-v},
			\end{aligned}\label{eq:odd:EmEm}\\
			&\begin{aligned}
				&\left[\mathcal{F}_m^{[1]}(u),\mathcal{F}_m^{[2]}(v)\right]+\left[\mathcal{F}_m^{[2]}(u),\mathcal{F}_m^{[1]}(v)\right]\\
				&=\frac{P_{m+1,m+1}-\frac{1}{2}Q_{m+1,m+1}}{u-v}\left\{\mathcal{F}_{m}^{[1]}(u)-\mathcal{F}_{m}^{[1]}(v),\mathcal{F}_{m}^{[2]}(u)-\mathcal{F}_{m}^{[2]}(v)\right\},\label{eq:odd:FmFm}
			\end{aligned}
		\end{align}
		where $\{A,B\}=AB+BA$.
	\end{lem}
	\begin{proof}
		According to the embedding theorem~\ref{embedding}, we only need to verify the case where $M=3$. In this situation, it is known from the block Gauss decomposition~\ref{blockGaussDec} that
		$$T_{11}(u)=\mathcal{D}_1(u), \quad T_{12}(u)=\mathcal{D}_1(u)\mathcal{E}_1(u),\quad  T_{13}(u)=\mathcal{D}_1(u)\mathcal{E}_{13}(u).$$
		
		We consider the block RTT relation \eqref{eq:blockRTT}:
		\begin{align*}
			&\left[T_{12}^{[1]}(u),T_{12}^{[2]}(v)\right]
			=\frac{1}{u-v}\left(P_{11}T_{12}^{[1]}(u)T_{12}^{[2]}(v)-T_{12}^{[2]}(v)T_{12}^{[1]}(u)P_{22}\right)\\
			&\qquad+\frac{1}{u-v-\kappa}\left(T_{13}^{[2]}(v)T_{11}^{[1]}(u)Q_{12}+T_{12}^{[2]}(v)T_{12}^{[1]}(u)Q_{22}+T_{11}^{[2]}(v)T_{13}^{[1]}(u)Q_{32}\right),
		\end{align*}
		whose left hand side can be written as
		$$\left[T_{12}^{[1]}(u), T_{11}^{[2]}(v)\right]\mathcal{E}_1^{[2]}(v)
		+\mathcal{D}_1^{[2]}(v)\mathcal{D}_1^{[1]}(u)\left[\mathcal{E}_1^{[1]}(u),\mathcal{E}_1^{[2]}(v)\right]
		+\mathcal{D}_1^{[2]}(v)\left[\mathcal{D}_1^{[1]}(u),\mathcal{E}_1^{[2]}(v)\right]\mathcal{E}_1^{[1]}(u).$$
		Then it follows from \eqref{eq:blockRTT} and \eqref{eq:DmEm_odd} that
		\begin{equation}
			\label{eq:odd:E1E1a}
			\begin{aligned}
				&\mathcal{D}_1^{[2]}(v)\mathcal{D}_1^{[1]}(u)\left[\mathcal{E}_1^{[1]}(u),\mathcal{E}_1^{[2]}(v)\right]\\
				=&\mathcal{D}_1^{[2]}(v)\mathcal{E}_1^{[2]}(v)\mathcal{D}_1^{[1]}(u)
				\left(\mathcal{E}_1^{[1]}(v)-\mathcal{E}_1^{[1]}(u)\right)\frac{P_{22}}{u-v}\\
				&-\mathcal{D}_1^{[2]}(v)\mathcal{D}_1^{[1]}(u)
				\left(\mathcal{E}_1^{[1]}(v)-\mathcal{E}_1^{[1]}(u)\right)\mathcal{E}_1^{[2]}(u)\frac{P_{22}}{u-v}\\
				&+\frac{\mathcal{D}_1^{[2]}(v)}{u-v-\kappa}\left(\mathcal{E}_{13}^{[2]}(v)\mathcal{D}_1^{[1]}(u)Q_{12}
				+\mathcal{E}_1^{[2]}(v)\mathcal{D}_1^{[1]}(u)\mathcal{E}_1^{[1]}(u)Q_{22}
				+\mathcal{D}_1^{[1]}(u)\mathcal{E}_{13}^{[1]}(u)Q_{32}\right).
			\end{aligned}
		\end{equation}
		
		Now, we use the block RTT relation \eqref{eq:blockRTT} for $\left[T_{11}^{[1]}(u), T_{13}^{[2]}(v)\right]$ to deduce that
		\begin{align*}
			& {\left[\mathcal{D}_1^{[1]}(u), \mathcal{E}_{13}^{[2]}(v)\right]=\frac{1}{u-v} \mathcal{D}_1^{[1]}(u)\left(\mathcal{E}_{13}^{[1]}(v)-\mathcal{E}_{13}^{[1]}(u)\right) P_{31}} \\
			& +\frac{1}{u-v-\kappa}\left(\mathcal{E}_{13}^{[2]}(v)\mathcal{D}_1^{[1]}(u)Q_{11}
			+\mathcal{E}_1^{[2]}(v)\mathcal{D}_1^{[1]}(u)\mathcal{E}_1^{[1]}(u)Q_{21}
			+\mathcal{D}_1^{[1]}(u)\mathcal{E}_{13}^{[1]}(u)Q_{31}\right).
		\end{align*}
		Right multiplying $Q_{12}$ on both sides of the above equality, we obtain that
		\begin{align*}
			\frac{u-v-\kappa_2}{u-v-\kappa}\mathcal{E}_{13}^{[2]}(v)\mathcal{D}_1^{[1]}(u)Q_{12}
			=&\mathcal{D}_1^{[1]}(u)\left(\mathcal{E}_{13}^{[2]}(v)
			+\frac{\mathcal{E}_{13}^{[1]}(u)-\mathcal{E}_{13}^{[1]}(v)}{u-v}\right)Q_{12}\\
			&-\frac{\nu_1}{u-v-\kappa}\left(\mathcal{D}_1^{[1]}(u)\mathcal{E}_{13}^{[1]}(u)Q_{32}+\mathcal{E}_1^{[2]}(v)\mathcal{D}_1^{[1]}(u)\mathcal{E}_1^{[1]}(u)Q_{22}\right),
		\end{align*}
		and thus
		\begin{align*}
			&\frac{1}{u-v-\kappa}\left(\mathcal{E}_{13}^{[2]}(v)\mathcal{D}_1^{[1]}(u)Q_{12}
			+\mathcal{E}_1^{[2]}(v)\mathcal{D}_1^{[1]}(u)\mathcal{E}_1^{[1]}(u)Q_{22}
			+\mathcal{D}_1^{[1]}(u)\mathcal{E}_{13}^{[1]}(u)Q_{32}\right)\\
			=&\mathcal{D}_1^{[1]}(u)
			\left(\mathcal{E}_{13}^{[2]}(v)\pm\mathcal{E}_{13}^{[1]}(u)
			+\frac{\mathcal{E}_{13}^{[1]}(u)-\mathcal{E}_{13}^{[1]}(v)}{u-v}\right)\frac{Q_{12}}{u-v-\kappa_2}\\
			&+\mathcal{E}_1^{[2]}(v)\mathcal{D}_1^{[1]}(u)\mathcal{E}_1^{[1]}(u)\frac{Q_{22}}{u-v-\kappa_2},
		\end{align*}
		Throughout the proof to this lemma, if the sign $\pm$ or $\mp$ appears, we always choose the upper one if $\mathfrak{g}_N$ is of type $B$ and choose the lower one if $\mathfrak{g}_N$ is of type $C$. We also deduce from the relation \eqref{eq:DmEm_odd} that
		$$\mathcal{E}_1^{[2]}(v)\mathcal{D}_1^{[1]}(u)=\mathcal{D}_1^{[1]}(u)\left(\mathcal{E}_1^{[2]}(v)-\frac{\mathcal{E}_1^{[1]}(v)-\mathcal{E}_1^{[1]}(u)}{u-v}P_{21}\right).$$
		Hence, the relation \eqref{eq:odd:E1E1a} can be further simplified as:
		\begin{align*}
			&\left[\mathcal{E}_1^{[1]}(u),\mathcal{E}_1^{[2]}(v)\right]\left(1-\frac{P_{22}}{u-v}\right)
			+\left[\mathcal{E}_1^{[1]}(v),\mathcal{E}_1^{[2]}(v)\right]\frac{P_{22}}{u-v}\\
			=&\left(\mathcal{E}_1^{[1]}(v)-\mathcal{E}_1^{[1]}(u)\right)\left(\mathcal{E}_1^{[2]}(v)-\mathcal{E}_1^{[2]}(u)\right)
			\frac{P_{22}}{u-v}\left(1-\frac{P_{22}}{u-v}\right)\\
			&+\left(\mathcal{E}_1^{[2]}(v)\mathcal{E}_1^{[1]}(u)
			\mp\frac{1}{u-v}\left(\mathcal{E}_1^{[1]}(v)-\mathcal{E}_1^{[1]}(u)\right)\mathcal{E}_1^{[2]}(u)\right)
			\frac{Q_{22}}{u-v-\kappa_2}\\
			&+\left(\mathcal{E}_{13}^{[2]}(v)\pm\mathcal{E}_{13}^{[1]}(u)
			+\frac{\mathcal{E}_{13}^{[1]}(u)-\mathcal{E}_{13}^{[1]}(v)}{u-v}\right)\frac{Q_{12}}{u-v-\kappa_2}.
		\end{align*}
		Right multiplying $1+\frac{P_{22}}{u-v}$ on both sides of the above equation and using Lemma~\ref{lem:PQ}, we compute that
		\begin{equation}
			\label{eq:odd:E1E1b}		
			\begin{aligned}
				&\frac{(u-v)^2-1}{(u-v)^2}\left[\mathcal{E}_1^{[1]}(u),\mathcal{E}_1^{[2]}(v)\right]
				+\left[\mathcal{E}_1^{[1]}(v),\mathcal{E}_1^{[2]}(v)\right]\frac{P_{22}}{u-v}\left(1+\frac{P_{22}}{u-v}\right)\\
				=&\frac{(u-v)^2-1}{(u-v)^2}\left(\mathcal{E}_1^{[1]}(v)-\mathcal{E}_1^{[1]}(u)\right)\left(\mathcal{E}_1^{[2]}(v)-\mathcal{E}_1^{[2]}(u)\right)\frac{P_{22}}{u-v}\\
				&+\frac{u-v\pm1}{u-v}\left(\mathcal{E}_1^{[2]}(v)\mathcal{E}_1^{[1]}(u)
				\mp\frac{1}{u-v}\left(\mathcal{E}_1^{[1]}(v)-\mathcal{E}_1^{[1]}(u)\right)\mathcal{E}_1^{[2]}(u)\right)
				\frac{Q_{22}}{u-v-\kappa_2}\\
				&+\frac{u-v\pm1}{u-v}\left(\mathcal{E}_{13}^{[2]}(v)\pm\mathcal{E}_{13}^{[1]}(u)
				+\frac{\mathcal{E}_{13}^{[1]}(u)-\mathcal{E}_{13}^{[1]}(v)}{u-v}\right)\frac{Q_{12}}{u-v-\kappa_2},
			\end{aligned}
		\end{equation}
		Let $u=v\mp 1$, we have
		\begin{equation}
			\label{eq:odd:E1E1c}
			\left[\mathcal{E}_1^{[1]}(v),\mathcal{E}_1^{[2]}(v)\right]\left(1\mp P_{22}\right)=\left(1\pm P_{11}\right)\left[\mathcal{E}_1^{[1]}(v),\mathcal{E}_1^{[2]}(v)\right]
			=0.
		\end{equation}
		Then the relation \eqref{eq:odd:E1E1b} implies that
		\begin{equation}
			\label{eq:odd:E1E1d}		
			\begin{aligned}
				&\frac{u-v\mp 1}{u-v}\left[\mathcal{E}_1^{[1]}(u),\mathcal{E}_1^{[2]}(v)\right]
				\pm\frac{1}{u-v}\left[\mathcal{E}_1^{[1]}(v),\mathcal{E}_1^{[2]}(v)\right]\\
				=&\frac{u-v\mp1}{u-v}\left(\mathcal{E}_1^{[1]}(v)-\mathcal{E}_1^{[1]}(u)\right)\left(\mathcal{E}_1^{[2]}(v)-\mathcal{E}_1^{[2]}(u)\right)\frac{P_{22}}{u-v}\\
				&+\left(\mathcal{E}_1^{[2]}(v)\mathcal{E}_1^{[1]}(u)
				\mp\frac{1}{u-v}\left(\mathcal{E}_1^{[1]}(v)-\mathcal{E}_1^{[1]}(u)\right)\mathcal{E}_1^{[2]}(u)\right)
				\frac{Q_{22}}{u-v-\kappa_2}\\
				&+\left(\mathcal{E}_{13}^{[2]}(v)\pm\mathcal{E}_{13}^{[1]}(u)
				+\frac{\mathcal{E}_{13}^{[1]}(u)-\mathcal{E}_{13}^{[1]}(v)}{u-v}\right)\frac{Q_{12}}{u-v-\kappa_2}.
			\end{aligned}
		\end{equation}
		Since $Q_{22}(1\mp P_{22})=0$ and $Q_{12}(1\mp P_{22})=0$ by Lemma~\ref{lem:PQ}, we deduce from \eqref{eq:odd:E1E1c} and \eqref{eq:odd:E1E1d} that
		\begin{equation}
			\label{eq:odd:E1E1e}
			\left[\mathcal{E}_1^{[1]}(u),\mathcal{E}_1^{[2]}(v)\right]\left(1\mp P_{22}\right)
			=\left(\mathcal{E}_1^{[1]}(v)-\mathcal{E}_1^{[1]}(u)\right)\left(\mathcal{E}_1^{[2]}(v)-\mathcal{E}_1^{[2]}(u)\right)
			\frac{P_{22}\mp 1}{u-v}.
		\end{equation}
		
		On the other hand, by left multiplying $1\pm P_{11}$ on both sides of \eqref{eq:odd:E1E1d}, we have
		\begin{equation}
			\label{eq:odd:E1E1f}		
			\begin{aligned}
				&\frac{u-v\mp1}{u-v}(1\pm P_{11})\left[\mathcal{E}_1^{[1]}(u),\mathcal{E}_1^{[2]}(v)\right]\\
				=&\frac{u-v\mp 1}{u-v}(1\pm P_{11}) \left(\mathcal{E}_1^{[1]}(v)-\mathcal{E}_1^{[1]}(u)\right)\left(\mathcal{E}_1^{[2]}(v)-\mathcal{E}_1^{[2]}(u)\right)
				\frac{P_{22}}{u-v}\\
				&+\frac{u-v\mp1}{u-v}\left(\mathcal{E}_1^{[1]}(v)\mathcal{E}_1^{[2]}(u)+\mathcal{E}_1^{[2]}(v)\mathcal{E}_1^{[1]}(u)\right)\frac{Q_{22}}{u-v-\kappa_2}\\
				&\pm(1\pm P_{11})\mathcal{E}_1^{[1]}(u)\mathcal{E}_1^{[2]}(u)\frac{Q_{22}}{(u-v)(u-v-\kappa_2)}\\
				&+\left(\frac{u-v\mp1}{u-v}\left(\mathcal{E}_{13}^{[2]}(v)\pm\mathcal{E}_{13}^{[1]}(v)\right)
				+\frac{u-v\pm1}{u-v}\left(\mathcal{E}_{13}^{[2]}(u)\pm\mathcal{E}_{13}^{[1]}(u)\right)\right)\frac{Q_{12}}{u-v-\kappa_2}.\end{aligned}
		\end{equation}
		Let $v=u\mp1$, we obtain that
		$$\left(1\pm P_{11}\right)
		\mathcal{E}_1^{[1]}(u)\mathcal{E}_1^{[2]}(u)Q_{22}
		+2\left(\mathcal{E}_{13}^{[2]}(u)\pm\mathcal{E}_{13}^{[1]}(u)\right)Q_{12}=0.$$
		Therefore, the relation \eqref{eq:odd:E1E1f} can be written as
		\begin{align*}
			&(1\pm P_{11})\left[\mathcal{E}_1^{[1]}(u),\mathcal{E}_1^{[2]}(v)\right]\\
			=&(1\pm P_{11}) \left(\mathcal{E}_1^{[1]}(v)-\mathcal{E}_1^{[1]}(u)\right)
			\left(\mathcal{E}_1^{[2]}(v)-\mathcal{E}_1^{[2]}(u)\right)\frac{P_{22}}{u-v}\\
			&-(1\pm P_{11})\left(\mathcal{E}_1^{[1]}(v)\mathcal{E}_1^{[2]}(v)+\mathcal{E}_1^{[1]}(u)\mathcal{E}_1^{[2]}(u)-2\mathcal{E}_1^{[1]}(v)\mathcal{E}_1^{[2]}(u)\right)\frac{Q_{22}}{2(u-v-\kappa_2)}.
		\end{align*}
		which is equivalent to
		\begin{equation}\label{eq:odd:E1E1g}
			\begin{aligned}
				&\left(1\pm P_{11}\right)\left[\mathcal{E}_1^{[1]}(u),\mathcal{E}_1^{[2]}(v)\right]\left(1+\frac{Q_{22}}{2(u-v-\kappa_2)}\right)\\
				=&(1\pm P_{11}) \left(\mathcal{E}_1^{[1]}(v)-\mathcal{E}_1^{[1]}(u)\right)
				\left(\mathcal{E}_1^{[2]}(v)-\mathcal{E}_1^{[2]}(u)\right)\frac{P_{22}}{u-v}\\
				&-\left(1\pm P_{11}\right)\left(\mathcal{E}_1^{[1]}(v)-\mathcal{E}_1^{[1]}(u)\right)
				\left(\mathcal{E}_1^{[2]}(v)-\mathcal{E}_1^{[2]}(u)\right)\frac{Q_{22}}{2(u-v-\kappa_2)}.
			\end{aligned}
		\end{equation}
		Right multiplying $Q_{22}$, we have
		\begin{align*}
			&\left(1\pm P_{11}\right)\left[\mathcal{E}_1^{[1]}(u),\mathcal{E}_1^{[2]}(v)\right]Q_{22}\left(1+\frac{\nu_2}{2(u-v-\kappa_2)}\right)\\
			=&\left(1\pm P_{11}\right)\left(\mathcal{E}_1^{[1]}(v)-\mathcal{E}_1^{[1]}(u)\right)
			\left(\mathcal{E}_1^{[2]}(v)-\mathcal{E}_1^{[2]}(u)\right)Q_{22}
			\left(\pm\frac{1}{u-v}-\frac{\nu_2}{2(u-v-\kappa_2)}\right).
		\end{align*}
		which yields that
		$$\left(1\pm P_{11}\right)\left[\mathcal{E}_1^{[1]}(u),\mathcal{E}_1^{[2]}(v)\right]Q_{22}\\
		=-\left(1\pm P_{11}\right)\left(\mathcal{E}_1^{[1]}(v)-\mathcal{E}_1^{[1]}(u)\right)
		\left(\mathcal{E}_1^{[2]}(v)-\mathcal{E}_1^{[2]}(u)\right)Q_{22}\frac{\kappa_2}{u-v}.$$
		Thus, the equation \eqref{eq:odd:E1E1g} is further reduced to
		\begin{equation}\label{eq:odd:E1E1h}
			\begin{aligned}
				&\left(1\pm P_{11}\right)\left[\mathcal{E}_1^{[1]}(u),\mathcal{E}_1^{[2]}(v)\right]\\
				=&\left(1\pm P_{11}\right)\left(\mathcal{E}_1^{[1]}(v)-\mathcal{E}_1^{[1]}(u)\right)
				\left(\mathcal{E}_1^{[2]}(v)-\mathcal{E}_1^{[2]}(u)\right)\frac{P_{22}-\frac{1}{2}Q_{22}}{u-v}.
			\end{aligned}
		\end{equation}
		
		Considering the difference of \eqref{eq:odd:E1E1e} and \eqref{eq:odd:E1E1h}, we have
		$$P_{11}\left[\mathcal{E}_1^{[1]}(u),\mathcal{E}_1^{[2]}(v)\right]+\left[\mathcal{E}_1^{[1]}(u),\mathcal{E}_1^{[2]}(v)\right]P_{22}
		=\left\{\mathcal{E}_1^{[1]}(v)-\mathcal{E}_1^{[1]}(u),\mathcal{E}_1^{[2]}(v)-\mathcal{E}_1^{[2]}(u)\right\}\frac{1\mp\frac{1}{2}Q_{22}}{u-v},$$
		which yields \eqref{eq:odd:EmEm} by Lemma \ref{lem:APQ} and Lemma \ref{lem:PQ}. The identity \eqref{eq:odd:FmFm} can be verified similarly.
	\end{proof}
	
	\begin{lem}
		\label{lem:odd:Em1Em}
		Let $M=2m+1\geqslant5$ be an odd integer and $\nu=(\nu_1,\nu_2,\ldots,\nu_M)$ be a symmetric composition of $N$. Then the following identities hold in $\X(\mathfrak{g}_{N})$:
		\begin{align}
			&\left[\mathcal{E}_l^{[1]}(u),\mathcal{E}_{m}^{[2]}(v)\right]=0,~~\left[\mathcal{F}_l^{[1]}(u),\mathcal{F}_{m}^{[2]}(v)\right]=0\label{eq:ElEm},
			~~\text{for $l<m-1$},\\
			&u\left[\mathring{\mathcal{E}}_{m-1}^{[1]}(u),\mathcal{E}_m^{[2]}(v)\right]
			-v\left[\mathcal{E}_{m-1}^{[1]}(u),\mathring{\mathcal{E}}_m^{[2]}(v)\right]
			=\mathcal{E}_{m-1}^{[1]}(u)\mathcal{E}_m^{[1]}(v)P_{m+1,m},\label{eq:Em-1Em_odd}\\
			&u\left[\mathring{\mathcal{F}}_{m-1}^{[1]}(u),\mathcal{F}_m^{[2]}(v)\right]
			-v\left[\mathcal{F}_{m-1}^{[1]}(u),\mathring{\mathcal{F}}_m^{[2]}(v)\right]
			=-\mathcal{F}_{m}^{[2]}(v)\mathcal{F}_{m-1}^{[2]}(u)P_{m,m-1}.\label{eq:Fm-1Fm_odd}
		\end{align}
		where $\mathring{\mathcal{E}}_a(u)=\sum\limits_{r\geqslant2}\mathcal{E}_a^{(r)} u^{-r}$ for $a=m-1,m$.
	\end{lem}
	\begin{proof}
		The identities \eqref{eq:ElEm} follow from \eqref{eq:TTpcomm}.
		
		In order to verify \eqref{eq:Em-1Em_odd}, we deduce from Lemma \ref{lem:EFTp} that
		$$\left[\mathcal{E}_{m-1}^{[1]}(u), \varPsi_{m-1}\left({T}_{m,m+1}^{[2]}(v)\right)\right]
		=\varPsi_{m-1}\left(T_{mm}^{[2]}(v)\right)\left(\mathcal{E}_{m-1,m+1}^{[1]}(v)-\mathcal{E}_{m-1,m+1}^{[1]}(u)\right)\frac{P_{m+1,m}}{u-v}.$$
		Since $\varPsi_{m-1}\left(T_{mm}(v)\right)=\mathcal{D}_m(v)$ and $\varPsi_{m-1}\left(T_{m,m+1}(v)\right)=\mathcal{D}_m(v)\mathcal{E}_m(v)$, we have
		\begin{align*}
			&\left[\mathcal{E}_{m-1}^{[1]}(u), \mathcal{D}_{m}^{[2]}(v)\right]\mathcal{E}_{m}^{[2]}(v)+
			\mathcal{D}_{m}^{[2]}(v)\left[\mathcal{E}_{m-1}^{[1]}(u), \mathcal{E}_{m}^{[2]}(v)\right]\\
			=&\mathcal{D}_{m}^{[2]}(v)\left(\mathcal{E}_{m-1,m+1}^{[1]}(v)-\mathcal{E}_{m-1,m+1}^{[1]}(u)\right)\frac{P_{m+1,m}}{u-v}.
		\end{align*}
		By \eqref{eq:A:DaEb}, it can be further simplified as following:
		\begin{equation}
			\label{eq:Em-1EmEm-1m+1}
			\begin{aligned}	
				\left[\mathcal{E}_{m-1}^{[1]}(u), \mathcal{E}_{m}^{[2]}(v)\right]=&
				\left(\mathcal{E}_{m-1,m+1}^{[1]}(v)-\mathcal{E}_{m-1,m+1}^{[1]}(u)\right)\frac{P_{m+1,m}}{u-v}\\
				&-\left(\mathcal{E}_{m-1}^{[1]}(v)-\mathcal{E}_{m-1}^{[1]}(u)\right)\mathcal{E}_{m}^{[1]}(v)
				\frac{P_{m+1,m}}{u-v}.
			\end{aligned}
		\end{equation}
		Multiplying $u$ on both sides of \eqref{eq:Em-1EmEm-1m+1} and let $u\rightarrow\infty$, we obtain that
		$$\left[\mathcal{E}_{m-1}^{(1),[1]}, \mathcal{E}_{m}^{[2]}(v)\right]
		=\left(\mathcal{E}_{m-1,m+1}^{[1]}(v)-\mathcal{E}_{m-1}^{[1]}(v)\mathcal{E}_m^{[1]}(v)\right)P_{m+1,m},$$
		where $\mathcal{E}_{m-1}^{(1)}$ is the coefficient of $u^{-1}$ in $\mathcal{E}_{m-1}(u)$. A similar consideration on $v$ yields that
		$$\left[\mathcal{E}_{m-1}^{[1]}(u),\mathcal{E}_m^{(1),[2]}\right]=\mathcal{E}_{m-1,m+1}^{[1]}(u)P_{m+1,m}.$$
		Consequently, \eqref{eq:Em-1EmEm-1m+1} is reduced to
		\begin{align*}
			(u-v)\left[\mathcal{E}_{m-1}^{[1]}(u), \mathcal{E}_{m}^{[2]}(v)\right]=&
			\left[\mathcal{E}_{m-1}^{(1),[1]}, \mathcal{E}_{m}^{[2]}(v)\right]
			-\left[\mathcal{E}_{m-1}^{[1]}(u),\mathcal{E}_m^{(1),[2]}\right]
			+\mathcal{E}_{m-1}^{[1]}(u)\mathcal{E}_m^{[2]}(v)P_{m+1,m},
		\end{align*}
		which implies \eqref{eq:Em-1Em_odd}. The identity \eqref{eq:Fm-1Fm_odd} can be checked similarly.
	\end{proof}
	
	\begin{lem}
		\label{prop:odd:Serre}
		Let $M=2m+1\geqslant5$ be an odd integer and $\nu=(\nu_1,\nu_2,\ldots,\nu_M)$ be a symmetric composition of $N$. Then the following identities hold in $\X(\mathfrak{g}_{N})$:
		\begin{align}
			&\left[\left[\mathcal{E}_{m}(u),\mathcal{E}_{m-1}(v)\right],\mathcal{E}_{m-1}(w)\right]
			+\left[\left[\mathcal{E}_{m}(u),\mathcal{E}_{m-1}(w)\right],\mathcal{E}_{m-1}(v)\right]
			=0,\label{eq:odd:SerreE1E1E2}\\
			&\sum\limits_{\sigma\in\mathfrak{S}_3}
			\left[\left[\mathcal{E}_m^{[\sigma(1)]}(u_{\sigma(1)}),\left[\mathcal{E}_m^{[\sigma(2)]}(u_{\sigma(2)}),
			\left[\mathcal{E}_m^{[\sigma(3)]}(u_{\sigma(3)}),\mathcal{E}_{m-1}^{[4]}(v)\right]\right]\right]\right]=0,
			\label{eq:odd:SerreE2E2E2E1}\\
			&\left[\left[\mathcal{F}_{m}(u),\mathcal{F}_{m-1}(v)\right],\mathcal{F}_{m-1}(w)\right]
			+\left[\left[\mathcal{F}_{m}(u),\mathcal{F}_{m-1}(w)\right],\mathcal{F}_{m-1}(v)\right]
			=0,\label{eq:odd:SerreF1F1F2}\\
			&\sum\limits_{\sigma\in\mathfrak{S}_3}
			\left[\left[\mathcal{F}_m^{[\sigma(1)]}(u_{\sigma(1)}),\left[\mathcal{F}_m^{[\sigma(2)]}(u_{\sigma(2)}),
			\left[\mathcal{F}_m^{[\sigma(3)]}(u_{\sigma(3)}),\mathcal{F}_{m-1}^{[4]}(v)\right]\right]\right]\right]=0,
			\label{eq:odd:SerreF2F2F2F1}
		\end{align}
		where
		$\mathfrak{S}_3$ is the symmetric group on the set $\{1,2,3\}$.
	\end{lem}
	\begin{proof}
		On one hand, we have already shown in the proof of Lemma~\ref{lem:odd:Em1Em} that
		$$\left[\mathcal{E}_{m-1}^{[1]}(u), \mathcal{E}_m^{(1),[2]}\right]=\mathcal{E}_{m-1,m+1}^{[1]}(u)P_{m+1,m}.$$
		On the other hand, we deduce from \eqref{eq:blockRTT} that
		\begin{align*}
			\left[T_{m-1,m}^{[1]}(u), T_{m-1,m+1}^{[2]}(v)\right]=\frac{P_{m-1,m-1}}{u-v}\left(T_{m-1,m}^{[1]}(u)T_{m-1,m+1}^{[2]}(v)-T_{m-1,m}^{[1]}(v)T_{m-1,m+1}^{[2]}(u)\right),
		\end{align*}
		which yields that
		$$\left[\mathcal{E}_{m-1}^{(1),[1]},\mathcal{E}_{m-1,m+1}^{[2]}(v)\right]=0.$$
		Therefore,
		$$\left[\mathcal{E}_{m-1}^{(1),[1]},\left[\mathcal{E}_{m-1}^{[2]}(v),\mathcal{E}_m^{(1),[3]}\right]\right]
		=\left[\mathcal{E}_{m-1}^{(1),[1]},\mathcal{E}_{m-1,m+1}^{[2]}(v)P_{m+1,m}^{[2],[3]}\right]=0.$$
		Then \eqref{eq:odd:SerreE1E1E2} can be verified by Levendorskii's trick developed in \cite{l:gd}.
		\medskip
		
		By using \eqref{eq:blockRTT} again, we have
		\begin{align*}
			\left[T_{m,m+1}^{[1]}(u),T_{m-1,m+1}^{[2]}(v)\right]
			=&\frac{P_{m,m-1}}{u-v}\left(T_{m-1,m+1}^{[1]}(u)T_{m,m+1}^{[2]}(v)-T_{m-1,m+1}^{[1]}(v)T_{m,m+1}^{[2]}(u)\right)\\
			&+\frac{1}{u-v-\kappa}\sum_{r=1}^{2m+1}T_{m-1,r'}^{[2]}(v)T_{m,r}^{[1]}(u)Q_{r,m+1}.
		\end{align*}
		We multiply $v$ on both sides of the equality, let $v\rightarrow\infty$, and deduce that
		$$\left[T_{m,m+1}^{[1]}(u),T_{m-1,m+1}^{(1), [2]}\right]
		=T_{m,m+3}^{[1]}(u)Q_{m+3,m+1}.$$
		A similar consideration on $\left[T_{m,m+1}^{[1]}(u), T_{m,m+3}^{[2]}(v)\right]$ also shows that
		$$\left[T_{m,m+1}^{[1]}(u),T_{m,m+3}^{(1),[2]}\right]=0.$$
		Therefore, by block Gauss decomposition we have
		$$\left[\mathcal{E}_m^{(1),[1]},\mathcal{E}_{m-1,m+1}^{(1),[2]}\right]=\mathcal{E}_{m,m+3}^{(1),[1]}Q_{m+3,m+1},
		\text{ and }\left[\mathcal{E}_m^{(1),[1]},\mathcal{E}_{m,m+3}^{(1),[2]}\right]=0.$$
		Combining \eqref{eq:sp4:serrelevel1}, we obtain that
		$$\left[\mathcal{E}_m^{(1),[1]},\left[\mathcal{E}_m^{(1),[2]},
		\left[\mathcal{E}_m^{(1),[3]},\mathcal{E}_{m-1}^{(1),[4]}\right]\right]\right]=0.$$
		Then \eqref{eq:odd:SerreE2E2E2E1} can also be verified by Levendorskii's trick developed in \cite{l:gd}.  The relations \eqref{eq:odd:SerreF1F1F2} and \eqref{eq:odd:SerreF2F2F2F1} can be verified similarly.
	\end{proof}
	
	\begin{thm}
		\label{thm:odd:para}
		Suppose that $\nu=(\nu_1,\ldots,\nu_m,\nu_{m+1},\nu_{m+2},\ldots,\nu_{2m+1})$ is an odd symmetric composition of $N$. Then the extended Yangian $\X(\mathfrak{g}_{N})$ is the associative algebra presented by the generators that are the coefficients of the series $\mathcal{D}_{a,ij}(u)$, $\mathcal{E}_{b,kl}(u)$ and $\mathcal{F}_{b,lk}(u)$ for $a=1,\ldots,m+1$, $b=1,\ldots, m$, $1\leqslant i, j\leqslant \nu_a$, $1\leqslant k\leqslant\nu_b$ and $1\leqslant l\leqslant \nu_{b+1}$, we set
		$$\mathcal{D}_a(u)=\left(\mathcal{D}_{a,ij}(u)\right)_{\nu_a\times\nu_a},\,\, \mathcal{E}_b(u)=\left(\mathcal{E}_{b,ij}(u)\right)_{\nu_b\times\nu_{b+1}},
		\text{ and  }\mathcal{F}_b(u)=\left(\mathcal{F}_{b,ij}(u)\right)_{\nu_{b+1}\times\nu_b},$$
		the mere relations are given as following:
		\begin{enumerate}
			\item For $1\leqslant a,b\leqslant m+1$,
			\begin{equation}
				\left[\mathcal{D}_a^{[1]}(u), \mathcal{D}_b^{[2]}(v)\right]
				=\delta_{ab}\left(K_a(u,v)\mathcal{D}_a^{[1]}(u)\mathcal{D}_a^{[2]}(v)-\mathcal{D}_a^{[2]}(v)\mathcal{D}_a^{[1]}(u)K_a(u,v)\right),\label{eq:Dodd:DaDb}
			\end{equation}
			where $K_a(u,v)=\frac{P_{aa}}{u-v}$ if $1\leqslant a\leqslant m$, and $K_{m+1}(u,v)=\frac{P_{m+1,m+1}}{u-v}-\frac{Q_{m+1,m+1}}{u-v-\kappa_{m+1}}$.
			\item For $1\leqslant s\leqslant m+1$ and $1\leqslant a,b\leqslant m$,
			\begin{align}
				\left[\mathcal{D}_s^{[1]}(u),\mathcal{E}_a^{[2]}(v)\right]
				=&\left[\mathcal{D}_s^{[1]}(u),\mathcal{F}_a^{[2]}(v)\right]
				=0,\quad\text{if }s\neq a,a+1,
				\label{eq:Dodd:DsEa}\\
				\left[\mathcal{D}_a^{[1]}(u),\mathcal{E}_a^{[2]}(v)\right]
				=&\mathcal{D}_a^{[1]}(u)\frac{P_{aa}}{u-v}
				\left(\mathcal{E}_a^{[2]}(v)-\mathcal{E}_a^{[2]}(u)\right),
				\label{eq:Dodd:DaEa}\\
				\left[\mathcal{D}_{a+1}^{[1]}(u),\mathcal{E}_a^{[2]}(v)\right]
				=&\mathcal{D}_{a+1}^{[1]}(u)\left(\mathcal{E}_a^{[2]}(u)-\mathcal{E}_a^{[2]}(v)\right)\frac{P_{a+1,a+1}}{u-v}\nonumber\\
				&-\delta_{am}\mathcal{D}_{m+1}^{[1]}(u)\left(\mathcal{E}_m^{[2]}(u+\kappa_{m+1})-\mathcal{E}_m^{[2]}(v)\right)
				\frac{Q_{m+1,m+1}}{u-v+\kappa_{m+1}},
				\label{eq:Dodd:Da+1Ea}\\
				\left[\mathcal{D}_a^{[1]}(u),\mathcal{F}_a^{[2]}(v)\right]
				=&\left(\mathcal{F}_a^{[2]}(u)-\mathcal{F}_a^{[2]}(v)\right)\frac{P_{aa}}{u-v}\mathcal{D}_a^{[1]}(u),
				\label{eq:Dodd:DaFa}\\
				\left[\mathcal{D}_{a+1}^{[1]}(u),\mathcal{F}_a^{[2]}(v)\right]
				&=\frac{P_{a+1,a+1}}{u-v}\left(\mathcal{F}_a^{[2]}(v)-\mathcal{F}_a^{[2]}(u)\right)\mathcal{D}_{a+1}^{[1]}(u)\nonumber\\
				&-\frac{\delta_{am}Q_{m+1,m+1}}{u-v+\kappa_{m+1}}
				\left(\mathcal{F}_m^{[2]}(v)-\mathcal{F}_m^{[2]}(u+\kappa_{m+1})\right)\mathcal{D}_{m+1}^{[1]}(u),
				\label{eq:Dodd:Da+1Fa}\\
				\left[\mathcal{E}_a^{[1]}(u),\mathcal{F}_b^{[2]}(v)\right]
				=&\frac{\delta_{ab}}{u-v}
				\left(\widetilde{\mathcal{D}}_a^{[1]}(u)P_{a,a+1}\mathcal{D}_{a+1}^{[1]}(u)
				-\mathcal{D}_{a+1}^{[2]}(v)P_{a,a+1}\widetilde{\mathcal{D}}_a^{[2]}(v)\right).
				\label{eq:Dodd:EaFb}
			\end{align}
			\item For $1\leqslant a,b\leqslant m$,
			\begin{align}
				&\left[\mathcal{E}_a^{[1]}(u), \mathcal{E}_b^{[2]}(v)\right]
				=0,\qquad \text{ if }|a-b|>1,
				\label{eq:Dodd:EaEb}\\
				&\left[\mathcal{E}_a^{[1]}(u), \mathcal{E}_a^{[2]}(v)\right]
				=\frac{P_{aa}}{u-v}\left(\mathcal{E}_a^{[1]}(u)-\mathcal{E}_a^{[1]}(v)\right)
				\left(\mathcal{E}_a^{[2]}(u)-\mathcal{E}_a^{[2]}(v)\right),\quad\text{ if }a<m,
				\label{eq:Dodd:EaEa}\\
				&\left[\mathcal{E}_m^{[1]}(u), \mathcal{E}_m^{[2]}(v)\right]
				+\left[\mathcal{E}_m^{[2]}(u), \mathcal{E}_m^{[1]}(v)\right]\nonumber\\
				&\qquad=\left\{\mathcal{E}_m^{[1]}(u)-\mathcal{E}_m^{[1]}(v),
				\mathcal{E}_m^{[2]}(u)-\mathcal{E}_m^{[2]}(v)\right\}
				\frac{P_{m+1,m+1}-\frac{1}{2}Q_{m+1,m+1}}{u-v},
				\label{eq:Dodd:EmEm}\\
				&u\left[\mathring{\mathcal{E}}_{a-1}^{[1]}(u), \mathcal{E}_a^{[2]}(v)\right]
				-v\left[\mathcal{E}_{a-1}^{[1]}(u), \mathring{\mathcal{E}}_a^{[2]}(v)\right]
				=\mathcal{E}_{a-1}^{[1]}(u)P_{aa}\mathcal{E}_a^{[2]}(v),
				\label{eq:Dodd:Ea-1Ea}\\
				&\sum\limits_{\sigma\in\mathfrak{S}_p}\left[\mathcal{E}_a^{[\sigma(1)]}(u_{\sigma(1)}),\cdots,\left[\mathcal{E}_a^{[\sigma(p)]}(u_{\sigma(p)}), \mathcal{E}_b^{[2]}(v)\right]\right]
				=0, \text{ if }|a-b|=1,
				\label{eq:Dodd:ESerre}
			\end{align}
			where $p=3$ if $(a,b)=(m,m-1)$, and $p=2$ for all other pairs $(a,b)$ with $|a-b|=1$, and $\mathfrak{S}_p$ is the symmetric group in $\{1,2,\ldots, p\}$.
			\item For $1\leqslant a,b\leqslant m$,
			\begin{align}
				&\left[\mathcal{F}_a^{[1]}(u), \mathcal{F}_b^{[2]}(v)\right]=0,\qquad \text{ if }|a-b|>1,
				\label{eq:Dodd:FaFb}\\
				&\left[\mathcal{F}_a^{[1]}(u), \mathcal{F}_a^{[2]}(v)\right]
				=-\left(\mathcal{F}_a^{[2]}(u)-\mathcal{F}_a^{[2]}(v)\right)
				\left(\mathcal{F}_a^{[1]}(u)-\mathcal{F}_a^{[1]}(v)\right)\frac{P_{aa}}{u-v},\,\text{if }a<m,
				\label{eq:Dodd:FaFa}\\
				&\left[\mathcal{F}_m^{[1]}(u), \mathcal{F}_m^{[2]}(v)\right]
				+\left[\mathcal{F}_m^{[2]}(u), \mathcal{F}_m^{[1]}(v)\right]\nonumber\\
				&\qquad=\frac{P_{m+1,m+1}-\frac{1}{2}Q_{m+1,m+1}}{u-v}
				\left\{\mathcal{F}_m^{[1]}(u)-\mathcal{F}_m^{[1]}(v),
				\mathcal{F}_m^{[2]}(u)-\mathcal{F}_m^{[2]}(v)\right\},
				\label{eq:Dodd:FmFm}\\
				&u\left[\mathring{\mathcal{F}}_{a-1}^{[1]}(u), \mathcal{F}_a^{[2]}(v)\right]
				-v\left[\mathcal{F}_{a-1}^{[1]}(u), \mathring{\mathcal{F}}_a^{[2]}(v)\right]
				=-\mathcal{F}_{a}^{[2]}(v)P_{aa}\mathcal{F}_{a-1}^{[1]}(u),
				\label{eq:Dodd:Fm-1Fm}\\
				&\sum\limits_{\sigma\in\mathfrak{S}_p}\left[\mathcal{F}_a^{[\sigma(1)]}(u_{\sigma(1)}),\cdots,\left[\mathcal{F}_a^{[\sigma(p)]}(u_{\sigma(p)}), \mathcal{F}_b^{[2]}(v)\right]\right]
				=0, \text{ if }|a-b|=1,
				\label{eq:Dodd:FSerre}
			\end{align}
			where $p$ and $\mathfrak{S}_p$ has the same meaning as in (iii).
		\end{enumerate}
	\end{thm}
	\begin{proof}
		Let $\X^{\nu}(\mathfrak{g}_N)$ be the abstract associative algebra presented by the generators and relations as stated in the theorem. We also denoted by $\X^D(\mathfrak{g}_N)$ and $\X^R(\mathfrak{g}_N)$ the extended Yangian with the Drinfeld presentation and the R-matrix presentation, respectively. We shall show that $\X^{\nu}(\mathfrak{g}_N)$ is isomorphic to $\X^R(\mathfrak{g}_N)$.
		
		Let $T(u)$ be generator matrix of the extended Yangian $\X^R(\mathfrak{g}_N)$.
		The block Gauss decomposition \eqref{blockGaussDec} yields the matrices $\mathcal{D}_a(u)$, $\mathcal{E}_b(u)$ and $\mathcal{F}_b(u)$ for $a=1,\ldots, m+1$ and $b=1,\ldots, m$. It has been shown in Proposition~\ref{prop:generators} that the coefficients of their entries generate the algebra $\X^R(\mathfrak{g}_N)$. According to Proposition~\ref{prop:typeA} and Lemmas~\ref{lem:DaDm+1_odd}-\ref{prop:odd:Serre}, these generators satisfies all relations \eqref{eq:Dodd:DaDb}-\eqref{eq:Dodd:FSerre}. Hence, there is a canonical surjective homorphism of algebras
		$$\pi_{\nu}: \X^{\nu}(\mathfrak{g}_N)\rightarrow\X^R(\mathfrak{g}_N)$$
		such that
		$$\pi_{\nu}\left(\mathcal{D}_a(u)\right)=\mathcal{D}_a(u), \quad \pi_{\nu}\left(\mathcal{E}_b(u)\right)=\mathcal{E}_b(u),\text{ and }\pi_{\nu}\left(\mathcal{F}_b(u)\right)=\mathcal{F}_b(u),$$
		for $a=1,\ldots, m+1$ and $b=1,\ldots, m$.
		
		It remains to show that the homomorphism $\pi_{\nu}$ is injective. Since an explicit isomorphism $\pi^D:\X^D(\mathfrak{g}_N)\rightarrow\X^R(\mathfrak{g}_N)$ from the Drinfeld presented Yangian $\X^D(\mathfrak{g}_N)$ to the $R$-matrix presented Yangian $\X^R(\mathfrak{g}_N)$ has been given in \cite[Theorem~5.14]{jlm:ib}, we shall construct a surjective homomorphism $\widetilde{\pi}_{\nu}: \X^D(\mathfrak{g}_N)\rightarrow\X^{\nu}(\mathfrak{g}_N)$ such that $\pi^D=\pi_{\nu}\circ\widetilde{\pi}_{\nu}$. Then the injectivity of $\pi_{\nu}$ follows from the injectivity of the isomorphism $\pi^D$.
		
		In order to define $\widetilde{\pi}_{\nu}$, we consider the Gauss decomposition for each $\mathcal{D}_a(u)$ in $\X^{\nu}(\mathfrak{g}_N)$. Let
		$$\mathcal{D}_a(u)=\mathcal{D}_a^-(u)\mathsf{H}_a(u)\mathcal{D}_a^+(u),$$
		where $\mathsf{H}_a(u)=\diag\left[\mathsf{H}_{a;1}(u),\ldots,\mathsf{H}_{a;\nu_a}(u)\right]$ and $\mathcal{D}_a^-(u)=\left(\mathcal{D}_{a,ij}^-(u)\right)$ (resp. $\mathcal{D}_a^+(u)=\left(\mathcal{D}_{a,ij}^+(u)\right)$)
		is a lower-triangular (resp. an upper-triangular) matrix with diagonal entries $1$. Then we set
		\begin{align*}
			\mathsf{h}_i(u)=&\mathsf{H}_{a; r}(u),&&\text{ if }i=\nu_1+\cdots+\nu_{a-1}+r, 1\leqslant r\leqslant\nu_a,\\
			\mathsf{e}_j(u)=&\mathcal{D}_{a; r,r+1}^+(u),&&\text{ if }j=\nu_1+\cdots+\nu_{a-1}+r, 1\leqslant r\leqslant\nu_{a}-1,\\
			\mathsf{e}_j(u)=&\mathcal{E}_{a; \nu_a,1}(u),&&\text{ if }j=\nu_1+\cdots+\nu_a,\\
			\mathsf{f}_j(u)=&\mathcal{D}_{a; r+1,r}^-(u),&&\text{ if }j=\nu_1+\cdots+\nu_{a-1}+r, 1\leqslant r\leqslant\nu_{a}-1,\\
			\mathsf{f}_j(u)=&\mathcal{F}_{a; 1,\nu_a}(u),&&\text{ if }j=\nu_1+\cdots+\nu_a,
		\end{align*}
		for $i=1,\cdots,n+1$ and $j=1,\ldots, n$, where $n=\left[\frac{N}{2}\right]$. We have to show that $\mathsf{h}_i(u), i=1,2,\ldots, n+1, \mathsf{e}_j(u)$, and $\mathsf{f}_j(u)$, $j=1,2,\ldots, n$ satisfy all the relation in the Drinfeld presentation $\X^D(\mathfrak{g}_N)$ as stated in \cite[Theorem~5.14]{jlm:ib}.
		
		Fix $p=\nu_1+\nu_2+\cdots+\nu_m$. Since $\mathcal{D}_a(u), a=1,2,\ldots, m$, $\mathcal{E}_b(u)$ and $\mathcal{F}_b(u)$, $b=1,2,\ldots,m-1$ satisfy all the relations in the Parabolic presented Yangian $\Y^{\bar{\nu}}(\mathfrak{gl}_{p})$ associated with $\bar{\nu}=(\nu_1, \nu_2, \ldots, \nu_{m})$, which is isomorphic to the Drinfeld presented Yangian $\Y(\mathfrak{gl}_{p})$ as shown in \cite{bk:pp}. Hence,  $\mathsf{h}_i(u), \mathsf{e}_j(u), \mathsf{f}_j(u)$ for $i=1,2,\ldots,p, j=1,2,\ldots,p-1$ satisfies the relations in the Drinfeld presented Yangian $\Y(\mathfrak{gl}_{p})$. On the other hand, the relation \eqref{eq:Dodd:DaDb} implies that $\mathcal{D}_{m+1}(u)$ satisfies the RTT relation of $\X^R(\mathfrak{g}_{\nu_{m+1}})$. Hence, it follows from \cite{jlm:ib} that $\mathsf{h}_i(u), \mathsf{e}_j(u), \mathsf{f}_j(u)$ for $i=p+1,\ldots,n+1, j=p+1,\ldots,n$ satisfy all the relations in the Drinfeld presented Yangian $\X(\mathfrak{g}_{\nu_{m+1}})$. It remains to verify the relations involving $\mathsf{e}_p(u)$ or $\mathsf{f}_p(u)$.
		
		The relations~\eqref{eq:Dodd:DsEa}, \eqref{eq:Dodd:EaEb}, and \eqref{eq:Dodd:EaFb} imply that
		$$[\mathsf{h}_i(u), \mathsf{e}_p(u)]=[\mathsf{e}_j(u), \mathsf{e}_p(u)]=[\mathsf{f}_j(u), \mathsf{e}_p(u)]=0,$$
		for $i\leqslant \nu_1+\nu_2+\cdots+\nu_{m-1}$ and $j<\nu_1+\nu_2+\cdots+\nu_{m-1}$.
		Then the relation~\eqref{eq:Dodd:DaEa} with $a=m$ is equivalent to
		$$\left[\mathcal{D}_{m;ij}(u), \mathcal{E}_{m;kl}(v)\right]
		=\frac{\delta_{jk}}{u-v}\sum\limits_{r=1}^{\nu_m}
		\mathcal{D}_{m;ir}(u)\left(\mathcal{E}_{m;rl}(v)-\mathcal{E}_{m;rl}(u)\right).$$
		Note that $\mathsf{e}_p(u)=\mathcal{E}_{m;\nu_1,1}(u)$. According to the Gauss decomposition $\mathcal{D}_m(u)=\mathcal{D}^-_m(u)\mathcal{H}_m(u)\mathcal{D}^+_m(u)$ implies that $\mathsf{h}_i(u), \mathsf{f}_j(u), \mathsf{e}_j(u)$ for $i=p-\nu_m+1,\ldots, p-1, j=p-\nu_m+1,\ldots,p-2$ are generated by $\mathcal{D}_{m;ij}(u)$ with $j<\nu_m$. Hence,
		$$[\mathsf{h}_i(u),\mathsf{e}_p(v)]
		=[\mathsf{f}_i(u),\mathsf{e}_p(v)]=[\mathsf{e}_j(u),\mathsf{e}_p(v)]=0,$$
		for $p-\nu_m<i<p$ and $p-\nu_m<j<p-1$.
		
		Since $\mathcal{D}_m(u)$ is invertible, the relation~\eqref{eq:Dodd:DaEa} with $a=m$ can also be written as
		$$\left[\widetilde{\mathcal{D}}_m^{[1]}(u),\mathcal{E}_m^{[2]}(v)\right]
		=\frac{P_{mm}}{u-v}\left(\mathcal{E}_m^{[2]}(u)-\mathcal{E}_m^{[2]}(v)\right)
		\widetilde{\mathcal{D}}_m^{[1]}(u),$$
		where $\widetilde{\mathcal{D}}_m(u)$ is the inverse of $\mathcal{D}_m(u)$.
		Equivalently,
		$$\left[\widetilde{\mathcal{D}}_{m;ij}(u),\mathcal{E}_{m;kl}(v)\right]=
		\frac{1}{u-v}\left(\mathcal{E}_{m;il}(u)-\mathcal{E}_{m;il}(v)\right)
		\widetilde{\mathcal{D}}_{m;kj}(u).$$
		Note that $\widetilde{\mathcal{D}}_{m;\nu_m,\nu_m}(u)=\mathsf{h}_p(u)^{-1}$ and $\mathcal{E}_{m;\nu_m,1}(v)=\mathsf{e}_p(v)$, we deduce that
		$$[\mathsf{h}_p(u)^{-1},\mathsf{e}_p(v)]
		=\frac{1}{u-v}(\mathsf{e}_p(u)-\mathsf{e}_p(v))h_p(u)^{-1},$$
		that is equivalent to
		$$[\mathsf{h}_p(u),\mathsf{e}_p(v)]
		=\frac{1}{u-v}\mathsf{h}_p(u)(\mathsf{e}_p(v)-\mathsf{e}_p(u)).$$
		If $\nu_m>1$, $\widetilde{\mathcal{D}}_{m;\nu_m-1,\nu_m}(u)=-\mathsf{e}_{p-1}(u)\mathsf{h}_p(u)^{-1}$. Then we deduce that
		$$[\mathsf{e}_{p-1}(u),\mathsf{e}_p(v)]=\frac{1}{u-v}\left(\mathsf{e}_{p-1}(u)\mathsf{e}_p(v)-\mathsf{e}_{p-1}(u)\mathsf{e}_p(u)+\mathcal{E}_{m;\nu_{m-1},1}(v)-\mathcal{E}_{m;\nu_m-1,1}(u)\right),$$
		which yields that
		$$u[\mathring{\mathsf{e}}_{p-1}(u), e_p(v)]-v[\mathsf{e}_{p-1}(u),\mathring{\mathsf{e}}_p(v)]=\mathsf{e}_{p-1}(u)\mathsf{e}_p(v).$$
		When $\nu_m=1$, the same relation follows from \eqref{eq:Dodd:Ea-1Ea} directly.
		
		Next, we deduce from the relation~\eqref{eq:Dodd:EmEm} that
		\begin{align*}
			&\left[\mathcal{E}_{m;ij}(u),\mathcal{E}_{m;kl}(v)\right]
			+\left[\mathcal{E}_{m;kl}(u),\mathcal{E}_{m;ij}(v)\right]\\
			=&\frac{1}{u-v}\left\{
			\mathcal{E}_{m,il}(u)-\mathcal{E}_{m,il}(v), \mathcal{E}_{m,kj}(u)-\mathcal{E}_{m;kj}(v)
			\right\}\\
			&-\frac{\delta_{jl^{\prime}}}{2(u-v)}\sum_{r=1}^{\nu_{m+1}}\epsilon^{m+1}_r\epsilon^{m+1}_j
			\left\{
			\mathcal{E}_{m,ir}(u)-\mathcal{E}_{m,ir}(v), \mathcal{E}_{m,kr^{\prime}}(u)-\mathcal{E}_{m;kr^{\prime}}(v)
			\right\},
		\end{align*}
		which yields that
		\begin{align*}
			\left[\mathsf{e}_p(u),\mathsf{e}_p(v)\right]=&\frac{(\mathsf{e}_p(u)-\mathsf{e}_p(v))^2}{(u-v)}, &&\text{if }\nu_{m+1}>1,\\
			\left[\mathsf{e}_p(u),\mathsf{e}_p(v)\right]=&\frac{(\mathsf{e}_p(u)-\mathsf{e}_p(v))^2}{2(u-v)}, &&\text{if }\nu_{m+1}=1.
		\end{align*}
		All other relations stated in \cite[Theorem~5.14]{jlm:ib} can be verified using \eqref{eq:Dodd:DaDb}-\eqref{eq:Dodd:FSerre} similarly.
		
		In order to prove the surjectivity, we need to show all $\mathcal{D}_a(u)$, $\mathcal{E}_{b}(u)$, $\mathcal{F}_b(b)$ for $a=1,...,m+1$, $b=1,...,m$ can be
		generated by $h_i(u)$, $e_j(u)$, $f_j(u)$ for $i=1,...,n+1$, $j=1,...,n$ where
		$n=\left[\frac{N}{2}\right]$. Here we only show $\mathcal{D}_{m}(u)$, $\mathcal{D}_{m+1}(u)$, $\mathcal{E}_{m}(u)$ and $\mathcal{F}_m(u)$ lie in the
		image of $\widetilde{\pi}_{\nu}$, and other cases can be proved similarly.
		By relation \eqref{eq:Dodd:DaDb}, we have that the subalgebra generated by $\mathcal{D}_{m}(u)$ in $\X^{\nu}(\mathfrak{g}_N)$ is isomorphic to $Y(\mathfrak{gl}_{\nu_m})$ and the subalgebra generated by $\mathcal{D}_{m+1}(u)$ is isomorphic to $X(\mathfrak{g}_{\nu_{m+1}})$. The result in \cite{df:it} implies that
		$\mathcal{D}_{m}(u)$ can be generated by $h_i(u)$, $e_{j}(u)$ and $f_j(u)$ where
		$i=\sum\limits_{p=1}^{m-1}\nu_p+1,...,\sum\limits_{p=1}^{m-1}\nu_p+\nu_m$, $j=\sum\limits_{p=1}^{m-1}\nu_p+1,...,\sum\limits_{p=1}^{m-1}\nu_p+\nu_m-1$.
		Similarly, by the result in \cite{jlm:ib} we can show that $\mathcal{D}_{m+1}(u)$
		can be generated by  $h_i(u)$, $e_{j}(u)$ and $f_j(u)$ where
		$i=\sum\limits_{p=1}^{m}\nu_p+1,...,n+1$, $j=\sum\limits_{p=1}^{m}\nu_p+1,...,n$.
		
		Next we show $\mathcal{E}_m(u)$ can be generated by $\mathcal{D}_m(u)$, $\mathcal{D}_{m+1}(u)$ and $e_{\nu_{1}+...+\nu_m}=\mathcal{E}_{m;\nu_m,1}(u)$.
		Using relation \eqref{eq:Dodd:DaEa}, we have
		\begin{equation*}
			(u-v)\left[\mathcal{D}_{m;i\nu_m}(u), \mathcal{E}_{m;\nu_m 1}(v)\right]
			=\sum\limits_{r=1}^{\nu_m}
			\mathcal{D}_{m;ir}(u)\left(\mathcal{E}_{m;r1}(v)-\mathcal{E}_{m;r1}(u)\right),
		\end{equation*}
		which implies that
		\begin{equation*}
			\left[\mathcal{D}_{m;i\nu_m}^{(1)}, \mathcal{E}_{m;\nu_m 1}(v)\right]
			=	\mathcal{E}_{m;i1}(v)
		\end{equation*}
		for $i=1,..., \nu_m-1$. By relation \eqref{eq:Dodd:Da+1Ea}, we have
		\begin{equation*}
			(u-v)\left[\mathcal{D}_{m+1;1j}(u), \mathcal{E}_{m;k1}(v)\right]
			=\mathcal{D}_{m+1;11}(u)\left(\mathcal{E}_{m;kj}(u)-\mathcal{E}_{m;kj}(v)\right)
		\end{equation*}
		which implies that
		$	\left[\mathcal{D}_{m+1;1j}^{(1)}, \mathcal{E}_{m;k1}(v)\right]
		=	-\mathcal{E}_{m;kj}(v)
		$
		for $1<j<\nu_{m+1}$.
		Furthermore, by By relation \eqref{eq:Dodd:Da+1Ea}, we also have
		\begin{equation*}
			(u-v)\left[\mathcal{D}_{m+1;j\nu_{m+1}}(u), \mathcal{E}_{m;kj}(v)\right]
			=\mathcal{D}_{m+1;jj}(u)\left(\mathcal{E}_{m;k\nu_{m+1}}(u)-\mathcal{E}_{m;k\nu_{m+1}}(v)\right)
		\end{equation*}
		which implies $	\left[\mathcal{D}_{m+1;j\nu_{m+1}}^{(1)}, \mathcal{E}_{m;kj}(v)\right]
		=	-\mathcal{E}_{m;k\nu_{m+1}}(v)
		$
		for $1<j<\nu_{m+1}$. Thus, we show $\mathcal{E}_m(u)$ can be generated by $\mathcal{D}_m(u)$, $\mathcal{D}_{m+1}(u)$ and $e_{\nu_{1}+...+\nu_m}=\mathcal{E}_{m;\nu_m,1}(u)$.
		The case of $\mathcal{F}_{m}(u)$ can be proved similarly.
		
		Hence, there is a surjective homomorphism $\widetilde{\pi}_{\nu}: \X^D(\mathfrak{g}_N)\rightarrow\X^{\nu}(\mathfrak{g}_N)$ such that $\pi^D=\pi_{\nu}\circ\widetilde{\pi}_{\nu}$. This completes the proof.
	\end{proof}
	
	\begin{remark}
		When $a=b$, the relation \eqref{eq:Dodd:DaDb} is equivalent to
		$$\left(I-K_a(u,v)\right)\mathcal{D}_a^{[1]}(u)\mathcal{D}_a^{[2]}(v)
		=\mathcal{D}_a^{[2]}(v)\mathcal{D}_a^{[1]}(u)\left(I-K_a(u,v)\right).$$
		It shows that, for a fixed $a$, the subalgebra of $\X(\mathfrak{g}_N)$ generated by the coefficients of $\mathcal{D}_{a, ij}(u), i,j=1,2,\ldots,\nu_a$ is isomorphic to the Yangian $\Y(\mathfrak{gl}_{\nu_a})$ if $a\leqslant m$, and is isomorphic to $\X(\mathfrak{g}_{\nu_{m+1}})$ if $a=m+1$.
	\end{remark}
	
	\begin{remark}
		The parabolic presentation of $\X(\mathfrak{g}_N)$ associated to the composition $\nu=(N)$ is reduced to the R-matrix presentation of $\X(\mathfrak{g}_N)$. If $\mathfrak{g}_N=\mathfrak{o}_{2n+1}$ is of type $B$, the parabolic presentation of $\X(\mathfrak{o}_{2n+1})$ associated to the composition $\nu=(1,1,\ldots,1)$ is reduced to the Drinfeld presentation of $\X(\mathfrak{o}_{2n+1})$ given in \cite{jlm:ib}.
	\end{remark}
	
	\section{Parabolic presentation associated with an even composition}
	\label{sec:para:even}
	
	In this section, we always assume that $\nu=(\nu_1,\ldots,\nu_m,\nu_{m+1},\ldots,\nu_{2m})$ is an even symmetric composition of $N$. In this situation, $N=2n$ is also an even number and $\mathfrak{g}_N=\mathfrak{sp}_{2n}$ is the symplectic Lie algebra.
	
	\begin{lem}
		\label{lem:DaDEFm:even}
		Let $M=2m\geqslant 4$ be an even integer and $\nu=(\nu_1,\nu_2,\ldots,\nu_M)$ be a symmetric composition of $N$. Then the following relations hold in $\X(\mathfrak{g}_N)$:
		\begin{equation}
			\left[\mathcal{D}_a^{[1]}(u), \mathcal{D}_{m+1}^{[2]}(v)\right]=0,\quad
			\left[\mathcal{D}_a^{[1]}(u), \mathcal{E}_m^{[2]}(v)\right]=0,\quad\text{ and }
			\left[\mathcal{D}_a^{[1]}(u), \mathcal{F}_m^{[2]}(v)\right]=0.
			\label{eq:DaDEFm:even}
		\end{equation}
		for $1\leqslant a\leqslant m-1$.
	\end{lem}
	\begin{proof}
		It follows from \eqref{eq:TTpcomm} that
		$$\left[\mathcal{D}_a^{[1]}(u), \varPsi_{m-1}\left(T_{bc}^{[2]}(v)\right)\right]=0$$ for $b,c=m,m+1$ since $a<m$ and $\varPsi_{a-1}\left(T_{aa}(u)\right)=\mathcal{D}_a(u)$. Note that
		\begin{align*}
			\varPsi_{m-1}\left(T_{m,m}(v)\right)=&\mathcal{D}_{m}(v), &
			\varPsi_{m-1}\left(T_{m+1,m+1}(v)\right)=&\mathcal{D}_{m+1}(v)+\mathcal{F}_m(v)\mathcal{D}_m(v)\mathcal{E}_m(v),\\
			\varPsi_{m-1}\left(T_{m,m+1}(v)\right)=&\mathcal{D}_{m}(v)\mathcal{E}_m(v),&
			\varPsi_{m-1}\left(T_{m+1,m}(v)\right)=&\mathcal{F}_m(v)\mathcal{D}_{m}(v).
		\end{align*}
		So all relations in \eqref{eq:DaDEFm:even} hold.
	\end{proof}
	
	\begin{lem}
		\label{lem:EaFm:even}
		Let $M=2m\geqslant 6$ be an even integer and $\nu=(\nu_1,\nu_2,\ldots,\nu_M)$ be a symmetric composition of $N$. Then the following relations hold in $\X(\mathfrak{g}_N)$:
		\begin{align}
			\left[\mathcal{E}_a^{[1]}(u),\mathcal{D}_{m+1}^{[2]}(v)\right]=&0,&
			\left[\mathcal{E}_a^{[1]}(u),\mathcal{E}_m^{[2]}(v)\right]=&0,&
			\left[\mathcal{E}_a^{[1]}(u),\mathcal{F}_m^{[2]}(v)\right]=&0,\\
			\left[\mathcal{F}_a^{[1]}(u),\mathcal{D}_{m+1}^{[2]}(v)\right]=&0,&
			\left[\mathcal{F}_a^{[1]}(u),\mathcal{E}_m^{[2]}(v)\right]=&0,&
			\left[\mathcal{F}_a^{[1]}(u),\mathcal{F}_m^{[2]}(v)\right]=&0,
		\end{align}
		for $1\leqslant a\leqslant m-2$.
	\end{lem}
	\begin{proof}
		All these relations can be easily verified using \eqref{eq:TTpcomm}.
	\end{proof}
	
	\begin{lem}
		\label{lem:DmDm:even}
		Let $M=2m\geqslant 2$ be an even integer and $\nu=(\nu_1,\nu_2,\ldots,\nu_M)$ be a symmetric composition of $N$. Then the following relations hold in $\X(\mathfrak{g}_N)$:
		\begin{align}
			\left[\mathcal{D}_m^{[1]}(u),\mathcal{D}_m^{[2]}(v)\right]
			&=\frac{P_{mm}}{u-v}\left(\mathcal{D}_m^{[1]}(u)\mathcal{D}_m^{[2]}(v)-
			\mathcal{D}_m^{[1]}(v)\mathcal{D}_m^{[2]}(u)\right),\label{eq:DmDm_even}\\
			\left[\mathcal{D}_{m+1}^{[1]}(u),\mathcal{D}_{m+1}^{[2]}(v)\right]
			&=\frac{P_{mm}}{u-v}\left(\mathcal{D}_{m+1}^{[1]}(u)\mathcal{D}_{m+1}^{[2]}(v)-
			\mathcal{D}_{m+1}^{[1]}(v)\mathcal{D}_{m+1}^{[2]}(u)\right)\label{eq:Dm1Dm1_even}\\
			\left[\mathcal{D}_m^{[1]}(u),\mathcal{D}_{m+1}^{[2]}(v)\right]&
			=-\frac{1}{u-v-\kappa_m}\left(Q_{mm}\mathcal{D}_m^{[1]}(u)\mathcal{D}_{m+1}^{[2]}(v)-
			\mathcal{D}_{m+1}^{[2]}(v)\mathcal{D}_m^{[1]}(u)Q_{mm}\right)\label{eq:DmDm1_even}
		\end{align}
	\end{lem}
	
	\begin{proof}
		We consider the embedding homomorphism $\varPsi_{m-1}:\X(\mathfrak{g}_{2\nu_m})\rightarrow\X(\mathfrak{g}_N)$. It satisfies
		$$\psi_{m-1}\left(T_{mm}(u)\right)=\mathcal{D}_m(u),\quad \varPsi_{m-1}\left(\widetilde{T}_{m+1,m+1}(u)\right)=\widetilde{\mathcal{D}}_{m+1}(u),$$
		where $\widetilde{\mathcal{D}}_{m+1}(u)$ is the inverse of $\mathcal{D}_{m+1}(u)$.  Then the relation \eqref{eq:DmDm_even} follows from the block RTT relation \eqref{eq:blockRTT} in $\X(\mathfrak{g}_{2\nu_m})$. While \eqref{eq:Dm1Dm1_even} can be deduced from  the following relation in $\X(\mathfrak{g}_{2\nu_m})$:
		$$\widetilde{T}^{[1]}(u)\widetilde{T}^{[2]}(v)R^{12}(u-v)
		=R^{12}(u-v)\widetilde{T}^{[2]}(v)\widetilde{T}^{[1]}(u).$$
		
		For \eqref{eq:DmDm1_even}, we apply the homomorphism to the relation \eqref{blockTTtilde} in $\X(\mathfrak{g}_{2\nu_m})$ and obtain that
		$$\left[\mathcal{D}_m^{[1]}(u),\widetilde{\mathcal{D}}_{m+1}^{[2]}(v)\right]
		=\frac{1}{u-v-\kappa_m}\left(\widetilde{\mathcal{D}}_{m+1}^{[2]}(v)Q_{mm}\mathcal{D}_m^{[1]}(u)-
		\mathcal{D}_m^{[1]}(u)Q_{mm}\widetilde{\mathcal{D}}_{m+1}^{[2]}(v)\right).$$ This proves \eqref{eq:DmDm1_even} since $\widetilde{\mathcal{D}}_{m+1}(v)$ is the inverse of $\mathcal{D}_m(v)$.
	\end{proof}
	
	\begin{lem}
		\label{lem:DmEm:even}
		Let $M=2m\geqslant 2$ be an even integer and $\nu=(\nu_1,\nu_2,\ldots,\nu_M)$ be a symmetric composition of $N$. Then the following relations hold in $\X(\mathfrak{g}_N)$:
		\begin{align}
			\left[\mathcal{D}_m^{[1]}(u),\mathcal{E}_m^{[2]}(v)\right]
			=&\mathcal{D}_m^{[1]}(u)\left(\mathcal{E}_m^{[1]}(v)-\mathcal{E}_m^{[1]}(u)\right)\frac{P_{mm}+Q_{mm}}{u-v},
			\label{eq:DmEm_even}\\
			\left[\mathcal{D}_m^{[1]}(u),\mathcal{F}_m^{[2]}(v)\right]
			=&\frac{P_{mm}+Q_{mm}}{u-v}\left(
			\mathcal{F}_m^{[1]}(u)-\mathcal{F}_m^{[1]}(v)\right)\mathcal{D}_m^{[1]}(u).
			\label{eq:DmFm:even}
		\end{align}
	\end{lem}
	
	\begin{proof}
		By the block RTT realtion \eqref{eq:blockRTT} and $Q_{mm}=-Q_{m+1,m}$, the following relation holds in $\X(\mathfrak{g}_{2\nu_m})$ associated with the composition $(\nu_m,\nu_m)$:
		\begin{equation*}
			\begin{aligned}
				\left[T_{m,m}^{[1]}(u), T_{m,m+1}^{[2]}(v)\right]
				=&\left(T_{m,m}^{[2]}(u)T_{m,m+1}^{[1]}(v)-T_{m,m}^{[2]}(v)T_{m,m+1}^{[1]}(u)\right)\frac{P_{m,m}}{u-v}\\
				&+\left(T_{m,m+1}^{[2]}(v)T_{m,m}^{[1]}(u)-T_{m,m}^{[2]}(v)T_{m,m+1}^{[1]}(u)\right)\frac{Q_{m,m}}{u-v-\kappa_m}.
			\end{aligned}
		\end{equation*}
		Applying the embedding homomorphism $\varPsi_{m-1}$, we obtain the following relation in $\X(\mathfrak{g}_N)$:
		\begin{align*}
			\left[\mathcal{D}_m^{[1]}(u), \mathcal{D}_m^{[2]}(v)\mathcal{E}_m^{[2]}(v)\right]
			=
			&\left(\mathcal{D}_m^{[2]}(u)\mathcal{D}_m^{[1]}(v)\mathcal{E}_m^{[1]}(v)
			-\mathcal{D}_m^{[2]}(v)\mathcal{D}_m^{[1]}(u)\mathcal{E}_m^{[1]}(u)\right)\frac{P_{mm}}{u-v}\\
			&+\left(\mathcal{D}_m^{[2]}(v)\mathcal{E}_m^{[2]}(v)\mathcal{D}_m^{[1]}(u)
			-\mathcal{D}_m^{[2]}(v)\mathcal{D}_m^{[1]}(u)\mathcal{E}_m^{[1]}(u)\right)\frac{Q_{mm}}{u-v-\kappa_m}.
		\end{align*}
		Note that
		\begin{align*}
			\left[\mathcal{D}_m^{[1]}(u), \mathcal{D}_m^{[2]}(v)\mathcal{E}_m^{[2]}(v)\right]
			=\left[\mathcal{D}_m^{[1]}(u), \mathcal{D}_m^{[2]}(v)\right]\mathcal{E}_m^{[2]}(v)
			+\mathcal{D}_m^{[2]}(v)\left[\mathcal{D}_m^{[1]}(u), \mathcal{E}_m^{[2]}(v)\right].
		\end{align*}
		Then we deduce from \eqref{eq:DmDm_even} and Lemma~\ref{lem:PQ} that
		\begin{align*}
			\mathcal{D}_m^{[2]}(v)\left[\mathcal{D}_m^{[1]}(u), \mathcal{E}_m^{[2]}(v)\right]
			=&\mathcal{D}_m^{[2]}(v)\mathcal{D}_m^{[1]}(u)\left(\mathcal{E}_m^{[1]}(v)
			-\mathcal{E}_m^{[1]}(u)\right)\frac{P_{mm}}{u-v}\\
			&+\mathcal{D}_m^{[2]}(v)\left(\mathcal{E}_m^{[2]}(v)\mathcal{D}_m^{[1]}(u)
			-\mathcal{D}_m^{[1]}(u)\mathcal{E}_m^{[1]}(u)\right)\frac{Q_{mm}}{u-v-\kappa_m}.
		\end{align*}
		It yields that
		\begin{equation}
			\label{eq:D1E1proof}
			\begin{aligned}
				\left[\mathcal{D}_m^{[1]}(u), \mathcal{E}_m^{[2]}(v)\right]
				=&\mathcal{D}_m^{[1]}(u)\left(\mathcal{E}_m^{[1]}(v)
				-\mathcal{E}_m^{[1]}(u)\right)\frac{P_{mm}}{u-v}\\
				&+\left(\mathcal{E}_m^{[2]}(v)\mathcal{D}_m^{[1]}(u)
				-\mathcal{D}_m^{[1]}(u)\mathcal{E}_m^{[1]}(u)\right)\frac{Q_{mm}}{u-v-\kappa_m}.
			\end{aligned}
		\end{equation}
		Since $Q_{mm}^2=\nu_mQ_{mm}$, $P_{mm}Q_{mm}=Q_{mm}$ and $\kappa_{m}=\nu_m+1$,  by right multiplying $Q_{mm}$, we obtain
		\begin{align*}
			&\mathcal{E}_m^{[2]}(v)\mathcal{D}_m^{[1]}(u)Q_{mm}\frac{u-v-1}{u-v-\kappa_m}\\
			=&\mathcal{D}_m^{[1]}(u)\mathcal{E}_m^{[2]}(v)Q_{mm}-\mathcal{D}_m^{[1]}(u)\mathcal{E}_m^{[1]}(v)\frac{Q_{mm}}{u-v}
			+\mathcal{D}_m^{[1]}(u)\mathcal{E}_m^{[1]}(u)Q_{mm}\frac{\kappa_m(u-v-1)}{(u-v)(u-v-\kappa_m)}.
		\end{align*}
		Let $u=v+1$, we get $\mathcal{E}_m^{[2]}(v)Q_{mm}=\mathcal{E}_m^{[1]}(v)Q_{mm}$, and thus
		\begin{align*}
			\mathcal{E}_m^{[2]}(v)\mathcal{D}_m^{[1]}(u)\frac{Q_{mm}}{u-v-\kappa_m}
			=\mathcal{D}_m^{[1]}(u)\mathcal{E}_m^{[1]}(v)\frac{Q_{mm}}{u-v}
			+\mathcal{D}_m^{[1]}(u)\mathcal{E}_m^{[1]}(u)\frac{\kappa_m Q_{mm}}{(u-v)(u-v-\kappa_m)}.
		\end{align*}
		Then the relation \eqref{eq:D1E1proof} is simplified to  \eqref{eq:DmEm_even}.
		
		The relation \eqref{eq:DmFm:even} follows from a similar argument.
	\end{proof}
	
	\begin{lem}
		\label{lem:Dm1Em:even}
		Let $M=2m\geqslant 2$ be an even integer and $\nu=(\nu_1,\nu_2,\ldots,\nu_M)$ be a symmetric composition of $N$. Then the following relations hold in $\X(\mathfrak{g}_N)$:
		\begin{align}
			\left[\mathcal{D}_{m+1}^{[1]}(u),\mathcal{E}_m^{[2]}(v)\right]
			=&\mathcal{D}_{m+1}^{[1]}(u)\frac{P_{mm}+Q_{mm}}{u-v}\left(\mathcal{E}_m^{[1]}(u)-\mathcal{E}_m^{[1]}(v)\right),
			\label{eq:Dm1Em_even}\\
			\left[\mathcal{D}_{m+1}^{[1]}(u), \mathcal{F}_m^{[2]}(v)\right]
			=&\left(\mathcal{F}_m^{[1]}(v)-\mathcal{F}_m^{[1]}(u)\right)\frac{P_{mm}+Q_{mm}}{u-v}\mathcal{D}_{m+1}^{[1]}(u),
			\label{eq:Dm1Fm_even}
		\end{align}
	\end{lem}
	
	\begin{proof}
		By the block RTT relation \eqref{blockTTtilde} in $\X(\mathfrak{g}_{2\nu_m})$ associated with the composition $(\nu_m,\nu_m)$ and $Q_{mm}=-Q_{m,m+1}$, we have
		\begin{align*}
			&\left[T_{m,m+1}^{[1]}(u),\widetilde{T}_{m+1,m+1}^{[2]}(v)\right]\\
			=&\frac{1}{u-v}\left(T_{m,m}^{[1]}(u)P_{m,m}\widetilde{T}_{m,m+1}^{[2]}(v)+T_{m,m+1}^{[1]}(u)P_{m,m}\widetilde{T}_{m+1,m+1}^{[2]}(v)\right)\\
			&+\frac{1}{u-v-\kappa_m}\left(\widetilde{T}_{m+1,m+1}^{[2]}(v)Q_{m,m}T_{m,m+1}^{[1]}(u)
			+T_{m,m}^{[1]}(u)Q_{m,m}\widetilde{T}_{m,m+1}^{[2]}(v)\right).
		\end{align*}
		Now the embedding homomorphism
		$\varPsi_{m-1}:\X(\mathfrak{g}_{2\nu_m})\rightarrow\X(\mathfrak{g}_N)$ satisfies
		\begin{align*}
			\varPsi_{m-1}\left(T_{m,m+1}(u)\right)=&\mathcal{D}_m(u)\mathcal{E}_m(u),\\ \varPsi_{m-1}\left(\widetilde{T}_{m,m+1}(v)\right)=&-\mathcal{E}_m(v)\widetilde{\mathcal{D}}_{m+1}(v),& \varPsi_{m-1}\left(\widetilde{T}_{m+1,m+1}(v)\right)=&\widetilde{\mathcal{D}}_{m+1}(v),
		\end{align*}
		where $\widetilde{D}_{m+1}(v)$ is the inverse of $\mathcal{D}_{m+1}(v)$.
		So we obtain the following relation in $\X(\mathfrak{g}_N)$:
		\begin{align*}
			&\left[\mathcal{D}_m^{[1]}(u)\mathcal{E}_m^{[1]}(u),\widetilde{\mathcal{D}}_{m+1}^{[2]}(v)\right]\\
			=&\frac{1}{u-v}\left(-\mathcal{D}_m^{[1]}(u)P_{mm}\mathcal{E}_m^{[2]}(v)\widetilde{\mathcal{D}}_{m+1}^{[2]}(v)
			+\mathcal{D}_m^{[1]}(u)\mathcal{E}_m^{[1]}(u)P_{mm}\widetilde{\mathcal{D}}_{m+1}^{[2]}(v)\right)\\
			&+\frac{1}{u-v-\kappa_m}\left(\widetilde{\mathcal{D}}_{m+1}^{[2]}(v)Q_{mm}\mathcal{D}_m^{[1]}(u)
			\mathcal{E}_m^{[1]}(u)
			-\mathcal{D}_m^{[1]}(u)Q_{mm}\mathcal{E}_m^{[2]}(v)\widetilde{\mathcal{D}}_{m+1}^{[2]}(v)\right).
		\end{align*}
		Note that
		$$\left[\mathcal{D}_m^{[1]}(u)\mathcal{E}_m^{[1]}(u),\widetilde{\mathcal{D}}_{m+1}^{[2]}(v)\right]
		=\mathcal{D}_m^{[1]}(u)\left[\mathcal{E}_m^{[1]}(u),\widetilde{\mathcal{D}}_{m+1}^{[2]}(v)\right]
		+\left[\mathcal{D}_m^{[1]}(u),\widetilde{\mathcal{D}}_{m+1}^{[2]}(v)\right]\mathcal{E}_m^{[1]}(u).$$
		It follows from \eqref{eq:DmDm1_even} that
		\begin{align*}
			&\mathcal{D}_m^{[1]}(u)\left[\mathcal{E}_m^{[1]}(u),\widetilde{\mathcal{D}}_{m+1}^{[2]}(v)\right]\\
			=&\frac{1}{u-v}\left(-\mathcal{D}_m^{[1]}(u)P_{mm}\mathcal{E}_m^{[2]}(v)\widetilde{\mathcal{D}}_{m+1}^{[2]}(v)
			+\mathcal{D}_m^{[1]}(u)\mathcal{E}_m^{[1]}(u)P_{mm}\widetilde{\mathcal{D}}_{m+1}^{[2]}(v)\right)\\
			&+\frac{1}{u-v-\kappa_m}\left(
			\mathcal{D}_m^{[1]}(u)Q_{mm}\widetilde{\mathcal{D}}_{m+1}^{[2]}(v)\mathcal{E}_m^{[1]}(u)
			-\mathcal{D}_m^{[1]}(u)Q_{mm}\mathcal{E}_m^{[2]}(v)\widetilde{\mathcal{D}}_{m+1}^{[2]}(v)\right).
		\end{align*}
		Then the invertibility of $\mathcal{D}_m(u)$ and Lemma~\ref{lem:PQ} imply that
		\begin{equation}
			\label{eq:sp2:E1D2:2}
			\begin{aligned}
				\left[\mathcal{E}_m^{[1]}(u),\widetilde{\mathcal{D}}_{m+1}^{[2]}(v)\right]
				=&\frac{P_{mm}}{u-v}\left(\mathcal{E}_m^{[2]}(u)-\mathcal{E}_m^{[2]}(v)\right)\widetilde{\mathcal{D}}_{m+1}^{[2]}(v)\\
				&+\frac{Q_{mm}}{u-v-\kappa_m}\left(\widetilde{\mathcal{D}}_{m+1}^{[2]}(v)\mathcal{E}_m^{[1]}(u)
				-\mathcal{E}_m^{[2]}(v)\widetilde{\mathcal{D}}_{m+1}^{[2]}(v)\right).
			\end{aligned}
		\end{equation}
		By left multiplying $Q_{mm}$, we obtain that
		\begin{align*}
			\frac{u-v-1}{u-v-\kappa_m}Q_{mm}\widetilde{\mathcal{D}}_{m+1}^{[2]}(v)\mathcal{E}_m^{[1]}(u)
			=&Q_{mm}\mathcal{E}_m^{[1]}(u)\widetilde{\mathcal{D}}_{m+1}^{[2]}(v)-\frac{Q_{mm}}{u-v}\mathcal{E}_m^{[2]}(u)\widetilde{\mathcal{D}}_{m+1}^{[2]}(v)\\
			&+\frac{\kappa_m(u-v-1)Q_{mm}}{(u-v)(u-v-\kappa_m)}\mathcal{E}_m^{[2]}(v)\widetilde{\mathcal{D}}_{m+1}^{[2]}(v).
		\end{align*}
		Let $v=u-1$, we get $Q_{mm}\mathcal{E}_m^{[1]}(u)=Q_{mm}\mathcal{E}_m^{[2]}(u)$ and thus
		$$\frac{Q_{mm}}{u-v-\kappa_m}\widetilde{\mathcal{D}}_{m+1}^{[2]}(v)\mathcal{E}_m^{[1]}(u)
		=\frac{Q_{mm}}{u-v}\mathcal{E}_m^{[2]}(u)\widetilde{\mathcal{D}}_{m+1}^{[2]}(v)
		+\frac{\kappa_m Q_{mm}}{(u-v)(u-v-\kappa_m)}\mathcal{E}_m^{[2]}(v)\widetilde{\mathcal{D}}_{m+1}^{[2]}(v).$$
		Then the relation \eqref{eq:sp2:E1D2:2} is reduced to \eqref{eq:Dm1Em_even}.
	\end{proof}
	
	\begin{lem}
		\label{lem:EmFm:even}
		Let $M=2m\geqslant 2$ be an even integer and $\nu=(\nu_1,\nu_2,\ldots,\nu_M)$ be a symmetric composition of $N$. Then the following relations hold in $\X(\mathfrak{g}_N)$:
		\begin{equation}
			\label{eq:EmFm:even}
			\left[\mathcal{E}_m^{[1]}(u), \mathcal{F}_m^{[2]}(v)\right]=\widetilde{\mathcal{D}}_m^{[1]}(u)\frac{P_{mm}+Q_{mm}}{u-v}\mathcal{D}_{m+1}^{[1]}(u)
			-\mathcal{D}_{m+1}^{[2]}(v)\frac{P_{mm}+Q_{mm}}{u-v}\widetilde{\mathcal{D}}_m^{[2]}(v)
		\end{equation}
	\end{lem}
	
	\begin{proof}
		By the block RTT relation \eqref{eq:blockRTT} and $Q_{m,m+1}=Q_{m+1,m}=-Q_{mm}$, the following relation holds in $\X(\mathfrak{g}_{2\nu_m})$ associated with the composition $(\nu_m,\nu_m)$:
		\begin{align*}
			&\left(1+\frac{Q_{mm}}{u-v-\kappa_m}\right)T_{m,m+1}^{[1]}(u)T_{m+1,m}^{[2]}(v)-T_{m+1,m}^{[2]}(v)T_{m,m+1}^{[1]}(u)\left(1+\frac{Q_{mm}}{u-v-\kappa_m}\right)\\
			=&\left(\frac{P_{mm}}{u-v}+\frac{Q_{mm}}{u-v-\kappa_m}\right)T_{m+1,m+1}^{[1]}(u)T_{mm}^{[2]}(v)
			-T_{m+1,m+1}^{[2]}(v)T_{mm}^{[1]}(u)\left(\frac{P_{mm}}{u-v}+\frac{Q_{mm}}{u-v-\kappa_m}\right).
		\end{align*}
		Applying the embedding homomorphism $\varPsi_{m-1}:\X(\mathfrak{g}_{2\nu_m})\rightarrow\X(\mathfrak{g}_N)$, we deduce that
		\begin{align*}
			&\left(1+\frac{Q_{mm}}{u-v-\kappa_m}\right)\mathcal{D}_m^{[1]}(u)\mathcal{E}_m^{[1]}(u)\mathcal{F}_m^{[2]}(v)\mathcal{D}_m^{[2]}(v)
			-\mathcal{F}_m^{[2]}(v)\mathcal{D}_m^{[2]}(v)\mathcal{D}_m^{[1]}(u)\mathcal{E}_m^{[1]}(u)\left(1+\frac{Q_{mm}}{u-v-\kappa_m}\right)\\
			=&\left(\frac{P_{mm}}{u-v}+\frac{Q_{mm}}{u-v-\kappa_m}\right)\mathcal{D}_{m+1}^{[1]}(u)\mathcal{D}_m^{[2]}(v)
			-\mathcal{D}_{m+1}^{[2]}(v)\mathcal{D}_m^{[1]}(u)\left(\frac{P_{mm}}{u-v}+\frac{Q_{mm}}{u-v-\kappa_m}\right)\\
			&+\left(\frac{P_{mm}}{u-v}+\frac{Q_{mm}}{u-v-\kappa_m}\right)
			\mathcal{F}_m^{[1]}(u)\mathcal{D}_m^{[1]}(u)\mathcal{E}_m^{[1]}(u)\mathcal{D}_m^{[2]}(v)\\
			&-\mathcal{F}_m^{[2]}(v)\mathcal{D}_m^{[2]}(v)\mathcal{E}_m^{[2]}(v)\mathcal{D}_m^{[1]}(u)
			\left(\frac{P_{mm}}{u-v}+\frac{Q_{mm}}{u-v-\kappa_m}\right).
		\end{align*}
		By a similar argument as the proof of \eqref{eq:DmEm_even}, we have
		\begin{equation*}
			\begin{aligned}
				\left(\frac{P_{mm}}{u-v}+\frac{Q_{mm}}{u-v-\kappa_m}\right)\mathcal{F}_m^{[1]}(u)\mathcal{D}_m^{[1]}(u)
				&=-\mathcal{F}_{m}^{[2]}(v)\mathcal{D}_{m}^{[1]}(u)
				+\frac{P_{mm}}{u-v}\mathcal{F}_{m}^{[1]}(v)\mathcal{D}_m^{[1]}(u)\\
				&+\left(1+\frac{Q_{mm}}{u-v-\kappa_m}\right)\mathcal{D}_{m}^{[1]}(u)\mathcal{F}_m^{[2]}(v).
			\end{aligned}
		\end{equation*}
		Applying $\varPsi_{m-1}$ to the following block RTT relation in $X(\mathfrak{g}_{2\nu_m})$
		\begin{equation*}
			\begin{aligned}
				\left[T_{m,m+1}^{[1]}(u),T_{mm}^{[2]}(v)\right]&=\frac{1}{u-v}\left(P_{mm}T_{m,m+1}^{[1]}(u)T_{mm}^{[2]}(v)-T_{m,m+1}^{[2]}(v)T_{mm}^{[1]}(u)P_{mm}\right)\\
				&-\left(T_{m,m+1}^{[2]}(v)T_{mm}^{[1]}(u)-T_{m,m}^{[2]}(v)T_{m,m+1}^{[1]}(u)\right)\frac{Q_{mm}}{u-v-\kappa_m}
			\end{aligned}
		\end{equation*}
		we can obtain that
		\begin{equation*}
			\begin{aligned}
				&\mathcal{D}_m^{[2]}(v)\mathcal{E}_m^{[2]}(v)\mathcal{D}_m^{[1]}(u)
				\left(\frac{P_{mm}}{u-v}+\frac{Q_{mm}}{u-v-\kappa_m}\right)=\\
				&\mathcal{D}_m^{[2]}(v)\mathcal{D}_m^{[1]}(u)\mathcal{E}_m^{[1]}(u)
				\left(1+\frac{Q_{mm}}{u-v-\kappa_m}\right)+\left(\frac{P_{mm}}{u-v}-1\right)
				\mathcal{D}_m^{[1]}(u)\mathcal{E}_m^{[1]}(u)\mathcal{D}_m^{[2]}(v)
			\end{aligned}
		\end{equation*}
		
		Then we can get
		\begin{equation}
			\label{sp2:E1F1:4}
			\begin{aligned}
				\left(1+\frac{Q_{mm}}{u-v-\kappa_m}\right)\mathcal{D}_m^{[1]}(u)&\left[\mathcal{E}_m^{[1]}(u),\mathcal{F}_m^{[2]}(v)\right]\mathcal{D}_m^{[2]}(v)
				=\left(\frac{P_{mm}}{u-v}+\frac{Q_{mm}}{u-v-\kappa_m}\right)\mathcal{D}_{m+1}^{[1]}(u)\mathcal{D}_m^{[2]}(v)\\
				&-\mathcal{D}_{m+1}^{[2]}(v)\mathcal{D}_m^{[1]}(u)\left(\frac{P_{mm}}{u-v}+\frac{Q_{mm}}{u-v-\kappa_m}\right).
			\end{aligned}
		\end{equation}
		Note that $1+\frac{Q_{mm}}{u-v-\kappa_m}$ has an inverse $1-\frac{Q_{mm}}{u-v-\kappa_m}$ and
		$$\left(I-\frac{Q_{mm}}{u-v-1}\right)\left(\frac{P_{mm}}{u-v}+\frac{Q_{mm}}{u-v-\kappa_m}\right)=\frac{P_{mm}+Q_{mm}}{u-v},$$
		we further obtain that
		$$\mathcal{D}_m^{[1]}(u)\left[\mathcal{E}_m^{[1]}(u),\mathcal{F}_m^{[2]}(v)\right]\mathcal{D}_m^{[2]}(v)\\
		=\frac{P_{mm}+Q_{mm}}{u-v}\mathcal{D}_{m+1}^{[1]}(u)\mathcal{D}_m^{[2]}(v)-
		\mathcal{D}_m^{[1]}(u)\mathcal{D}_{m+1}^{[2]}(v)\frac{P_{mm}+Q_{mm}}{u-v},$$
		which implies \eqref{eq:EmFm:even} since $\mathcal{D}_m(u)$ is invertible.
	\end{proof}
	
	\begin{lem}
		\label{lem:EmEm:even}
		Let $M=2m\geqslant 2$ be an even integer and $\nu=(\nu_1,\nu_2,\ldots,\nu_M)$ be a symmetric composition of $N$. Then the following relations hold in $\X(\mathfrak{g}_N)$:
		\begin{align}
			\label{eq:EmEm_even}
			\left[\mathcal{E}_m^{[1]}(u), \mathcal{E}_m^{[2]}(v)\right]
			=&\left(\mathcal{E}_m^{[1]}(u)-\mathcal{E}_m^{[1]}(v)\right)
			\frac{P_{mm}+Q_{mm}}{u-v}\left(\mathcal{E}_m^{[1]}(u)-\mathcal{E}_m^{[1]}(v)\right),\\
			\label{eq:FmFm_even}
			\left[\mathcal{F}_m^{[1]}(u), \mathcal{F}_m^{[2]}(v)\right]
			=&\left(\mathcal{F}_m^{[1]}(u)-\mathcal{F}_m^{[1]}(v)\right)
			\frac{P_{mm}+Q_{mm}}{u-v}\left(\mathcal{F}_m^{[1]}(u)-\mathcal{F}_m^{[1]}(v)\right).
		\end{align}
	\end{lem}
	\begin{proof}
		It follows from \eqref{eq:blockRTT} that the following relation holds in $\X(\mathfrak{g}_{2\nu_m})$:
		\begin{align*}
			\left[T_{m,m+1}^{[1]}(u), T_{m,m+1}^{[2]}(v)\right]
			=\frac{1}{u-v}\left(P_{mm}T_{m,m+1}^{[1]}(u)T_{m,m+1}^{[2]}(v)-T_{m,m+1}^{[2]}(v)T_{m,m+1}^{[1]}(u)P_{mm}\right).
		\end{align*}
		Applying the embedding homomorphism $\psi_{m-1}:\X(\mathfrak{g}_{2\nu_m})\rightarrow\X(\mathfrak{g}_N)$, we obtain that
		$$\left(I-\frac{P_{mm}}{u-v}\right)\mathcal{D}_m^{[1]}(u)\mathcal{E}_m^{[1]}(u)\mathcal{D}_m^{[2]}(v)\mathcal{E}_m^{[2]}(v)
		=\mathcal{D}_m^{[2]}(v)\mathcal{E}_m^{[2]}(v)\mathcal{D}_m^{[1]}(u)\mathcal{E}_m^{[1]}(u)\left(I-\frac{P_{mm}}{u-v}\right).$$
		Using \eqref{eq:DmEm_even} and Lemma \ref{lem:PQ}, we deduce that
		\begin{align*}
			&\left(I-\frac{P_{mm}}{u-v}\right)\mathcal{D}_m^{[1]}(u)\mathcal{D}_m^{[2]}(v)\left(\mathcal{E}_m^{[1]}(u)
			+\left(\mathcal{E}_m^{[2]}(u)-\mathcal{E}_m^{[2]}(v)\right)\frac{P_{mm}+Q_{mm}}{u-v}\right)\mathcal{E}_m^{[2]}(v)\\
			=&\mathcal{D}_m^{[2]}(v)\mathcal{D}_m^{[1]}(u)\left(\mathcal{E}_m^{[2]}(v)-\left(\mathcal{E}_m^{[1]}(v)-\mathcal{E}_m^{[1]}(u)\right)\frac{P_{mm}+Q_{mm}}{u-v}\right)\mathcal{E}_m^{[1]}(u)\left(I-\frac{P_{mm}}{u-v}\right).
		\end{align*}
		Since $\mathcal{D}_m(u)$ is invertible, it follows from \eqref{eq:DmDm_even} that
		\begin{align*}
			&\left(I-\frac{P_{mm}}{u-v}\right)\left(\mathcal{E}_m^{[1]}(u)
			+\left(\mathcal{E}_m^{[2]}(u)-\mathcal{E}_m^{[2]}(v)\right)\frac{P_{mm}+Q_{mm}}{u-v}\right)\mathcal{E}_m^{[2]}(v)\\
			=&\left(\mathcal{E}_m^{[2]}(v)-\left(\mathcal{E}_m^{[1]}(v)-\mathcal{E}_m^{[1]}(u)\right)\frac{P_{mm}+Q_{mm}}{u-v}\right)\mathcal{E}_m^{[1]}(u)\left(I-\frac{P_{mm}}{u-v}\right).
		\end{align*}
		Note that $Q_{mm}\mathcal{E}_m^{[1]}(u)=Q_{mm}\mathcal{E}_m^{[2]}(u)$ and $\mathcal{E}_m^{[1]}(u)Q_{mm}=\mathcal{E}_m^{[2]}(u)Q_{mm}$ as shown in Lemmas~\ref{lem:DmEm:even} and~\ref{lem:Dm1Em:even}, respectively, we obtain that
		\begin{equation}
			\label{sp2:E1E1:2}
			\begin{aligned}
				&\left[\mathcal{E}_m^{[1]}(u),\mathcal{E}_m^{[2]}(v)\right]\left(I-\frac{P_{mm}}{u-v}\right)\\
				=&\left(\mathcal{E}_m^{[1]}(u)-\mathcal{E}_m^{[1]}(v)\right)\frac{P_{mm}}{u-v}\left(\mathcal{E}_m^{[1]}(u)-\mathcal{E}_m^{[1]}(v)\right)\\
				&+\left(\mathcal{E}_m^{[1]}(u)-\mathcal{E}_m^{[1]}(v)\right)\left(\frac{Q_{mm}}{u-v}-\frac{I+Q_{mm}}{(u-v)^2}\right)\left(\mathcal{E}_m^{[2]}(u)-\mathcal{E}_m^{[2]}(v)\right).
			\end{aligned}
		\end{equation}
		The \eqref{eq:EmEm_even} follows since $I-\frac{P_{mm}}{u-v}$ has the inverse $\left(I+\frac{P_{mm}}{u-v}\right)\frac{(u-v)^2}{(u-v)^2-1}$.
		
		The relation \eqref{eq:FmFm_even} follows from a similar argument.
	\end{proof}
	
	\begin{lem}
		\label{lem:Dm1Em1:even}
		Let $M=2m\geqslant 4$ be an even integer and $\nu=(\nu_1,\nu_2,\ldots,\nu_M)$ be a symmetric composition of $N$. Then the following relations hold in $\X(\mathfrak{g}_N)$:
		\begin{align}
			&\left[\mathcal{D}_{m+1}^{[1]}(u), \mathcal{E}_{m-1}^{[2]}(v)\right]
			=\mathcal{D}_{m+1}^{[1]}(u)\left(\mathcal{E}_{m-1}^{[2]}(v)-\mathcal{E}_{m-1}^{[2]}(u+\kappa_m)\right)
			\frac{Q_{mm}}{u-v+\kappa_{m}},\label{eq:Dm1Em1:even}\\
			&\left[\mathcal{D}_{m+1}^{[1]}(u), \mathcal{F}_{m-1}^{[2]}(v)\right]				=-\frac{Q_{mm}}{u-v+\kappa_m}\left(\mathcal{F}_{m-1}^{[2]}(v)-\mathcal{F}_{m-1}^{[2]}(u+\kappa_m)\right)\mathcal{D}_{m+1}^{[1]}(u).\label{eq:Dm1Fm1:even}
		\end{align}
	\end{lem}
	
	\begin{proof}
		Using the equivalent RTT relation
		$$\widetilde{T}^{[1]}(u)\widetilde{T}^{[2]}(v)R^{12}(u-v)
		=R^{12}(u-v)\widetilde{T}^{[2]}(v)\widetilde{T}^{[1]}(u)$$
		in $\X(\mathfrak{g}_{2(\nu_{m-1}+\nu_m)})$ and the embedding homomorphism $\varPsi_{m-2}:\X(\mathfrak{g}_{2(\nu_{m-1}+\nu_m)})\rightarrow\X(\mathfrak{g}_N)$, we obtain that
		\begin{equation*}
			\left[\mathcal{D}_{m+1}^{[1]}(u),\mathcal{E}_{m+1}^{[2]}(v)\right]
			=\frac{1}{u-v}\mathcal{D}_{m+1}^{[1]}(u)\left(\mathcal{E}_{m+1}^{[1]}(v)-\mathcal{E}_{m+1}^{[1]}(u)\right)P_{m+2,m+1}.
		\end{equation*}
		Applying the transpose $t$ to the second tensor factor, we obtain that
		\begin{equation*}
			\left[\mathcal{D}_{m+1}^{[1]}(u),(\mathcal{E}_{m+1}^{[2]}(v))^{t_2}\right]
			=\frac{1}{u-v}\mathcal{D}_{m+1}^{[1]}(u)\left(\mathcal{E}_{m+1}^{[1]}(v)-\mathcal{E}_{m+1}^{[1]}(u)\right)Q_{m+2,m+1}.
		\end{equation*}
		By Proposition \ref{prop:trans}, $\mathcal{E}_{m+1}^t(v)=-\mathcal{E}_{m-1}(v+\kappa_m)$, and it follows from Lemma~\ref{lem:PQ} that
		$$\mathcal{E}_{m+1}^{[1]}(u)Q_{m+2,m+1}=\mathcal{E}_{m+1}^{t,[2]}(u)Q_{m+1,m+1}=-\mathcal{E}_{m-1}^{[2]}(u+\kappa_m)Q_{mm}.$$
		Hence,
		\begin{equation*}
			\left[\mathcal{D}_{m+1}^{[1]}(u),\mathcal{E}_{m-1}^{[2]}(v+\kappa_m)\right]
			=\frac{1}{u-v}\mathcal{D}_{m+1}^{[1]}(u)\left(\mathcal{E}_{m-1}^{[2]}(v+\kappa_m)-\mathcal{E}_{m-1}^{[2]}(u+\kappa_m)\right)Q_{mm},
		\end{equation*}
		which implies \eqref{eq:Dm1Em1:even}. The relation \eqref{eq:Dm1Fm1:even} can be verified similarly.
	\end{proof}
	
	\begin{lem}
		\label{lem:Em1Em:even}
		Let $M=2m\geqslant 4$ be an even integer and $\nu=(\nu_1,\nu_2,\ldots,\nu_M)$ be a symmetric composition of $N$. Then the following relations hold in $\X(\mathfrak{g}_N)$:
		\begin{align}
			u\left[\mathring{\mathcal{E}}_{m-1}^{[1]}(u),\mathcal{E}_m^{[2]}(v)\right]
			-v\left[\mathcal{E}_{m-1}^{[1]}(u),\mathring{\mathcal{E}}_m^{[2]}(v)\right]
			=&\mathcal{E}_{m-1}^{[1]}(u)\mathcal{E}_m^{[1]}(v)\left(P_{mm}+Q_{mm}\right),
			\label{eq:Em1Em:even}\\
			u\left[\mathring{\mathcal{F}}_{m-1}^{[1]}(u),\mathcal{F}_m^{[2]}(v)\right]
			-v\left[\mathcal{F}_{m-1}^{[1]}(u),\mathring{\mathcal{F}}_m^{[2]}(v)\right]
			=&-(P_{mm}+Q_{mm})\mathcal{F}_m^{[1]}(v)\mathcal{F}_{m-1}^{[1]}(u),
			\label{eq:Fm1Fm:even}
		\end{align}
		where $\mathring{\mathcal{X}}_a(u)=\sum\limits_{r\geqslant2}\mathcal{X}_a^{(r)} u^{-r}$ if $\mathcal{X}_a=\mathcal{E}_{m-1},\mathcal{E}_m, \mathcal{F}_{m-1}$ and $\mathcal{F}_m$.
	\end{lem}
	\begin{proof}
		We only verify \eqref{eq:Em1Em:even}. The relation \eqref{eq:Fm1Fm:even} can be verified similarly.
		
		According to the embedding theorem~\ref{embedding}, we only need to consider the case where $M=4$. In this situation, $T_{12}(u)=\mathcal{D}_1(u)\mathcal{E}_1(u)$ and  $T_{23}(v)=\mathcal{D}_2(v)\mathcal{E}_2(v)+\mathcal{F}_1(v)\mathcal{D}_1(v)\mathcal{E}_{13}(v)$.
		
		By the block RTT relation \eqref{eq:blockRTT}, we have
		\begin{align*}
			\left[T_{12}^{[1]}(u), T_{23}^{[2]}(v)\right]
			=&\left(T_{22}^{[2]}(u)T_{13}^{[1]}(v)-T_{22}^{[2]}(v)T_{13}^{[1]}(u)\right)\frac{P_{32}}{u-v}\\
			&+\frac{1}{u-v-\kappa}\sum_{r=1}^4T_{2r'}^{[2]}(v)T_{1r}^{[1]}(u)Q_{r2}.
		\end{align*}
		Meanwhile, we compute that
		\begin{align*}
			&\left[T_{12}^{[1]}(u), \mathcal{D}_2^{[2]}(v)\mathcal{E}_2^{[2]}(v)\right]
			=\mathcal{D}_1^{[1]}(u)\mathcal{D}_2^{[2]}(v)\left[\mathcal{E}_1^{[1]}(u),\mathcal{E}_2^{[2]}(v)\right]\\
			&\qquad\qquad+\mathcal{D}_1^{[1]}(u)\left[\mathcal{E}_1^{[1]}(u),\mathcal{D}_2^{[2]}(v)\right]\mathcal{E}_2^{[2]}(v)
			+\left[\mathcal{D}_1^{[1]}(u),\mathcal{D}_2^{[2]}(v)\mathcal{E}_2^{[2]}(v)\right]\mathcal{E}_1^{[1]}(u),\\
			&\left[T_{12}^{[1]}(u),\mathcal{F}_1^{[2]}(v)\mathcal{D}_1^{[2]}(v)\mathcal{E}_{13}^{[2]}(v)\right]
			=\left[T_{12}^{[1]}(u), T_{21}^{[2]}(v)\right]\mathcal{E}_{13}^{[2]}(v)\\
			&\qquad\qquad-\mathcal{F}_1^{[2]}(v)\left[T_{12}^{[1]}(u),T_{11}^{[2]}(v)\right]\mathcal{E}_{13}^{[2]}(v)
			+\mathcal{F}_1^{[2]}(v)\left[T_{12}^{[1]}(u), T_{13}^{[2]}(v)\right].
		\end{align*}
		It follows from \eqref{eq:A:DaDb}, \eqref{eq:A:DaEb}, and \eqref{eq:blockRTT} that
		\begin{equation}
			\label{eq:sp4:E1E2a}
			\begin{aligned}
				&\mathcal{D}_1^{[1]}(u)\mathcal{D}_2^{[2]}(v)\left[\mathcal{E}_1^{[1]}(u),\mathcal{E}_2^{[2]}(v)\right]\\
				=&\mathcal{D}_1^{[1]}(u)\mathcal{D}_2^{[2]}(v)\left(\mathcal{E}_1^{[1]}(u)\mathcal{E}_2^{[1]}(v)-\mathcal{E}_1^{[1]}(v)\mathcal{E}_2^{[1]}(v)-\mathcal{E}_{13}^{[1]}(u)+\mathcal{E}_{13}^{[1]}(v)\right)\frac{P_{32}}{u-v}\\
				&+\mathcal{D}_1^{[1]}(u)\mathcal{D}_2^{[2]}(v)\mathcal{E}_2^{[2]}(v)\mathcal{E}_1^{[1]}(u)\frac{Q_{22}}{u-v-\kappa}
				+\mathcal{D}_1^{[1]}(u)\mathcal{D}_2^{[2]}(v)\mathcal{E}_{13}^{[1]}(u)\frac{Q_{32}}{u-v-\kappa}\\
				&+\mathcal{D}_2^{[2]}(v)\mathcal{E}_{24}^{[2]}(v)\mathcal{D}_1^{[1]}(u)\frac{Q_{12}}{u-v-\kappa}.
			\end{aligned}
		\end{equation}
		
		Now, we consider the block RTT relation \eqref{eq:blockRTT}:
		\begin{align*}
			\left[T_{11}^{[1]}(u),T_{24}^{[2]}(v)\right]
			=&\left(T_{21}^{[2]}(u)T_{14}^{[1]}(v)-T_{21}^{[2]}(v)T_{14}^{[1]}(u)\right)\frac{P_{41}}{u-v}\\
			&+\frac{1}{u-v-\kappa}\sum_{r=1}^4T_{2r'}^{[2]}(v)T_{1r}^{[1]}(u)Q_{r1},
		\end{align*}
		which yields that
		\begin{align*}
			&\left[\mathcal{D}_1^{[1]}(u),\mathcal{E}_{24}^{[2]}(v)\right]\\
			&=\frac{1}{u-v-\kappa}
			\left(\mathcal{E}_{24}^{[2]}(v)\mathcal{D}_1^{[1]}(u)Q_{11}
			+\mathcal{D}_1^{[1]}(u)\mathcal{E}_2^{[2]}(v)\mathcal{E}_1^{[1]}(u)Q_{21}+\mathcal{D}_1^{[1]}(u)\mathcal{E}_{13}^{[1]}(u)Q_{31}\right).
		\end{align*}
		By right-multiplying $Q_{12}$ on the above equation, we obtain that
		\begin{align*}
			\mathcal{E}_{24}^{[2]}(v)\mathcal{D}_1^{[1]}(u)Q_{12}
			=&\frac{u-v-\kappa}{u-v-\kappa_2}\mathcal{D}_1^{[1]}(u)\mathcal{E}_{24}^{[2]}(v)Q_{12}
			-\frac{\nu_1}{u-v-\kappa_2}\mathcal{D}_1^{[1]}(u)\mathcal{E}_2^{[2]}(v)\mathcal{E}_1^{[1]}(u)Q_{22}\\
			&-\frac{\nu_1}{u-v-\kappa_2}\mathcal{D}_1^{[1]}(u)\mathcal{E}_{13}^{[1]}(u)Q_{32}.
		\end{align*}
		Note that $\mathcal{D}_1(u)$ and $\mathcal{D}_2(v)$ are invertible, we deduce from \eqref{eq:sp4:E1E2a} that
		\begin{equation}
			\label{eq:sp4:E1E2b}
			\begin{aligned}
				\left[\mathcal{E}_1^{[1]}(u),\mathcal{E}_2^{[2]}(v)\right]
				=&\left(\left(\mathcal{E}_1^{[1]}(u)-\mathcal{E}_1^{[1]}(v)\right)\mathcal{E}_2^{[1]}(v)-\mathcal{E}_{13}^{[1]}(u)+\mathcal{E}_{13}^{[1]}(v)\right)\frac{P_{32}}{u-v}\\
				&+\frac{1}{u-v-\kappa_2}
				\left(\mathcal{E}_{24}^{[2]}(v)Q_{12}+\mathcal{E}_2^{[2]}(v)\mathcal{E}_1^{[1]}(u)Q_{22}
				+\mathcal{E}_{13}^{[1]}(u)Q_{32}
				\right).
			\end{aligned}
		\end{equation}
		We further deduce by right-multiplying $Q_{22}$ on the equation  \eqref{eq:sp4:E1E2b} that
		\begin{align*}
			&\frac{u-v-1}{u-v-\kappa_2}\mathcal{E}_2^{[2]}(v)\mathcal{E}_1^{[1]}(u)Q_{22}
			=\mathcal{E}_1^{[1]}(u)\mathcal{E}_2^{[2]}(v)Q_{22}
			-\left(\mathcal{E}_1^{[1]}(u)-\mathcal{E}_1^{[1]}(v)\right)\mathcal{E}_2^{[1]}(v)\frac{Q_{32}}{u-v}\\
			&\qquad\qquad+\left(\mathcal{E}_{13}^{[1]}(u)-\mathcal{E}_{13}^{[1]}(v)\right)\frac{Q_{32}}{u-v}
			-\frac{\nu_2}{u-v-\kappa_2}\left(\mathcal{E}_{13}^{[1]}(u)Q_{32}+\mathcal{E}_{24}^{[2]}(v)Q_{12}\right).
		\end{align*}
		Then the relation \eqref{eq:sp4:E1E2b} is simplified as
		\begin{align*}
			\left[\mathcal{E}_1^{[1]}(u),\mathcal{E}_2^{[2]}(v)\right]
			=&\left(\left(\mathcal{E}_1^{[1]}(u)-\mathcal{E}_1^{[1]}(v)\right)\mathcal{E}_2^{[1]}(v)-\mathcal{E}_{13}^{[1]}(u)+\mathcal{E}_{13}^{[1]}(v)\right)\frac{P_{22}+Q_{22}}{u-v}\\
			&+\frac{1}{u-v-1}
			\left(\mathcal{E}_{24}^{[2]}(v)Q_{12}+\mathcal{E}_1^{[1]}(v)\mathcal{E}_2^{[2]}(v)Q_{22}
			+\mathcal{E}_{13}^{[1]}(v)Q_{32}
			\right).
		\end{align*}
		By multiplying $u-v-1$ on both sides of the equality and set $u=v+1$, we observe that
		$$\mathcal{E}_{24}^{[2]}(v)Q_{12}+\mathcal{E}_1^{[1]}(v)\mathcal{E}_2^{[2]}(v)Q_{22}
		+\mathcal{E}_{13}^{[1]}(v)Q_{32}=0,$$
		and hence,
		\begin{equation}
			\label{eq:sp4:E1E2c}
			\left[\mathcal{E}_1^{[1]}(u),\mathcal{E}_2^{[2]}(v)\right]
			=\left(\left(\mathcal{E}_1^{[1]}(u)-\mathcal{E}_1^{[1]}(v)\right)\mathcal{E}_2^{[1]}(v)-\mathcal{E}_{13}^{[1]}(u)+\mathcal{E}_{13}^{[1]}(v)\right)\frac{P_{22}+Q_{22}}{u-v}.
		\end{equation}
		This implies \eqref{eq:Em1Em:even} via a similar argument as in Lemma~\ref{lem:odd:Em1Em}.
	\end{proof}
	
	\begin{lem}
		\label{lem:serre:even}
		Let $M=2m\geqslant 6$ be an even integer and $\nu=(\nu_1,\nu_2,\ldots,\nu_M)$ be a symmetric composition of $N$. Then the following relations hold in $\X(\mathfrak{g}_N)$:
		\begin{align}
			&\left[\left[\mathcal{E}_{m-1}(u),\mathcal{E}_m(v)\right],\mathcal{E}_m(w)\right]
			+\left[\left[\mathcal{E}_{m-1}(u),\mathcal{E}_m(v)\right],\mathcal{E}_m(w)\right]
			=0,\label{eq:sp4:SerreE1E2E2}\\
			&\sum\limits_{\sigma\in\mathfrak{S}_3}
			\left[\left[\mathcal{E}_{m-1}^{[\sigma(1)]}(u_{\sigma(1)}),\left[\mathcal{E}_{m-1}^{[\sigma(2)]}(u_{\sigma(2)}),
			\left[\mathcal{E}_{m-1}^{[\sigma(3)]}(u_{\sigma(3)}),\mathcal{E}_m^{[4]}(v)\right]\right]\right]\right]=0,
			\label{eq:sp4:SerreE1E1E1E2}
		\end{align}
		where $\mathring{\mathcal{E}}_a(u)=\sum\limits_{r\geqslant2}\mathcal{E}_a^{(r)} u^{-r}$, $a=m-1,m$ and $\mathfrak{S}_3$ is the symmetric group on the set $\{1,2,3\}$.
	\end{lem}
	\begin{proof}
		Multiplying $v$ on both sides of \eqref{eq:sp4:E1E2c} and taking $v\rightarrow\infty$, we obtain
		$$\left[\mathcal{E}_{m-1}^{[1]}(u), \mathcal{E}_m^{(1),[2]}\right]=\mathcal{E}_{m-1,m+1}^{[1]}(u)\left(P_{mm}+Q_{mm}\right),$$
		where $\mathcal{E}_m^{(1)}$ is the coefficient of $v^{-1}$ in $\mathcal{E}_m(v)$. It shows that
		\begin{equation}
			\left[\mathcal{E}_{m-1}^{(1),[1]},\mathcal{E}_m^{(1),[2]}\right]=\mathcal{E}_{m-1,m+1}^{(1),[1]}\left(P_{mm}+Q_{mm}\right).
			\label{eq:sp4:serrelevel1}
		\end{equation}
		
		On the other hand, we have the following relation in $\X(\mathfrak{g}_{2(\nu_{m-1}+\nu_m)})$:
		\begin{align*}
			\left[T_{m,m+1}^{[1]}(u), T_{m-1,m+1}^{[2]}(v)\right]=\frac{P_{m-1,m}}{u-v}\left(T_{m-1,m+1}^{[1]}(u)T_{m,m+1}^{[2]}(v)-T_{m-1,m+1}^{[1]}(v)T_{m,m+1}^{[2]}(u)\right),
		\end{align*}
		which yields $\left[T_{m,m+1}^{[1]}(u), T_{m-1,m+1}^{(1),[2]}\right]=0$. Then we deduce that
		$$\left[\mathcal{E}_m^{(1),[1]},\mathcal{E}_{m-1,m+1}^{(1),[2]}\right]=0.$$
		This shows
		$$\left[\left[\mathcal{E}_{m-1}^{(1),[1]},\mathcal{E}_m^{(1),[2]}\right],\mathcal{E}_m^{(1),[3]}\right]
		=\left[\mathcal{E}_{m-1,m+1}^{(1),[1]}\left(P_{mm}^{[1][2]}+Q_{mm}^{[1][2]}\right),\mathcal{E}_m^{(1),[3]}\right]=0.$$
		
		By using \eqref{eq:blockRTT} again, we have the following relation in $\X(\mathfrak{g}_{2(\nu_{m-1}+\nu_m)})$:
		\begin{align*}
			\left[T_{m-1,m}^{[1]}(u),T_{m-1,m+1}^{[2]}(v)\right]
			=&\frac{P_{m-1,m-1}}{u-v}\left(T_{m-1,m}^{[1]}(u)T_{m-1,m+1}^{[2]}(v)-T_{m-1,m}^{[1]}(v)T_{m-1,m+1}^{[2]}(u)\right)\\
			&+\frac{1}{u-v-\kappa_{m-1}}\sum_{r=m-1}^{m+2}T_{m-1,r^{\prime}}^{[2]}(v)T_{m-1,r}^{[1]}(u)Q_{rm}.
		\end{align*}
		We multiply $v$ on both sides of the equality, let $v\rightarrow\infty$, and deduce that
		$$\left[T_{m-1,m}^{[1]}(u),T_{m-1,m+1}^{(1), [2]}\right]
		=T_{m-1,m+2}^{[1]}(u)Q_{m+2,m}.$$
		A similar consideration on $\left[T_{m-1,m}^{[1]}(u), T_{m-1,m+2}^{[2]}(v)\right]$ also shows that
		$$\left[T_{m-1,m}^{[1]}(u),T_{m-1,m+2}^{(1),[2]}\right]=0.$$
		Hence,
		$$\left[\mathcal{E}_{m-1}^{(1),[1]},\mathcal{E}_{m-1,m+1}^{(1),[2]}\right]=\mathcal{E}_{m-1,m+2}^{(1),[1]}Q_{m+2,m},
		\text{ and }\left[\mathcal{E}_{m-1}^{(1),[1]},\mathcal{E}_{m-1,m+2}^{(1),[2]}\right]=0.$$
		Combining \eqref{eq:sp4:serrelevel1}, we obtain that
		$$\left[\mathcal{E}_{m-1}^{(1),[1]},\left[\mathcal{E}_{m-1}^{(1),[2]},
		\left[\mathcal{E}_{m-1}^{(1),[3]},\mathcal{E}_m^{(1),[4]}\right]\right]\right]=0.$$
		Then the Serre relations \eqref{eq:sp4:SerreE1E2E2} and \eqref{eq:sp4:SerreE1E1E1E2} can be verified by Levendorskii's trick developed in \cite{l:gd}.
	\end{proof}
	
	\begin{thm}
		\label{thm:even:para}
		Suppose that $\nu=(\nu_1,\ldots,\nu_m,\nu_{m+1},\ldots,\nu_{2m})$ is an even symmetric composition of $N$. Then the extended Yangian $X(\mathfrak{g}_{N})$ is the associative algebra presented by the generators that are the coefficients of the series $\mathcal{D}_{a,ij}(u)$, $\mathcal{E}_{b,kl}(u)$ and $\mathcal{F}_{b,lk}(u)$ for $a=1,\ldots,m+1$, $b=1,\ldots, m$, $1\leqslant i, j\leqslant \nu_a$, $1\leqslant k\leqslant\nu_b$ and $1\leqslant l\leqslant \nu_{b+1}$, we set
		$$\mathcal{D}_a(u)=\left(\mathcal{D}_{a,ij}(u)\right)_{\nu_a\times\nu_a},\,\, \mathcal{E}_b(u)=\left(\mathcal{E}_{b,ij}(u)\right)_{\nu_b\times\nu_{b+1}},
		\text{ and  }\mathcal{F}_b(u)=\left(\mathcal{F}_{b,ij}(u)\right)_{\nu_{b+1}\times\nu_b},$$
		the mere relations are given as following:
		\begin{enumerate}
			\item For $1\leqslant a\leqslant m+1$ and $1\leqslant b\leqslant m$ such that $a\neq b$ and $(a,b)\neq(m,m+1), (m+1,m)$,
			\begin{align}
				\left[\mathcal{D}_a^{[1]}(u), \mathcal{D}_b^{[2]}(v)\right]=&0, \\
				\left[\mathcal{D}_a^{[1]}(u), \mathcal{D}_a^{[2]}(v)\right]
				=&\frac{P_{aa}}{u-v}
				\left(\mathcal{D}_a^{[1]}(u)\mathcal{D}_a^{[2]}(v)-\mathcal{D}_a^{[1]}(v)\mathcal{D}_a^{[2]}(u)\right),\\
				\left[\mathcal{D}_m^{[1]}(u), \mathcal{D}_{m+1}^{[2]}(v)\right]
				=&-\frac{Q_{mm}\mathcal{D}_m^{[1]}(u)\mathcal{D}_{m+1}^{[2]}(v)
					-\mathcal{D}_{m+1}^{[2]}(v)\mathcal{D}_{m}^{[1]}(u)Q_{mm}}{u-v-\kappa_m}.
			\end{align}
			\item For $1\leqslant r\leqslant m+1$ and $1\leqslant a,b\leqslant m$ such that $r\neq a,a+1$ and $(r,a)\neq(m+1,m-1)$,
			\begin{align}
				\left[\mathcal{D}_r^{[1]}(u),\mathcal{E}_a^{[2]}(v)\right]
				=&\left[\mathcal{D}_r^{[1]}(u),\mathcal{F}_a^{[2]}(v)\right]=0,\\
				\left[\mathcal{D}_a^{[1]}(u),\mathcal{E}_a^{[2]}(v)\right]
				=&\mathcal{D}_a^{[1]}(u)\left(\mathcal{E}_a^{[1]}(v)-\mathcal{E}_a^{[1]}(u)\right)\frac{P_{a+1,a}+\delta_{am}Q_{mm}}{u-v}\\
				\left[\mathcal{D}_{a+1}^{[1]}(u),\mathcal{E}_a^{[2]}(v)\right]
				=&\mathcal{D}_{a+1}^{[1]}(u)\frac{P_{a+1,a}+\delta_{am}Q_{mm}}{u-v}
				\left(\mathcal{E}_a^{[1]}(u)-\mathcal{E}_a^{[1]}(v)\right),\\
				\left[\mathcal{D}_{m+1}^{[1]}(u),\mathcal{E}_{m-1}^{[2]}(v)\right]
				=&\mathcal{D}_{m+1}^{[1]}(u)\left(\mathcal{E}_{m-1}^{[2]}(v)-\mathcal{E}_{m-1}^{[2]}(u+\kappa_m)\right)
				\frac{Q_{mm}}{u-v+\kappa_m},\\
				\left[\mathcal{D}_a^{[1]}(u),\mathcal{F}_a^{[2]}(v)\right]
				=&\frac{P_{a,a+1}+\delta_{am}Q_{mm}}{u-v}\left(\mathcal{F}_a^{[1]}(u)-\mathcal{F}_a^{[1]}(v)\right)
				\mathcal{D}_a^{[1]}(u),\\
				\left[\mathcal{D}_{a+1}^{[1]}(u),\mathcal{F}_a^{[2]}(v)\right]
				=&\left(\mathcal{F}_a^{[1]}(u)-\mathcal{F}_a^{[1]}(v)\right)\frac{P_{a,a+1}+\delta_{am}Q_{mm}}{u-v}\mathcal{D}_{a+1}^{[1]}(u),\\
				\left[\mathcal{D}_{m+1}^{[1]}(u),\mathcal{F}_{m-1}^{[2]}(v)\right]
				=&\frac{Q_{mm}}{u-v+\kappa_m}
				\left(\mathcal{F}_{m-1}^{[2]}(u+\kappa_m)-\mathcal{F}_{m-1}^{[2]}(v)\right)\mathcal{D}_{m+1}^{[1]}(u).
			\end{align}
			\item For $1\leqslant a,b\leqslant m$ such that $a\neq b$,
			\begin{align}
				\left[\mathcal{E}_a^{[1]}(u),\mathcal{F}_b^{[2]}(v)\right]=&0,\\
				\left[\mathcal{E}_a^{[1]}(u),\mathcal{F}_a^{[2]}(v)\right]
				=&\left(
				\widetilde{\mathcal{D}}_a^{[1]}(u)\frac{P_{a,a+1}}{u-v}\mathcal{D}_{a+1}^{[1]}(u)
				-\mathcal{D}_{a+1}^{[2]}(v)\frac{P_{a,a+1}}{u-v}\widetilde{\mathcal{D}}_a^{[2]}(v)\right)\nonumber\\
				&+\delta_{am}\left(\widetilde{\mathcal{D}}_m^{[1]}(u)\frac{Q_{mm}}{u-v}\mathcal{D}_{m+1}^{[1]}(u)
				-\mathcal{D}_{m+1}^{[2]}(v)\frac{Q_{mm}}{u-v}\widetilde{\mathcal{D}}_m^{[2]}(v)\right).
			\end{align}
			\item For $1\leqslant a,b\leqslant m$,
			\begin{align}
				&\left[\mathcal{E}_a^{[1]}(u), \mathcal{E}_b^{[2]}(v)\right]=0,\quad\text{ if }|a-b|>2,\\
				&\left[\mathcal{E}_a^{[1]}(u), \mathcal{E}_a^{[2]}(v)\right]
				=\left(\mathcal{E}_a^{[1]}(u)-\mathcal{E}_a^{[1]}(v)\right)\frac{P_{a+1,a}+\delta_{am}Q_{mm}}{u-v}
				\left(\mathcal{E}_a^{[1]}(u)-\mathcal{E}_a^{[1]}(v)\right),\\
				&u\left[\mathring{\mathcal{E}}_{a-1}^{[1]}(u), \mathcal{E}_a^{[2]}(v)\right]
				-v\left[\mathcal{E}_{a-1}^{[1]}(u), \mathring{\mathcal{E}}_a^{[2]}(v)\right]
				=\mathcal{E}_{a-1}^{[1]}(u)\mathcal{E}_a^{[1]}(v)\left(P_{a+1,a}+\delta_{am}Q_{mm}\right),\\
				&\sum\limits_{\sigma\in\mathfrak{S}_p}
				\left[\mathcal{E}_a^{[\sigma(1)]}(u_{\sigma(1)}),\cdots,
				\left[\mathcal{E}_a^{[\sigma(p)]}(u_{\sigma(p)}), \mathcal{E}_b^{[2]}(v)\right]\right]
				=0, \text{ if }|a-b|=1,
			\end{align}
			where $p=3$ if $(a,b)=(m-1,m)$ and $p=2$ for all other pairs $(a,b)$ with $|a-b|=1$, $\mathfrak{S}_p$ is the symmetric group on $\{1,2,\ldots,p\}$.
			\item For $1\leqslant a,b\leqslant m$,
			\begin{align}
				&\left[\mathcal{F}_a^{[1]}(u), \mathcal{F}_b^{[2]}(v)\right]=0,\quad\text{ if }|a-b|>2,\\
				&\left[\mathcal{F}_a^{[1]}(u), \mathcal{F}_a^{[2]}(v)\right]
				=-\left(\mathcal{F}_a^{[1]}(u)-\mathcal{F}_a^{[1]}(v)\right)\frac{P_{a,a+1}+\delta_{am}Q_{mm}}{u-v}
				\left(\mathcal{F}_a^{[1]}(u)-\mathcal{F}_a^{[1]}(v)\right),\\
				&u\left[\mathring{\mathcal{F}}_{a-1}^{[1]}(u), \mathcal{F}_a^{[2]}(v)\right]
				-v\left[\mathcal{F}_{a-1}^{[1]}(u), \mathring{\mathcal{F}}_a^{[2]}(v)\right]
				=-\left(P_{a,a+1}+\delta_{am}Q_{mm}\right)\mathcal{F}_{a}^{[1]}(v)\mathcal{F}_{a-1}^{[1]}(u),\\
				&\sum\limits_{\sigma\in\mathfrak{S}_p}\left[\mathcal{F}_a^{[\sigma(1)]}(u_{\sigma(1)}),\cdots,\left[\mathcal{F}_a^{[\sigma(p)]}(u_{\sigma(p)}), \mathcal{F}_b^{[2]}(v)\right]\right]
				=0, \text{ if }|a-b|=1,
			\end{align}
			where $p$ and $\mathfrak{S}_p$ have the same meaning as in (iii).
		\end{enumerate}
	\end{thm}
	\begin{proof}
		Let $\widehat{\X}(\mathfrak{g}_N)$ be the abstract associative algebra with the presentation given in the theorem. In the extended Yangian $\X(\mathfrak{g}_N)$, the block Gauss decomposition \eqref{blockGaussDec} gives us the matrices $\mathcal{D}_a(u)$, $\mathcal{E}_b(u)$ and $\mathcal{F}_b(u)$ for $a=1,\ldots, m+1$ and $b=1,\ldots, m$. We have shown in Proposition~\ref{prop:generators} that the coefficients of their entries generate the algebra $\X(\mathfrak{g}_N)$.
		
		We have also verified in Proposition~\ref{prop:typeA} and Lemmas~\ref{lem:DaDEFm:even}-\ref{lem:serre:even} that all these generators satisfy the relations in (i)-(v). Hence, there is a canonical surjective homorphism of algebras
		$$\pi: \widehat{\X}(\mathfrak{g}_N)\rightarrow\X(\mathfrak{g}_N)$$
		such that
		$$\pi\left(\mathcal{D}_a(u)\right)=\mathcal{D}_a(u), \quad \pi\left(\mathcal{E}_b(u)\right)=\mathcal{E}_b(u),\text{ and }\pi\left(\mathcal{F}_b(u)\right)=\mathcal{F}_b(u),$$
		for $a=1,\ldots, m+1$ and $b=1,\ldots, m$.
		
		The injectivity of $\pi$ can be proved using the same strategy as in Theorem~\ref{thm:odd:para}. We omit the details here.
	\end{proof}
	
	\begin{remark}
		If $\nu=(1,\ldots,1,1,\ldots,1)$, the parabolic presentation of $\X(\mathfrak{sp}_{2n})$ associated to $\nu$ is reduced to  the Drinfeld presentation as stated in \cite{jlm:ib}.
	\end{remark}

	\section{Centers in terms of parabolic generators}
	\label{sec:center}
	
	It is known in \cite{amr:rp} that the generator matrix $T(u)$ in the RTT presented Yangian $\X(\mathfrak{g}_N)$ gives $T^t(u+\kappa)T(u)=z(u)I_N$, in which the coefficients of $z(u)$ generate the center of $\X(\mathfrak{g}_N)$. In this section, we formulate an alternative expression for $z(u)$ in terms of the parabolic generators of $\X(\mathfrak{g}_N)$.
	
	\begin{lem}
		\label{lem:ct:FDE}
		The following identity holds in $\X(\mathfrak{g}_N)$:
		\begin{equation}
			\begin{aligned}
				&\mathcal{E}_a^t(u+\kappa_{a+1})\mathcal{D}_a^t(u+\kappa_a)\mathcal{F}_a^t(u+\kappa_{a+1})
				-\left(\mathcal{F}_a(u+\kappa_a)\mathcal{D}_a(u+\kappa_a)\mathcal{E}_a(u+\kappa_a)\right)^t\\
				=&\mathcal{D}_{a+1}^t(u+\kappa_a)-\frac{\mathrm{tr}\left(\widetilde{\mathcal{D}}_a(u+\kappa_{a+1})\mathcal{D}_a(u+\kappa_a)\right)}{\nu_a}\mathcal{D}_{a+1}^{t}(u+\kappa_{a+1}),
			\end{aligned}
		\end{equation}
		for $1\leqslant a<\frac{M}{2}$.
	\end{lem}
	\begin{proof}
		Since the embedding theorem~\ref{embedding} holds, it suffices to consider the case where $a=1$.
		
		It is known from the proof of Proposition~\ref{prop:trans} that
		\begin{align*}
			\left(\mathcal{D}_1(u+\kappa_1)\mathcal{E}_1(u+\kappa_1)\right)^t=&\mathcal{E}_1^t(u+\kappa_2)\mathcal{D}_1^t(u+\kappa_1),\\
			\left(\mathcal{F}_1(u+\kappa_1)\mathcal{D}_1(u+\kappa_1)\right)^t=&\mathcal{D}_1^t(u+\kappa_1)\mathcal{F}_1^t(u+\kappa_2),
		\end{align*}
		where $\kappa_1=\kappa$ and $\kappa_2=\kappa-\nu_1$.
		Then, it follows from Lemma~\ref{lem:APQ} that
		\begin{align*}
			&\mathcal{E}_1^{t,[1]}(u+\kappa_2)\mathcal{D}_1^{t,[1]}(u+\kappa_1)\mathcal{F}_1^{t,[1]}(u+\kappa_2)Q_{22}\\
			=&\mathcal{E}_1^{t,[1]}(u+\kappa_2)\left(\mathcal{F}_1(u+\kappa_1)\mathcal{D}_1(u+\kappa_1)\right)^{t,[1]}Q_{22}\\
			=&\pm\mathcal{E}_1^{t,[1]}(u+\kappa_2)\mathcal{F}_1^{[2]}(u+\kappa_1)\mathcal{D}_1^{[2]}(u+\kappa_1)Q_{12},
		\end{align*}
		and
		\begin{align*}
			&\left(\mathcal{F}_1(u+\kappa_1)\mathcal{D}_1(u+\kappa_1)\mathcal{E}_1(u+\kappa_1)\right)^{t,[1]}Q_{22}\\
			=&\mathcal{F}_1^{[2]}(u+\kappa_1)\mathcal{D}_1^{[2]}(u+\kappa_1)\mathcal{E}_1^{[2]}(u+\kappa_1)Q_{22}\\
			=&\pm\mathcal{F}_1^{[2]}(u+\kappa_1)\left(\mathcal{D}_1(u+\kappa_1)\mathcal{E}_1(u+\kappa_1)\right)^{t,[1]}Q_{12}\\
			=&\pm\mathcal{F}_1^{[2]}(u+\kappa_1)\mathcal{E}_1^{t,[1]}(u+\kappa_2)\mathcal{D}_1^{t,[1]}(u+\kappa_1)Q_{12}\\
			=&\pm\mathcal{F}_1^{[2]}(u+\kappa_1)\mathcal{E}_1^{t,[1]}(u+\kappa_2)\mathcal{D}_1^{[2]}(u+\kappa_1)Q_{12},
		\end{align*}
		where the negative sign appear merely in the special case where $\mathfrak{g}_N$ is of type $C$ and $M=3$.
		Hence,
		\begin{equation}
			\label{eq:ct:edfQ}
			\begin{aligned}
				&\left(\mathcal{E}_1^t(u+\kappa_2)\mathcal{D}_1^t(u+\kappa_1)\mathcal{F}_1^t(u+\kappa_2)
				-\left(\mathcal{F}_1(u+\kappa_1)\mathcal{D}_1(u+\kappa_1)\mathcal{E}_1(u+\kappa_1)\right)^t\right)^{[1]}Q_{22}\\
				=&\pm\left[\mathcal{E}_1^{t,[1]}(u+\kappa_2),\mathcal{F}_1^{[2]}(u+\kappa_1)\right]\mathcal{D}_1^{[2]}(u+\kappa_1)Q_{12},
			\end{aligned}
		\end{equation}
		where the negative sign appear merely in the special case where $\mathfrak{g}_N$ is of type $C$ and $M=3$.
		Note also that the relation \eqref{eq:Dodd:EaFb} implies that
		$$\left[\mathcal{E}_1^{[1]}(u),\mathcal{F}_1^{[2]}(v)\right]
		=\frac{1}{u-v}\left(
		\widetilde{\mathcal{D}}_1^{[1]}(u)P_{12}\mathcal{D}_2^{[1]}(u)
		-\mathcal{D}_2^{[2]}(v)P_{12}\widetilde{\mathcal{D}}_1^{[2]}(v)\right).$$
		Since $\mathcal{D}_1(u)$ commutes with $\mathcal{D}_2(v)$ if $a=1<\frac{M}{2}$, we apply the transpose $t$ to the first tensor factor and obtain that
		$$\left[\mathcal{E}_1^{t,[1]}(u),\mathcal{F}_1^{[2]}(v)\right]
		=\frac{1}{u-v}\left(
		\mathcal{D}_2^{t,[1]}(u)P_{12}^{t_1}\widetilde{\mathcal{D}}_2^{t,[1]}(u)
		-\mathcal{D}_2^{[2]}(v)P_{12}^{t_1}\widetilde{\mathcal{D}}_1^{[2]}(v)\right).$$
		Now, we know from Lemma~\ref{lem:PQ} that $P_{12}^{t_1}=\pm Q_{21}$, where the negative sign appear merely in the special case where $\mathfrak{g}_N$ is of type $C$ and $M=3$. Hence,
		\begin{align*}
			\left[\mathcal{E}_1^{t,[1]}(u),\mathcal{F}_1^{[2]}(v)\right]
			=&\pm\frac{1}{u-v}\left(
			\mathcal{D}_2^{t,[1]}(u)Q_{21}\widetilde{\mathcal{D}}_1^{t,[1]}(u)
			-\mathcal{D}_2^{[2]}(v)Q_{21}\widetilde{\mathcal{D}}_1^{[2]}(v)\right)\\
			=&\pm\frac{1}{u-v}\left(
			\mathcal{D}_2^{[2]}(u)Q_{21}\widetilde{\mathcal{D}}_1^{[2]}(u)
			-\mathcal{D}_2^{[2]}(v)Q_{21}\widetilde{\mathcal{D}}_1^{[2]}(v)\right).
		\end{align*}
		
		Now, we conclude by \eqref{eq:ct:edfQ} and Lemmas~\ref{lem:PQ}, \ref{lem:APQ} that
		\begin{align*}
			&\left(\mathcal{E}_1^t(u+\kappa_2)\mathcal{D}_1^t(u+\kappa_1)\mathcal{F}_1^t(u+\kappa_2)
			-\left(\mathcal{F}_1(u+\kappa_1)\mathcal{D}_1(u+\kappa_1)\mathcal{E}_1(u+\kappa_1)\right)^t\right)^{[1]}Q_{22}\\
			=&\frac{1}{\nu_1}\mathcal{D}_2^{[2]}(u+\kappa_1)Q_{21}\widetilde{\mathcal{D}}_1^{[2]}(u+\kappa_1)
			\mathcal{D}_1^{[2]}(u+\kappa_1)Q_{12}\\
			&-\frac{1}{\nu_1}\mathcal{D}_2^{[2]}(u+\kappa_2)Q_{21}\widetilde{\mathcal{D}}_1^{[2]}(u+\kappa_2)\mathcal{D}_1^{[2]}(u+\kappa_1)Q_{12}\\
			=& \left(\mathcal{D}_2^{[2]}(u+\kappa_1)
			-\frac{\mathrm{tr}\left(\widetilde{\mathcal{D}}_1(u+\kappa_2)\mathcal{D}_1(u+\kappa_1)\right)}{\nu_1}
			\mathcal{D}_2^{[2]}(u+\kappa_2)\right)Q_{22}\\
			=& \left(\mathcal{D}_2^t(u+\kappa_1)
			-\frac{\mathrm{tr}\left(\widetilde{\mathcal{D}}_1(u+\kappa_2)\mathcal{D}_1(u+\kappa_1)\right)}{\nu_1}
			\mathcal{D}_2^t(u+\kappa_2)\right)^{[1]}Q_{22}.
		\end{align*}
		It yields that
		\begin{align*}
			&\mathcal{E}_1^t(u+\kappa_2)\mathcal{D}_1^t(u+\kappa_1)\mathcal{F}_1^t(u+\kappa_2)
			-\left(\mathcal{F}_1(u+\kappa_1)\mathcal{D}_1(u+\kappa_1)\mathcal{E}_1(u+\kappa_1)\right)^t\\
			=& \mathcal{D}_2^t(u+\kappa_1)
			-\frac{\mathrm{tr}\left(\widetilde{\mathcal{D}}_1(u+\kappa_2)\mathcal{D}_1(u+\kappa_1)\right)}{\nu_1}
			\mathcal{D}_2^t(u+\kappa_2)
		\end{align*}
		by Lemma \ref{lem:PQtrs}, which completes the proof.
	\end{proof}

	\begin{thm}
		Let $\mathcal{D}_a(u), a=1,2,\ldots, M$ be the parabolic generator matrices of $\X(\mathfrak{g}_N)$ associated to the symmetric composition $\nu=(\nu_1, \nu_2, \ldots,\nu_M)$ of $N$.
		\begin{enumerate}
			\item If $M=2m+1$ is odd, then
			\begin{equation}
				z(u)=\prod\limits_{a=1}^{m}\frac{\mathrm{tr}\left(\widetilde{\mathcal{D}}_a(u+\kappa_{a+1})\mathcal{D}_a(u+\kappa_a)\right)}{\nu_a}\cdot \frac{\mathrm{tr}\left(\mathcal{D}_{m+1}^t(u+\kappa_{m+1})\mathcal{D}_{m+1}(u)\right)}{\nu_{m+1}}.
				\label{eq:ct:Bodd}
			\end{equation}
			\item If $\mathfrak{g}_N$ is of type $C$ and $M=2m$ is even, then
			\begin{equation}
				\label{eq:ct:Ceven}
				z(u)=\prod\limits_{a=1}^{m-1}\frac{\mathrm{tr}\left(\widetilde{\mathcal{D}}_a(u+\kappa_{a+1})\mathcal{D}_a(u+\kappa_a)\right)}{\nu_a}
				\cdot \frac{\mathrm{tr}\left(\mathcal{D}_m^t(u+\kappa_m)\mathcal{D}_{m+1}(u)\right)}{\nu_m}.
			\end{equation}
		\end{enumerate}
	\end{thm}
	
	\begin{proof}
		Recall from \eqref{zcenter} the series $z(u)$ is given by $T^t(u+\kappa)T(u)=z(u)I$.
		Set $z^{[i]}(u)$ be the corresponding central series in $X(\mathfrak{g}_{N-2(\nu_1+...,\nu_{i-1})})$. Thus $z^{[1]}(u)=z(u)$.
		
		Since  $T^t(u+\kappa)=z^{[1]}(u)\widetilde{T}(u)$, we have
		$$T_{11}^t(u+\kappa)=z^{[1]}(u)\widetilde{T}_{MM}(u), \text{ and } T_{22}^t(u+\kappa)=z^{[1]}(u)\widetilde{T}_{M-1,M-1}(u).$$
		The block Gauss decomposition yields that
		\begin{align*}
			T_{11}(u)=&\mathcal{D}_1(u),
			&T_{22}(u)=&\mathcal{D}_2(u)+\mathcal{F}_1(u)\mathcal{D}_1(u)\mathcal{E}_1(u),\\
			\widetilde{T}_{MM}(u)=&\widetilde{\mathcal{D}}_M(u),
			&\widetilde{T}_{M-1,M-1}(u)=&\widetilde{\mathcal{D}}_{M-1}(u)
			+\mathcal{E}_{M-1,M}(u)\widetilde{\mathcal{D}}_M(u)\mathcal{F}_{M,M-1}(u).
		\end{align*}
		Then we deduce by Lemma~\ref{prop:trans} that
		\begin{align*}
			&\mathcal{D}_2^t(u+\kappa)+\left(\mathcal{F}_1(u+\kappa)\mathcal{D}_1(u+\kappa)\mathcal{E}_1(u+\kappa)\right)^t\\
			=&z^{[1]}(u)\widetilde{\mathcal{D}}_{M-1}(u)+z^{[1]}(u)\mathcal{E}_{M-1,M}(u)\widetilde{\mathcal{D}}_M(u)\mathcal{F}_{M,M-1}(u)\\
			=&z^{[1]}(u)\widetilde{\mathcal{D}}_{M-1}(u)+\mathcal{E}_1^t(u+\kappa_2)\mathcal{D}_1^t(u+\kappa)\mathcal{F}_1^t(u+\kappa_2).
		\end{align*}
		Now, Lemma~\ref{lem:ct:FDE} implies that
		$$\mathcal{D}_2^t(u+\kappa)-z^{[1]}(u)\widetilde{\mathcal{D}}_{M-1}(u)
		=\mathcal{D}_2^t(u+\kappa)-\frac{\mathrm{tr}\left(\widetilde{\mathcal{D}}_1(u+\kappa_2)\mathcal{D}_1(u+\kappa)\right)}{\nu_1}\mathcal{D}_2^t(u+\kappa_2),$$
		which implies that
		$$z^{[1]}(u)I_{\nu_2}=\frac{\mathrm{tr}\left(\widetilde{\mathcal{D}}_1(u+\kappa_2)\mathcal{D}_1(u+\kappa)\right)}{\nu_1}\mathcal{D}_2^t(u+\kappa_2)\mathcal{D}_{M-1}(u).$$
		
		Now, the embedding $\Psi_1:\X(\mathfrak{g}_{N-2\nu_1})\rightarrow\X(\mathfrak{g}_N)$ yields that
		$$\Psi_1(T_{22}(u))=\mathcal{D}_2(u),\text{ and }\Psi_1(\widetilde{T}_{M-1,M-1}(u))=\widetilde{\mathcal{D}}_{M-1}(u).$$
		The center series $z^{[2]}(u)$ in $\X(\mathfrak{g}_{N-2\nu_1})$ satisfies
		$$T_{22}^t(u+\kappa_2)=z^{[2]}(u)\widetilde{T}_{M-1,M-1}(u).$$
		Hence,
		$$\mathcal{D}_2^t(u+\kappa_2)\mathcal{D}_{M-1}(u)=\Psi_1\left(T_{22}^t(u+\kappa_2)\widetilde{T}_{M-1,M-1}(u)^{-1}\right)=\Psi_1(z^{[2]}(u))I_{\nu_2}.$$
		We further obtain that
		$$z^{[1]}(u)=\frac{\mathrm{tr}\left(\widetilde{\mathcal{D}}_1(u+\kappa_2)\mathcal{D}_1(u+\kappa)\right)}{\nu_1}\Psi_1(z^{[2]}(u)).$$
		Then by the above recurrence formula, we get
		\begin{equation}
			z(u)=\prod\limits_{a=1}^{m}\frac{\mathrm{tr}\left(\widetilde{\mathcal{D}}_a(u+\kappa_{a+1})\mathcal{D}_a(u+\kappa_a)\right)}{\nu_a}\cdot \varPsi_{m}(z^{[m]}(u)),
		\end{equation}
		for $M=2m+1$ is odd and
		\begin{equation}
			z(u)=\prod\limits_{a=1}^{m-1}\frac{\mathrm{tr}\left(\widetilde{\mathcal{D}}_a(u+\kappa_{a+1})\mathcal{D}_a(u+\kappa_a)\right)}{\nu_a}\cdot \varPsi_{m-1}(z^{[m-1]}(u)),
		\end{equation}
		for $M=2m$.
		
		By Theorem \ref{embedding} and \eqref{zcenter}, we can derive
		\begin{equation*}
			\varPsi_{m}(z^{[m]}(u))=\frac{\mathrm{tr}\left(\mathcal{D}_{m+1}^t(u+\kappa_{m+1})\mathcal{D}_{m+1}(u)\right)}{\nu_{m+1}},
		\end{equation*}
		for $M=2m+1$, and
		\begin{equation*}
			\varPsi_{m-1}(z^{[m-1]}(u))=\frac{\mathrm{tr}\left(\mathcal{D}_m^t(u+\kappa_m)\mathcal{D}_{m+1}(u)\right)}{\nu_m},
		\end{equation*}
		for $M=2m$. Then we complete the proof.
	\end{proof}
	
	\begin{remark}
		\begin{enumerate}
			\item For the composition $\nu=(1,\ldots,1)$ of $2n+1$, we set $\mathsf{h}_i(u)=\mathcal{D}_i(u)$ in $\X(\mathfrak{o}_{2n+1})$. Then $\kappa_i=n+\frac{1}{2}-i$ and the central series of $\X(\mathfrak{o}_{2n+1})$ can be written as
			$$z(u)=\prod\limits_{i=1}^{n}\mathsf{h}_i\left(u+n-i-\frac{1}{2}\right)^{-1}\prod_{i=1}^{n+1}\mathsf{h}_i\left(u+n-i+\frac{1}{2}\right)\mathsf{h}_{n+1}(u),$$
			which has been given in \cite[Theorem~5.8]{jlm:ib}.
			\item For the composition $\nu=(1,\ldots,1)$ of $2n$, we set $\mathsf{h}_i(u)=\mathcal{D}_i(u)$ in $\X(\mathfrak{sp}_{2n})$. Then $\kappa_i=n-i$ and the central series of $\X(\mathfrak{sp}_{2n})$ can be written as
			$$z(u)=\prod\limits_{i=1}^{n-1}\mathsf{h}_i(u+n-i+1)^{-1}\prod_{i=1}^n\mathsf{h}_i(u+n-i+2)\mathsf{h}_{n+1}(u),$$
			that is the formula given in \cite[Theorem~5.8]{jlm:ib}.
		\end{enumerate}
	\end{remark}
	
	\section*{Acknowledgements}
	
	Z. Chang thanks the support of the National Natural Science Foundation of China (No. 12071150) and Science and Technology Planning Project of Guangzhou (No. 202102021204). Jing acknowledges the partial support of
	Simons Foundation grant MP-TSM-00002518 and the National Natural Science Foundation of China (No. 12171303). Liu acknowledges the National
	Natural Science Foundation of China grant no. 12471026.

\noindent
Z.C.:\newline
School of Mathematics\\
South China University of Technology\\
Guangzhou, Guangdong 510640, China\\
mazhhchang@scut.edu.cn
\bigskip

\noindent
N.J.:\\
Department of Mathematics\\
North Carolina State University, Raleigh, NC 27695, USA\\
jing@ncsu.edu
\bigskip

\noindent
M.L.:\\
School of Artificial Intelligence\\
Jianghan University, Wuhan 430056, Hubei, China\\
ming.l1984@gmail.com
\bigskip

\noindent
H.M.:\\
College of Mathematical Science,\\
Harbin Engineering University\\
Harbin, Heilongjiang 150001, China\\
hmamath@hrbeu.edu.cn

\end{document}